\def\DATE{\relax}
\magnification=1100
\baselineskip=14truept
\voffset=.75in
\hoffset=1truein
\hsize=4.5truein
\newdimen\hsizeGlobal
\hsizeGlobal=4.5truein%
\vsize=7.75truein
\parindent=.166666in
\pretolerance=500 \tolerance=1000 \brokenpenalty=5000

\footline={\vbox{\hsize=\hsizeGlobal\hfill{\rm\the\pageno}\hfill\llap{\sevenrm\DATE}}\hss}

\def\note#1{%
  \hfuzz=50pt%
  \vadjust{%
    \setbox1=\vtop{%
      \hsize 3cm\parindent=0pt\eightpoints\baselineskip=9pt%
      \rightskip=4mm plus 4mm\raggedright#1\hss%
      }%
    \hbox{\kern-4cm\smash{\box1}\hss\par}%
    }%
  \hfuzz=0pt
  }
\def\note#1{\relax}

\def\anote#1#2#3{\smash{\kern#1in{\raise#2in\hbox{#3}}}%
  \nointerlineskip}     

\newcount\equanumber
\equanumber=0
\newcount\sectionnumber
\sectionnumber=0
\newcount\subsectionnumber
\subsectionnumber=0
\newcount\snumber  
\snumber=0

\def\section#1{%
  \subsectionnumber=0%
  \snumber=0%
  \equanumber=0%
  \advance\sectionnumber by 1%
  \noindent{\bf \the\sectionnumber .~#1.~}%
}%
\def\subsection#1{%
  \advance\subsectionnumber by 1%
  \snumber=0%
  \equanumber=0%
  \noindent{\bf \the\sectionnumber .\the\subsectionnumber .~#1.~}%
}%
\def\prevs{\the\sectionnumber .\the\subsectionnumber .\the\snumber }

\long\def\Corollary#1{%
  \global\advance\snumber by 1%
  \bigskip
  \noindent{\bf Corollary~\prevs .}%
  \quad{\it#1}%
}%
\long\def\Lemma#1{%
  \global\advance\snumber by 1%
  \bigskip
  \noindent{\bf Lemma~\prevs .}%
  \quad{\it#1}%
}%
\def\Proof{\noindent{\bf Proof.~}}
\long\def\Proposition#1{%
  \advance\snumber by 1%
  \bigskip
  \noindent{\bf Proposition~\prevs .}%
  \quad{\it#1}%
}%
\long\def\Remark#1{%
  \bigskip
  \noindent{\bf Remark.~}#1%
}%
\long\def\Theorem#1{%
  \advance\snumber by 1%
  \bigskip
  \noindent{\bf Theorem~\prevs .}%
  \quad{\it#1}%
}%
\long\def\Statement#1{%
  \advance\snumber by 1%
  \bigskip
  \noindent{\bf Statement~\prevs .}%
  \quad{\it#1}%
}%
\def\ifundefined#1{\expandafter\ifx\csname#1\endcsname\relax}
\def\labeldef#1{\global\expandafter\edef\csname#1\endcsname{\prevs}}
\def\labelref#1{\expandafter\csname#1\endcsname}
\def\label#1{\ifundefined{#1}\labeldef{#1}\note{$<$#1$>$}\else\labelref{#1}\fi}

\def\preveq{(\the\sectionnumber .\the\subsectionnumber .\the\equanumber)}
\def\neq{\global\advance\equanumber by 1\eqno{\preveq}}

\def\ifundefined#1{\expandafter\ifx\csname#1\endcsname\relax}

\def\equadef#1{\global\advance\equanumber by 1%
  \global\expandafter\edef\csname#1\endcsname{\preveq}%
  \preveq}

\def\equaref#1{\expandafter\csname#1\endcsname}

\def\equa#1{%
  \ifundefined{#1}%
    \equadef{#1}%
  \else\equaref{#1}\fi}

\font\eightrm=cmr8%
\font\sixrm=cmr6%

\font\eightsl=cmsl8%

\font\eightbf=cmb8%

\font\eighti=cmmi8%
\font\sixi=cmmi6%

\font\eightsy=cmsy8%
\font\sixsy=cmsy6%

\font\eightex=cmex8%
\font\sixex=cmex6%
\font\fiveex=cmex5%

\font\eightit=cmti8%

\font\eighttt=cmtt8%

\font\tenbb=msbm10%
\font\eightbb=msbm8%
\font\sevenbb=msbm7%
\font\sixbb=msbm6%
\font\fivebb=msbm5%
\newfam\bbfam  \textfont\bbfam=\tenbb  \scriptfont\bbfam=\sevenbb  \scriptscriptfont\bbfam=\fivebb%

\font\tenbbm=bbm10

\font\tencmssi=cmssi10%
\font\sevencmssi=cmssi7%
\font\fivecmssi=cmssi5%
\newfam\ssfam  \textfont\ssfam=\tencmssi  \scriptfont\ssfam=\sevencmssi  \scriptscriptfont\ssfam=\fivecmssi%
\def\ssi{\fam\ssfam\tencmssi}%

\font\tenfrak=cmfrak10%
\font\eightfrak=cmfrak8%
\font\sevenfrak=cmfrak7%
\font\sixfrak=cmfrak6%
\font\fivefrak=cmfrak5%
\newfam\frakfam  \textfont\frakfam=\tenfrak  \scriptfont\frakfam=\sevenfrak  \scriptscriptfont\frakfam=\fivefrak%
\def\frak{\fam\frakfam\tenfrak}%

\font\tenmsam=msam10%
\font\eightmsam=msam8%
\font\sevenmsam=msam7%
\font\sixmsam=msam6%
\font\fivemsam=msam5%

\def\bb{\fam\bbfam\tenbb}%

\def\hexdigit#1{\ifnum#1<10 \number#1\else%
  \ifnum#1=10 A\else\ifnum#1=11 B\else\ifnum#1=12 C\else%
  \ifnum#1=13 D\else\ifnum#1=14 E\else\ifnum#1=15 F\fi%
  \fi\fi\fi\fi\fi\fi}
\newfam\msamfam  \textfont\msamfam=\tenmsam  \scriptfont\msamfam=\sevenmsam  \scriptscriptfont\msamfam=\fivemsam%
\def\msam{\msamfam\tenmsam}%
\mathchardef\leq"3\hexdigit\msamfam 36%
\mathchardef\geq"3\hexdigit\msamfam 3E%

\font\tentt=cmtt11%
\font\seventt=cmtt9%
\textfont\ttfam=\tentt
\scriptfont7=\seventt%
\def\tt{\fam\ttfam\tentt}%

\def\eightpoints{%
\def\rm{\fam0\eightrm}%
\textfont0=\eightrm   \scriptfont0=\sixrm   \scriptscriptfont0=\fiverm%
\textfont1=\eighti    \scriptfont1=\sixi    \scriptscriptfont1=\fivei%
\textfont2=\eightsy   \scriptfont2=\sixsy   \scriptscriptfont2=\fivesy%
\textfont3=\eightex   \scriptfont3=\sixex   \scriptscriptfont3=\fiveex%
\textfont\itfam=\eightit  \def\it{\fam\itfam\eightit}%
\textfont\slfam=\eightsl  \def\sl{\fam\slfam\eightsl}%
\textfont\ttfam=\eighttt  \def\tt{\fam\ttfam\eighttt}%
\textfont\bffam=\eightbf  \def\bf{\fam\bffam\eightbf}%

\textfont\frakfam=\eightfrak  \scriptfont\frakfam=\sixfrak \scriptscriptfont\frakfam=\fivefrak  \def\frak{\fam\frakfam\eightfrak}%
\textfont\bbfam=\eightbb      \scriptfont\bbfam=\sixbb     \scriptscriptfont\bbfam=\fivebb      \def\bb{\fam\bbfam\eightbb}%
\textfont\msamfam=\eightmsam  \scriptfont\msamfam=\sixmsam \scriptscriptfont\msamfam=\fivemsam  \def\msam{\msamfam\eightmsam}

\rm%
}

\def\poorBold#1{\setbox1=\hbox{#1}\wd1=0pt\copy1\hskip.25pt\box1\hskip .25pt#1}

\mathchardef\lsim"3\hexdigit\msamfam 2E%
\mathchardef\gsim"3\hexdigit\msamfam 26%

\def\d{\,{\rm d}}
\def\D{{\rm D}}
\def\ds{\displaystyle}
\def\dsk{d_{\rm Sk}}
\long\def\DoNotPrint#1{\relax}
\def\E{{\rm E}}
\def\eqd{\buildrel {\rm d}\over =}
\def\finetune#1{#1}
\def\fixedref#1{#1\note{fixedref$\{$#1$\}$}}
\def\g#1{g_{[0,#1)}}
\def\Id{\hbox{\rm Id}}
\def\limeps{\lim_{\epsilon\to 0}}

\def\limk{\lim_{k\to\infty}}
\def\limn{\lim_{n\to\infty}}
\def\limsupn{\limsup_{n\to\infty}}

\def\limt{\lim_{t\to\infty}}
\def\ltrunc#1{\One_{(b_n,\infty)}(#1)}
\def\LTwo{{\rm L}^2[\,0,1\,]}
\def\okn{\overline k_n}
\def\oA{\overline A{}}
\def\oF{\overline F{}}

\def\oSS{\overline\SS}
\def\osn{\overline s_n}

\def\Prob{\hbox{\rm P}}
\def\qed{{\vrule height .9ex width .8ex depth -.1ex}}
\def\ss{\scriptstyle}
\def\trunc#1{\One_{[a_n,b_n]}(#1)}
\def\Utruncn{\One\Bigl\{\,{W_i\over W_{n+1}}< {m_n\over n}\,\Bigr\}}
\def\Var{{\rm Var}}

\def\boc{\note{{\bf BoC}\hskip-11pt\setbox1=\hbox{$\Bigg\downarrow$}%
         \dp1=0pt\ht1=0pt\ht1=0pt\leavevmode\raise -20pt\box1}}
\def\eoc{\note{{\bf EoC}\hskip-11pt\setbox1=\hbox{$\Bigg\uparrow$}%
         \dp1=0pt\ht1=0pt\ht1=0pt\leavevmode\raise 20pt\box1}}

\def\One{\hbox{\tenbbm 1}}

\def\ttL{{\tt L}}
\def\ttT{{\tt T}}

\def\calF{{\cal F}}

\def\LL{{\bb L}}
\def\MM{{\bb M}\kern .4pt}
\def\NN{{\bb N}\kern .5pt}
\def\RR{{\bb R}}
\def\SS{{\bb S}}
\def\TT{{\bb T}}

\def\M{{\ssi M}}


\pageno=1
\centerline{\bf INVARIANCE PRINCIPLES FOR SOME FARIMA}
\centerline{\bf AND NONSTATIONARY LINEAR PROCESSES}
\centerline{\bf IN THE DOMAIN OF A STABLE DISTRIBUTION}

\bigskip

\centerline{Ph.\ Barbe$^{(1)}$ and W.P.\ McCormick$^{(2)}$}
\centerline{${}^{(1)}$CNRS {\sevenrm(UMR {\eightrm 8088})}, ${}^{(2)}$University of Georgia}
 
{\narrower
\baselineskip=9pt\parindent=0pt\eightpoints

\bigskip

{\bf Abstract.} We prove some invariance principles for processes which 
generalize FARIMA processes, when the innovations are in the domain of
attraction of a nonGaussian stable distribution. The limiting processes are
extensions of the fractional L\'evy processes. The technique used is interesting
in itself; it extends an older idea of splitting a sample into a central part
and an extreme one, analyzing each part with different techniques, and then
combining the results. This technique seems to have the potential to be useful
in other problems in the domain of nonGaussian stable distributions.

\bigskip


\bigskip

\hfuzz=3pt
\noindent{\bf Keywords:} invariance principle, nonGaussian stable distribution,
FARIMA processes, fractional L\'evy stable process, fractional Brownian motion,
generalized integrated process.

\hfuzz=0pt
}

\bigskip

{\eightpoints

\obeylines\baselineskip=9pt

{\bf Contents}

\smallskip

{\sectionnumber=0
\def\section#1{\advance\sectionnumber by 1\subsectionnumber=1\line{\hskip\parindent\the\sectionnumber. #1\hfill}}
\def\subsection#1{\line{\hskip\parindent\quad\the\sectionnumber.\the\subsectionnumber. #1\hss}\advance\subsectionnumber by 1}

\section{Introduction and main result}
\section{Gaussian theory in nonGaussian domains}
\subsection{Convergence of finite dimensional distributions}%
\subsection{Tightness}%
\section{L\'evy stable limits and putting things together}%
\subsection{L\'evy stable limits}%
\subsection{Putting things together}%
\subsection{Proof of Theorem \fixedref{1.1}}%
\subsection{Proof of Corollary \fixedref{1.2}}%
\section{Point process limit and putting things together}%
\subsection{Point process limit}%
\subsection{Putting things together}%
\subsection{Proof of Theorem \fixedref{1.3}}%
\section{Application to $(g,F)$-processes}%
\subsection{Unbounded coefficients}%
\subsection{Bounded coefficients}%
\section{Application to FARIMA processes}%

References
\par
}
}

\bigskip


\def\prevs{\the\sectionnumber .\the\snumber }
\def\preveq{(\the\sectionnumber .\the\equanumber)}

\section{Introduction and main result}%
A different and as telling title for this paper could have been `Gaussian
theory in nonGaussian domains of attractions'.
Broadly speaking, we want to develop the old idea that a slight truncation
of $n$ independent random variables in the domain of attraction of a stable
distribution results in a limiting Gaussian behavior as $n$ tends to infinity. 
This idea allows one to consider
the bulk of the sample and the extreme values separately. Both parts of the
sample are amenable to different techniques and can be combined to study
functions of the whole sample. To our knowledge, this idea has been used only 
for studying sums of random variables in connection with the central limit 
theorem ---~see LePage, Woodroofe and Zinn (1981), Cs\"org\H o, 
Horv\'ath and Mason (1986), the latter having inspired the 
technique developed in the current
paper. The present paper suggests that the scope of this technique can 
be vastly and usefully extended. In a different direction, Chakrabarty and 
Samorodnitsky (2010) extends some of the earlier univariate results in 
studying the effect of a truncation on a multivariate mean when the summands
are in the domain of attraction of a multivariate stable distribution.

While this overall idea is general enough to be potentially applicable to
many situations, this paper focuses on its use to obtain
invariance principles for a class of processes which encompass
some nonstationary fractional autoregressive integrated moving
average (FARIMA) processes, when the innovations are in the domain of
attraction of a nonGaussian stable distribution. In particular, we
prove that in a certain range of the parameters, these processes
suitably rescaled converge to fractional stable L\'evy processes. The
usefulness of such result stems in part from the important
consequences of such invariance principle in statistics, the ubiquity
of FARIMA models in time series analysis, and the current interest in
modeling processes with innovations having an infinite variance. In a
somewhat different flavor, our own interest in the topic arises from
considerations in queueing and risk theory, where the results of this paper
are a first step toward proving a heavy traffic approximation in the 
context of FARIMA models. On a purely theoretical
basis, the results presented in this paper also shed light on an aspect
of the sharpness of the known convergence (Prokhorov, 1956)
of the rescaled random walk to a L\'evy stable process.

Since this paper is as much about the technique used as it is about its
application to FARIMA processes, it is noteworthy that a technical 
difficulty that one faces when dealing with many linear
or more complicated processes is that most standard tools relying on 
stationarity or a Markovian structure are not available. In the domain of
attraction of nonGaussian stable distributions, the point process approach,
as described in Resnick's (1987, 2007) books, has the 
potential to overcome some of the shortcomings of the more classical approaches.
However, in the context of FARIMA processes, our attempt to use the
point process methodology failed, because, as it is apparent in the proof
of the convergence of the partial sum process to a L\'evy stable one, that
methodology requires a form of maximal inequality to show that the bulk
of the innovations can be neglected ---~ see e.g.\ step 5 in the proof of 
Theorem 7.1 in Resnick (2007, \S 7.2.1). In contrast, the approach used in 
this paper does not require a direct use of such inequality, and seems to have 
a wider range of applicability, though, arguably, at the price of longer 
but perhaps conceptually simpler arguments.

\bigskip

\noindent{\bf Assumptions, notation, and processes of interests.}
Through\-out this paper we consider a distribution function $F$ on the real
line. We assume that
\setbox1=\vbox{\hsize=2.6in \par\noindent
$F$ is in the domain of attraction of a stable distribution of
   index $\alpha$ less than $2$.}
$$
  \box1
  \eqno{\equa{AssumptionStable}}
$$

Whenever $G$ is a cumulative distribution function, we write
$\overline G$ for $1-G$. Assumption \AssumptionStable\ implies
that both tails of $F$ are regularly varying of index $-\alpha$
and that $F$ is tail balanced, meaning the following. Write
$\M_{-1}F$ for the distribution function of $-X$. Then $|X|$
has tail distribution $\oF_*$ which coincides with $\oF
+\overline{\M_{-1}F}$ on the positive half-line when $F$ is
continuous. The tail balance condition is that $\oF\sim p\oF_*$
and $\overline{\M_{-1}F}\sim q\oF_*$ at infinity where $p$ and
$q$ are nonnegative numbers which add to $1$.

To the distribution function $F$ is associated the L\'evy measure $\nu$
whose density with respect to the Lebesgue measure $\lambda$ is
$$
  {\d\nu\over \d \lambda}(x)
  = p\alpha x^{-\alpha-1}\One_{(0,\infty)}(x)
  +q\alpha |x|^{-\alpha-1}\One_{(-\infty,0)}(x) \, .
$$
This L\'evy measure gives rise to a L\'evy process $L$; this process is 
determined by the requirements that
 
{\leftskip=\parindent \parindent=-\parindent\par
--~its increments are independent,

--~for any $0\leq s\leq t$, the increment $L(t)-L(s)$ has the same distribution
as $L(t-s)$, and

--~$L(t)$ has a stable distribution determined by
$$
  \E e^{i\theta L(t)}
  =\exp\Bigl( t\int 
  \bigl( e^{i\theta x}-1-i\theta x\One_{(-q^{1/\alpha},p^{1/\alpha})}(x)\bigr)
  \d\nu (x)\Bigr) \, .
  \eqno{\equa{cfL}}
$$
\par
}
Note that if either $p$ or $q$ vanishes, the integrals involved are still well
defined. 

Throughout this paper $\Id$ denotes the identity function on the real line.

Introducing the process
$$
  \tilde L
  =\cases{ L+{\ds\alpha\over\ds\alpha-1} (p-q-p^{1/\alpha}+q^{1/\alpha})\Id
              & if $\alpha\not= 1$,\cr
           \noalign{\vskip 3pt}
           L+(p\log p-q\log q)\Id
              & if $\alpha=1$,\cr}
$$
one can check that we have the more usual looking characteristic function
$$
  \E e^{i\theta\tilde L(t)}=\exp\Bigl( t\int \bigl( e^{i\theta x}-1-i\theta x
  \One_{(-1,1)}(x)\bigr) \d\nu(x)\Bigr) \, . 
$$
We refer to Bertoin (1996) for a general presentation of L\'evy processes.

When $\alpha$ exceeds $1$, the process $L$ has a drift whenever $p$ and $q$ are
different. Defining the process
$$
  L_0=\cases{ L-{\ds\alpha\over\ds\alpha-1} (p^{1/\alpha}-q^{1/\alpha}) \Id
                 & if $\alpha\not=1$, \cr\noalign{\vskip 3pt}
                 L-(p-q)\Id\log & if $\alpha=1$,\cr}
  \eqno{\equa{DefLZero}}
$$
one can check that, whenever it has a finite mean, $L_0$ is centered. Moreover,
since $L$, and therefore $L_0$, have independent increments, comparison of the
characteristic functions of $L_0(t)$ and $L_0(1)$ shows that for all values of
$\alpha$ in $(0,2)$, we have the scaling property
$$
  L_0(t\,\cdot\,) \eqd t^{1/\alpha}L_0(\,\cdot\,)
  \eqno{\equa{LZeroScaling}}
$$
for any nonnegative $t$.

Given a continuous function $k$ on $[\,0,1\,]$, we define
the stochastic process
$$
  k\star \dot L(t)=\int_0^t k(t-s)\d L(s) \, ;
$$
this stochastic integral is well defined because $|k|^\alpha$ is integrable
on $[\,0,1\,]$ ---~see e.g.\ Samorodnitsky and Taqqu, 1994, exercise 3.6.2.
Related and more general continuous time integral of L\'evy processes have been
studied by Marquardt (2006), Bender and Marquardt (2008), and Magdziarz (2008).

We write $(X_i)_{i\geq 1}$ for a sequence of independent random variables,
all having distribution function $F$.

Throughout this paper, $(k_n)_{n\geq 1}$ is a sequence of piecewise linear
functions supported on $[\,0,1\,]$; each function 
$$
  k_n \hbox{ is linearly interpolated between the points }
 (i/n)_{0\leq i\leq n} \, .
  \eqno{\equa{HypKnLin}}
$$ 
We write $\okn$ for the c\`adl\`ag
function which coincides with $k_n$ on $(i/n)_{0\leq i\leq n}$ and is constant
on the intervals $[\,i/n,(i+1)/n)$. It is convenient to view $k$ and $k_n$ 
as defined on the 
whole real line and vanishing outside $[\,0,1\,]$, even though this may
introduce a discontinuity at $0$ or $1$. To such sequence we
associate the stochastic process
$$
  \SS_n(t)
  =\sum_{1\leq i\leq nt} k_n\Bigl(t-{i\over n}\Bigr) X_i
  =\sum_{1\leq i\leq n} k_n\Bigl(t-{i\over n}\Bigr) X_i \, , 
  \eqno{\equa{procXCont}}
$$
defined for $t$ nonnegative, as well as the piecewise constant one,
$$
  \oSS_n(t)= \sum_{1\leq i\leq n} \okn\Bigl(t-{i\over n}\Bigr) X_i \, .
  \eqno{\equa{procXStep}}
$$

Even though $k_n$ is a linear interpolation of $\okn$, the process $\SS_n$ is 
not a linear interpolation of $\oSS_n$ in general; indeed, for any integer $p$
between $0$ and $n$,
$$
  \SS_n\Bigl({p\over n}-\Bigr)
  = \sum_{1\leq i\leq p-1} k_n\Bigl({p-i\over n}\Bigr) X_i
  = \oSS_n\Bigl({p\over n}\Bigr) -k_n(0) X_p \, .
$$
For $\SS_n$ to be continuous and indeed a linear interpolation of $\oSS_n$, we
must have
$$
  k_n(0)=0 \, .
  \eqno{\equa{HypKnZero}}
$$

We define the c\`agl\`ad quantile function
$$
  F^\leftarrow(u)=\inf\{\, x \,:\, F(x)\geq u\,\} \, .
$$
We then introduce the centerings
$$
  s_n(t)=\sum_{1\leq i\leq n} k_n\Bigl(t-{i\over n}\Bigr)
  \int_{F^\leftarrow(1/n)}^{F^\leftarrow(1-1/n)} x \d F(x) 
$$
and
$$
  \osn(t)=\sum_{1\leq i\leq n} \okn\Bigl(t-{i\over n}\Bigr)
  \int_{F^\leftarrow(1/n)}^{F^\leftarrow(1-1/n)} x \d F(x) \, .
$$
The theoretical aspect of this paper presents a technique to derive functional
limit theorems for the processes $\SS_n-s_n$ and $\oSS_n-\osn$ under 
various assumptions. The results are then specialized to FARIMA processes.

For the processes $\SS_n$ and $\oSS_n$ to have some chance to
converge, it is quite natural to assume that either the sequence
$(k_n)$ or $(\overline k_n)$ converges to a limit $k$. The following 
assumption seems adequate for some applications as well as for the
theoretical developments,
\setbox1=\vbox{\hsize=3.2in \par\noindent
$(k_n)$ converges uniformly to a continuous function $k$ 
        on $[\,0,1\,]$.}
$$
  \box1
  \eqno{\equa{HypKnCvUnif}}
$$
For the distributions of the processes to be tight, we assume that
$(k_n)$ is uniformly H\"olderian in the following sense. Consider the
discrete modulus of continuity
$$
  \overline\omega_{k_n,r}(\delta)=\sup_{0\leq p\leq q\leq p+n\delta\leq n} 
  {1\over n} \sum_{1\leq i\leq n} 
  \Bigl| k_n\Bigl({q-i\over n}\Bigr)-k_n\Bigl({p-i\over n}\Bigr)\Bigr|^r
  \, .
$$
We will assume that there exists some positive $c$ and $\epsilon$ as well
as an even
positive integer $r$ such that for any $n$ large enough and for any 
positive $\delta$ less than $1$,
$$
  \overline\omega_{k_n,r}(\delta)\leq c \delta^{1+\epsilon} \, .
  \eqno{\equa{HypKnModulus}}
$$
Note that $\overline\omega_{k_n,r}=\overline\omega_{\okn,r}$.

To prove the continuity of the limiting process, we need a continuous analogue
of $\overline\omega_{k_n,r}$, namely
$$
  \omega_{k,r}(\delta)=\sup_{0\leq s\leq t\leq s+\delta\leq 1}\int_0^1 
  \bigl|k(t-u)-k(s-u)\bigr|^r \d u \, .
$$
We assume that for some positive real number $r$ and $\epsilon$ and 
any $\delta$ positive less than $1$,
$$
  \omega_{k,r}(\delta)\leq c\delta^{1+\epsilon} \, .
  \eqno{\equa{HypKModulus}}
$$

\bigskip

\noindent{\bf Main results.}
In order to state the main theoretical results of this paper, we write
${\rm C}[\,0,1\,]$ for the space of all continuous functions on $[\,0,1\,]$
and ${\rm D}[\,0,1\,]$ for the space of all c\`adl\`ag functions 
on $[\,0,1\,]$ --- see e.g.\ Billingsley (1968).

Our first result is an invariance principle for $\SS_n$ when the sequence
$(k_n)$ converges in ${\rm C}[\,0,1\,]$.

\Theorem{\label{LimitCD}%
  Assume that \AssumptionStable, \HypKnLin\ and \HypKnZero--\HypKModulus\ 
  hold. Then, the distributions of the processes 
  $(\SS_n-s_n)/F_*^\leftarrow(1-1/n)$ and 
  $(\oSS_n-\osn)/F_*^\leftarrow(1-1/n)$ converge 
  to that of $k\star \dot L$, respectively in ${\rm C}[\,0,1\,]$ 
  and ${\rm D}[\,0,1\,]$, both equipped with the supremum metric.
}

\bigskip

For our application to FARIMA processes, assumption \HypKnZero\ is 
violated. The following extension of Theorem \LimitCD\ will be used.

\Corollary{\label{LimitCDCorollary}%
  If $\bigl(k_n-k_n(0)\bigr)\One_{[0,1]}$ satisfies the assumptions of 
  Theorem \LimitCD\ and \HypKnCvUnif\ holds for $(k_n)$ as well, then
  the distribution of the processes $(\SS_n-s_n)/F_*^\leftarrow(1-1/n)$ and 
  $(\oSS_n-\osn)/F_*^\leftarrow(1-1/n)$ converge to that of $k\star \dot L$
  in ${\rm D}[\,0,1\,]$ equipped with the Skorokhod topology. The limiting
  process is continuous if and only if $k(0)=0$, and, in this case, the
  convergence holds in ${\rm D}[\,0,1\,]$ equipped with the topology of 
  uniform convergence.
}

\Remark We see from Theorem \LimitCD\ and its corollary, that if $k_n(0)=0$
ultimately, then, for $n$ large, we can view $\SS_n$ in ${\rm C}[\,0,1\,]$; 
if $k_n(0)$ does
not vanish ultimately but only tends to $0$, then Corollary \LimitCDCorollary\
shows that we can view our processes in ${\rm D}[\,0,1\,]$ and still use the
supremum metric and the limiting process is still continuous. In contrast,
if $k_n(0)$ converges to a nonzero limit, then the processes must be viewed
in ${\rm D}[\,0,1\,]$ and one must use the Skorokhod topology.

\Remark If instead of using the centering 
$$
  \E X\One_{[F^\leftarrow(1/n),F^\leftarrow(1-1/n)]}(X)
$$ 
in $s_n$ we use 
$$
  \E X\One_{[-F_*^\leftarrow(1-1/n),F_*^\leftarrow(1-1/n)]}(X)
$$
then we should 
replace the process $L$ by $\tilde L$ in Theorem \LimitCD. Our normalization 
is the same as in Cs\"org\H o, Horv\`ath and Mason (1986) and appears 
more natural in the 
proofs, while the alternative centering makes the result more compatible with 
Resnick (2007) when Theorem \LimitCD\ is formally specialized to the partial 
sum process.

\bigskip

In some situations of interest, the pointwise limit of $(k_n)$ vanishes
except perhaps at the origin, and,
more strongly, $(k_n)$ converges locally uniformly to $0$ on $(0,1)$. This 
happens
for instance when there is a sequence $(\kappa_i)$ converging to $0$ such that
$$
  \overline k_n(i/n)=\kappa_i \, .
  \eqno{\equa{HypKnKappa}}
$$
In such instance, we assume that $(\kappa_i)$ tends to $0$ at least at a 
polynomial rate, in the sense that
$$
  \sum_{i\geq 0} |\kappa_i|^r <\infty\quad\hbox{ for some positive $r$.}
  \eqno{\equa{HypKappaDecay}}
$$
In this case, our next result shows that the distribution of process 
$(\oSS_n-\osn)/F_*^\leftarrow(1-1/n)$ cannot converge to that of a 
c\`adl\`ag one.
Instead, we need to use a weaker mode of convergence, identifying the process
with the point process of its trajectories, and using convergence in 
distribution of random measures. 

In the next theorem, we will use
the following notation. We write $(V_i)_{i\geq 1}$ and $(V_i')_{i\geq 1}$, 
for independent sequences of independent
random variables all having a uniform distribution on $[\,0,1\,]$. We write
$(w_i)_{i\geq 1}$ and $(w_i')_{i\geq 1}$ for two independent sequences of 
independent random variables, independent of
$(V_i)$ and $(V_i')$, and having an exponential distribution with unit mean. 
Finally we write $W_j=w_1+\cdots+w_j$ 
and $W'_j=w_1'+\cdots+w_j'$ for the corresponding partial sums. We consider
the set $M([\,0,1\,]\times \RR\setminus\{\,0\,\})$ of all measures on
$[\,0,1\,]\times \RR\setminus\{\,0\,\}$, this set being equipped with the
topology of vague convergence.

\Theorem{\label{limitPP}
  Assume that \AssumptionStable, \HypKnKappa\ and \HypKappaDecay\ hold. 
  The sequence of distributions of the point processes
  $$
    \sum_{1\leq i\leq n} \delta_{\bigl(i/n,\,{\ss \oSS_n-\osn\over
     \ss F_*^\leftarrow(1-1/n)}(i/n)\bigr)} \, ,
  $$
  converges weakly$*$ to that of the point process
  $$
    \sum_{j\geq 1} \sum_{i\geq 0} 
    \delta_{(V_j,\,p^{1/\alpha}\kappa_iW_j^{-1/\alpha})} 
    + \delta_{(V_j',\,-q^{1/\alpha}\kappa_i{W_j'}^{-1/\alpha})} \, .
  $$
}

While the form of the limiting process is similar to that obtained if $\oSS_n$
were a stationary moving average (compare with, say, Proposition 4.27 in 
Resnick, 1987), two essential differences are that Theorem \limitPP\ does not
require the sequence $(\kappa_i)$ to be summable but only to 
satisfy \HypKappaDecay,
and that the centering $s_n$ is involved in the point process.
These distinctions are essential in some applications, such as the one to 
FARIMA models which we will develop, and requires us to design new 
theoretical arguments.

\bigskip

\noindent{\bf Organization of the paper.} In the second section,
we show that if we slightly truncate the $X_i$'s in $\SS_n$, then it
is possible to center and rescale the process to obtain a limiting
Gaussian process. In the third section, we study the part of the
process left off by the truncation in the second section, namely that
where the innovations have large absolute values. We prove that the
process based on the large innovations converges to a convolution of
$k$ with a L\'evy stable process.  The combination of the second and
third section yields Theorem \LimitCD\ and Corollary \LimitCDCorollary.  
In the fourth section we
assume that \HypKnKappa\ holds. We prove that a point process based on
the part of $\oSS_n$ containing the largest $X_i$'s
converges. Combined with the result of section \fixedref{2}, this allows us to
prove Theorem \limitPP. In the fifth section, we specialize our results
to the so-called $(g,F)$-processes; this is the natural setting to
derive results which can be directly applied to FARIMA processes in
the sixth and last section.

\bigskip

Throughout the paper, we write $c$ for a generic constant whose value 
may change from occurrence to occurrence.

\bigskip


\section{Gaussian theory in a nonGaussian domain}%
In this section we show how a slight truncation of the random variables $X_i$
results in a limiting Gaussian behavior. It will appear
that some of our arguments break down outside the domain of attraction of
a nonGaussian stable law. Thus, in some sense, we are dealing with the
part of the Gaussian limiting theory which can be retained only when it fails
as a whole.

\bigskip

To truncate the random variables, we consider two deterministic sequences
of real numbers, $(a_n)$ and $(b_n)$, with
$$
  \limn a_n=-\infty \qquad\hbox{and}\qquad
  \limn b_n=+\infty \, .
  \eqno{\equa{abDiverge}}
$$
Since the belonging of $F$ to a domain of attraction of a nonGaussian stable
distribution forces a tail balance condition, it is convenient to assume that
either both sequences $(a_n)$ and $(b_n)$ grow at the same rate, or, for 
simplicity, that $(b_n)$ grows faster than $(a_n)$, that is,
$$
  \limn b_n/(-a_n) \, \hbox{ is positive or infinite.}
  \eqno{\equa{abEquiv}}
$$
Those sequences determine the truncated variables $X_i\trunc{X_i}$, whose 
means and variances are
$$
  \mu_n=\E X_i\trunc{X_i}\qquad\hbox{and}\qquad
  \sigma_n^2=\E \bigl(X_i\trunc{X_i}-\mu_n\bigr)^2 \, .
$$
This allows us to define the centered and standardized variables,
$$
  Z_{i,n}={X_i\trunc{X_i}-\mu_n\over \sigma_n} \, .
$$
They give rise to the processes containing the middle part of the innovations,
$$
  \MM_n(t)=n^{-1/2} \sum_{1\leq i\leq n} k_n\Bigl(t- {i\over n}\Bigr)Z_{i,n} 
  \, , \qquad
  0\leq t\leq 1\, .
$$
This process is continuous if and only if $k_n$ vanishes at the origin.

The main result of this section shows that a slight
truncation of the variables imposed by the divergence 
of $\bigl(n\oF(b_n)\bigr)$ and $\bigl( nF(a_n)\bigr)$ to infinity at 
any rate, arbitrarily slow, results in
normal limiting behavior. This result holds
under weaker conditions than those stated in section \fixedref{1}. In 
particular, in this section, in order to obtain a well defined limiting 
behavior, we will only assume, instead of \HypKnCvUnif, that
$$
  \hbox{ $(\okn)$ converges to $k$ in $\LTwo$.}
  \eqno{\equa{HypKnCvLTwo}}
$$
Recall that by \HypKnLin\ we assumed that $k_n$ is linearly interpolated 
between the points $(i/n)_{0\leq i\leq n}$. We write
$$
  \Delta \okn(i/n)=\okn\Bigl({i+1\over n}\Bigr)-\okn\Bigl({i\over n}\Bigr)
$$
for the jump of $\okn$ at $(i+1)/n$, or, equivalently here, up to a sign,
the
variation of $k_n$ over the interval $[\,i/n,(i+1)/n\,]$. We assume that 
the jumps are not too large in quadratic average, in the sense that
$$
  \limn n^{-1}\sum_{0\leq i< n} \bigl(\Delta \okn(i/n)\bigr)^2 = 0 \, .
  \eqno{\equa{HypKnQuadr}}
$$
It follows from Lemma \fixedref{2.1.4} below that under \HypKnQuadr, 
\HypKnCvLTwo\ is equivalent to the convergence of $(k_n)$ to $k$ in $\LTwo$.

In the next theorem, we write $W$ for a standard Wiener process and we use the
integrated process,
$$
  k\star\dot W(t)=\int k(t-s)\d W(s) \, .
$$
This integral is well defined since $k$ is supported on $[\,0,1\,]$ 
and belongs to $\LTwo$.

\Theorem{\label{MiddleCv}
  Assume that \AssumptionStable, \HypKnLin, \HypKnModulus, 
  \abDiverge--\HypKnQuadr\ hold. 
  If $\limn n\oF(b_n)=\infty$, then the distribution of the process $\MM_n$,
  viewed in ${\rm C}[\,0,1\,]$ or ${\rm D}[\,0,1\,]$ according to whether
  $k_n(0)$ vanishes or not, converges to that of $k\star \dot W$.
}

\bigskip

Since the limiting process $k\star \dot W$ is continuous, the 
conclusion of the theorem is valid when ${\rm D}[\,0,1\,]$ is equipped
with either the Skorokhod or the uniform topology.

\bigskip

The proof of Theorem \MiddleCv\ has the usual two steps, namely proving that
the finite dimensional distributions converge to the proper limiting
ones, and that the sequence of measures is tight. These two steps are
carried out in the next two subsections.

\bigskip


\def\prevs{\the\sectionnumber .\the\subsectionnumber .\the\snumber }
\def\preveq{(\the\sectionnumber .\the \subsectionnumber .\the\equanumber)}

\subsection{Proof of the finite dimensional convergence}
Throughout this subsection we assume without further notice 
that \AssumptionStable\ holds. This subsection has a similar spirit as parts
of Chakrabarty and Samorodnitsky (2010); however, we deal with weighted 
univariate sums, while they focus on multivariate unweighted sums. Clearly both
approaches can be combined to yield results on multivariate weighted sums.

Our first lemma gives an estimate 
for the mean and variance of the truncated random 
variable $X_i\trunc{X_i}$. Recall that $\alpha$ is less than $2$.

\Lemma{\label{meanVar}%
  If \abDiverge\ and \abEquiv\ hold, then $\mu_n=o(\sigma_n)$ 
  and $\sigma_n^2 \sim cb_n^2\oF(b_n)$ as $n$ tends to infinity.
}

\bigskip

The proof shows that without \abEquiv, as $n$ tends to infinity,
$$
  \sigma_n^2\sim {\alpha\over 2-\alpha} \bigl( a_n^2F(a_n)+ b_n^2\oF(b_n)\bigr)
  \, .
$$

\bigskip

\Proof Karamata's theorem (see e.g.\ Bingham, Goldie and Teugels, 1989, Theorem
1.6.4) implies that
$$
  \int_0^b x^2\d F(x)\sim {\alpha\over 2-\alpha} b^2\oF(b)
$$
as $b$ tends to infinity. Thus,
$$
  \E X^2\One_{[a_n,b_n]}(x)
  \sim {\alpha\over 2-\alpha}\bigl(a_n^2 F(a_n)+b_n^2\oF(b_n)\bigr) \, .
  \eqno{\equa{meanVarB}}
$$
Since $\alpha$ is less than $2$, both $a_n^2F(a_n)$ and $b_n^2 \oF(b_n)$
tend to infinity with $n$, so that both sides of \meanVarB\ tend to infinity
with $n$.

We now claim that, as $n$ tends to infinity,
$$
  \mu_n^2=o \bigl(b_n^2\oF(b_n)\bigr)\, .
  \eqno{\equa{meanVarA}}
$$
We concentrate on the part of the proof dealing with the upper tail, the lower
tail being handled in the same way.
If $F$ has a finite expectation, it is clear that \meanVarA\ is true.
If $\alpha$ is less than $1$, \meanVarA\ holds since, as $b$ tends to infinity,
$$
  \int_0^bx\d F(x)\sim {\alpha\over 1-\alpha} b \oF(b)
  = o\Bigl( b\sqrt{\oF(b)}\Bigr) \, .
$$
If $\alpha$ is $1$ and $F$ has infinite expectation, then for any 
fixed $\epsilon$ and as $b$ tends to infinity,
$$
  \int_0^bx\d F(x)=O(b^\epsilon) \, ,
$$
and \meanVarA\ still holds.

Combining \meanVarA\ and \meanVarB\ yields the estimates on $\mu_n$ and
$\sigma_n^2$ given in the lemma.\hfill\qed

\bigskip

Our next lemma is of independent interest. It is a Lindeberg-Feller
type result, obtained by applying the theorem it mimics. It is one of
those results we alluded to in the introduction, whose proof relies
explicitly on the distribution $F$ being in the domain of attraction
of a nonGaussian stable law by using Lemma \meanVar.

\Lemma{\label{Lindeberg}
  Suppose that \abDiverge\ and \abEquiv\ hold.
  Let $A_{i,n}$, $1\leq i\leq n$, $n\geq 1$, be a triangular array of 
  constants, such that
  $$
    \sum_{1\leq i\leq n} A_{i,n}^2=1 \, .
    \eqno{\equa{AVar}}
  $$
  Then, $\sum_{1\leq i\leq n} A_{i,n}Z_{i,n}$ has a limiting normal 
  distribution if and only if for any positive $\lambda$,
  $$
    \limn \sum_{1\leq i\leq n} A_{i,n}^2
    \One\{\, |A_{i,n}|>\lambda \sqrt{\oF(b_n)} \,\} = 0 \, .
    \eqno{\equa{LFA}}
  $$
}

\Proof The Lindeberg-Feller theorem (see e.g.\ Lo\`eve, 1960, \S 21.2) for 
triangular arrays asserts that $\sum_{1\leq i\leq n} A_{i,n}Z_{i,n}$ has
a normal limiting distribution if and only if for any positive $\lambda$,
$$
  \limn \sum_{1\leq i\leq n} \E A_{i,n}^2Z_{i,n}^2
  \One\{\, |A_{i,n}Z_{i,n}|>\lambda\,\} = 0 \, . 
  \eqno{\equa{LFB}}
$$
Note that $|A_{i,n}Z_{i,n}|$ exceeds $\lambda$ if and only if
$$
  X_i\trunc{X_i}>\mu_n + \lambda {\sigma_n\over |A_{i,n}|} \, 
$$
or
$$
  X_i\trunc{X_i}<\mu_n - \lambda {\sigma_n\over |A_{i,n}|} \, .
$$
Condition \AVar\ forces the largest $|A_{i,n}|$ to be at most $1$. Therefore,
Lemma \meanVar\ implies that
$$
  \min_{1\leq i\leq n} \lambda\sigma_n/|A_{i,n}|\gg |\mu_n|
  \eqno{\equa{LFBB}}
$$
as $n$ tends to infinity. It follows that condition \LFB\ is
equivalent to require that for any positive $\lambda$
$$
  \limn \sum_{1\leq i\leq n} A_{i,n}^2 \E Z_{i,n}^2
  \One\{\, |X_i|\trunc{X_i}>\lambda\sigma_n/|A_{i,n}|\,\} =0 \, .
  \eqno{\equa{LFC}}
$$
If $\lambda\sigma_n/|A_{i,n}|$ is at least $b_n\vee (-a_n)$, then 
the corresponding indicator function in \LFC\ vanishes. If 
instead $\lambda\sigma_n/|A_{i,n}|$ is less than $b_n\vee (-a_n)$, then 
the inequality 
$$
  |X_i|\trunc{X_i}>\lambda\sigma_n/|A_{i,n}|
  \eqno{\equa{LFCC}}
$$ 
implies
$$
  X_i^2\trunc{X_i} 
  >\lambda{\sigma_n\over |A_{i,n}|} |X_i|\trunc{X_i}
  \gg |\mu_n|\, |X_i|\trunc{X_i} \, ,
$$
the last comparison coming from \LFBB;
this asymptotic domination holds uniformly in $i$ such that \LFCC\ holds. 
We also have, under the same condition,
$$
  X_i^2\trunc{X_i}
  >\lambda^2{\sigma_n^2\over A_{i,n}^2}
  \gg \mu_n^2 \, .
$$
Hence, whenever $\lambda\sigma_n/|A_{i,n}|\leq b_n\vee (-a_n)$ 
and \LFCC\ holds,
$$
  Z_{i,n}^2={X_i^2\trunc{X_i}-2\mu_nX_i\trunc{X_i}+\mu_n^2\over\sigma_n^2}
  \sim {X_i^2\trunc{X_i}\over\sigma_n^2}
$$
as $n$ tends to infinity, uniformly in $i$ at most $n$ and such that \LFCC\ 
holds. This implies that 
condition \LFC\ is equivalent to the requirement that for any positive 
$\lambda$,
$$\displaylines{\hfill
  \limn \sum_{0\leq i\leq n} {A_{i,n}^2\over\sigma_n^2} \E X_i^2\trunc{X_i}
  \One\bigl\{\, |X_i|\trunc{X_i}>\sigma_n\lambda/|A_{i,n}|\,\bigr\}
  \cr\hfill
  \One\bigl\{\, \sigma_n\lambda < |A_{i,n}| \bigl(b_n\vee (-a_n)\bigr)\,\bigr\}
  = 0 \, .
  \qquad
  \equa{LFD}\cr
  }
$$
Since Lemma \meanVar\ implies that $\E X_i^2\trunc{X_i}\sim\sigma_n^2$ and,
using \abEquiv, that
$\sigma_n/\bigl(b_n\vee(-a_n)\bigr)\sim c\sqrt{\oF(b_n)}$, bounding 
$$
  \One\bigl\{\, |X_i|\trunc{X_i}>\sigma_n\lambda/|A_{i,n}|\,\bigr\}
$$ 
by $1$ shows that condition \LFD\ is implied by \LFA.

Conversely, suppose that \LFA\ fails. There exists a subsequence 
$(n_k)_{k\geq 1}$, a positive $\lambda_0$ and a positive $c$ such that
$$
  \limk \sum_{1\leq i\leq n_k} A_{i,n_k}^2\One\bigl\{\, |A_{i,n_k}|>\lambda_0
  \sqrt{\oF(b_{n_k})}\,\bigr\} \geq c \, . 
$$
Given Lemma \meanVar, this implies that there exists $\lambda_1$ ---~in fact, 
any $\lambda_1$ less than $\limn \lambda_0 b_n\sqrt{\oF(b_n)}/\sigma_n$ 
would do~--- such that
$$
  \limk \sum_{1\leq i\leq n_k} A_{i,n_k}^2
  \One\{\, |A_{i,n_k}|b_{n_k}>\lambda_1 \sigma_{n_k} \,\}
  \geq c \, .
  \eqno{\equa{LFE}}
$$
In particular, for $k$ large enough, there exists some $i$ such that the
indicator function involved in \LFE\ does not vanish. In what follows, we
consider such an $i$. Fix a positive 
$\lambda$ less than $\lambda_1$. For such $i$, referring to the
corresponding summand in \LFD,
$$
  \One_{[a_{n_k},b_{n_k}]}(X_i)
  \geq \One_{[(\lambda/\lambda_1)b_{n_k},b_{n_k}]}(X_i) \, ,
$$
so that
$$\displaylines{\qquad
  \One\Bigl\{\, |X_i|\One_{[a_{n_k},b_{n_k}]}(X_i) 
                > {\sigma_{n_k}\lambda\over |A_{i,n_k}|}
      \,\Bigr\}
  \hfill\cr\noalign{\vskip 3pt}\hfill
  \eqalign{
  {}\geq{}& \One\Bigl\{\, |X_i|\One_{[(\lambda/\lambda_1)b_{n_k},b_{n_k}]}(X_i)
            > {\sigma_{n_k}\lambda\over |A_{i,n_k}|}\,\Bigr\} \cr
  {}\geq{}& \One\Bigl\{\, {\lambda\over\lambda_1} b_{n_k} 
                          \One_{[(\lambda/\lambda_1)b_{n_k},b_{n_k}]}(X_i)>
            {\sigma_{n_k}\lambda\over |A_{i,n_k}|}\,\Bigr\} \, .\cr
  }
  \qquad\cr}
$$
Thus, the corresponding summand in \LFD\ is at least
$$
  {A_{i,n_k}^2\over \sigma_{n_k}^2}
  \E X_i^2 \One_{[(\lambda/\lambda_1)b_{n_k},b_{n_k}]}(X_i)
  \One\{\, \sigma_{n_k}\lambda_1 <|A_{i,n_k}| b_{n_k} \,\} \, .
  \eqno{\equa{LFF}}
$$
As in the proof of Lemma \meanVar, Karamata's theorem implies
$$
  \E X_i^2 \One_{[(\lambda/\lambda_1)b,b]}(X_i)
  \sim c b^2 \oF(b)
  \Bigl(1-\Bigl({\lambda\over\lambda_1}\Bigr)^{2-\alpha}\Bigr)
$$
as $b$ tends to infinity. Thus, for $k$ large enough, for any $i$ such
that the indicator in \LFE\ does not vanish, \LFF\ is at least a
positive constant times $A_{i,n_k}^2\One\{\, \sigma_{n_k}\lambda_1<|A_{i,n_k}|
b_{n_k}\,\}$. Combined with \LFE, this implies that
\LFD\ does not hold for any $\lambda$ small enough; thus \LFB\ does not
hold either.\hfill\qed

\bigskip

We will use the following measure theoretic lemma, where $\lambda$ is the
Lebesgue measure on $[\,0,1\,]$. We write $|f|_p$ for the 
${\rm L}^p[\,0,1\,]$-norm of a function $f$.

\Lemma{\label{LTwoZero}
  Let $(K_n)_{n\geq 1}$ be a sequence of functions in $\LTwo$ which converges
  in $\LTwo$ to a function $K$. Then, for any sequence $(u_n)_{n\geq 1}$
  diverging to $+\infty$,
  $$
    \limn \int K_n^2 \One\{\, |K_n|>u_n\,\} \d\lambda = 0 \, .
  $$
}

\Proof Under the assumption of the lemma, using the Cauchy-Schwartz and
triangle inequalities,
$$\eqalign{
  |K_n^2-K^2|_1
  &{}= |(K_n-K)(K_n+K)|_1 \cr
  &{}\leq |K_n-K|_2 |K_n+K|_2 \cr
  &{}\leq |K_n-K|_2 ( |K_n-K|_2+2|K|_2) \, .\cr
  }
$$
This implies that $(K_n^2)$ converges to $K^2$ in ${\rm L}^1[\,0,1\,]$. 
In particular the sequence $(K_n^2)$ is uniformly integrable, that is
$$
  \lim_{u\to\infty}\sup_{n\geq 1} \int K_n^2\One\{\, K_n^2>u\,\} \d\lambda
  = 0 \, ,
$$
and the result follows.\hfill\qed

\bigskip

The purpose of the next result is to show that under \HypKnQuadr, 
assumption \HypKnCvLTwo\ implies that $(k_n)$ converges to $k$ in $\LTwo$ as 
well, and that some Riemann sums related to the covariance structure of the
process $\MM_n$ converge to the corresponding integral related to the 
covariance of the process $k\star \dot W$.

\Lemma{\label{LTwoCv}
  Let $(k_n)_{n\geq 1}$ be a sequence of functions linearly interpolated 
  between points $i/n$, $0\leq i\leq n$. Assume that \HypKnCvLTwo\ and
  \HypKnQuadr\ hold.
  Then, $(k_n)$ converges to $k$ in $\LTwo$ and for any $s$ and $t$
  in $[\,0,1\,]$,
  $$
    \limn {1\over n}\sum_{1\leq i\leq n} k_n\Bigl(s-{i\over n}\Bigr)
    k_n\Bigl(t-{i\over n}\Bigr)
    =\int_0^1 k(s-u)k(t-u)\d u \, .
  $$
}

\Proof 
We write $\lceil\,\cdot\,\rceil$ for the smallest integer strictly larger
than the argument, while $\lfloor\,\cdot\,\rfloor$ denotes the largest integer
at most equal to the argument. Thus, 
$\lceil\,\cdot\,\rceil=\lfloor\,\cdot\,\rfloor+1$,
and, for instance, $\lceil 3\rceil=4$ while $\lfloor
3\rfloor = 3$. Extending the notation used in condition \HypKnQuadr, we define
$$
  \Delta \okn(s)=\okn\Bigl({\lceil ns\rceil\over n}\Bigr) 
  -\okn\Bigl({\lfloor ns\rfloor\over n}\Bigr) \, .
$$
If $s$ is in an interval $[\,i/n,(i+1)/n)$, then $\Delta \okn(s)$ is the size
of the jump of $\okn$ at $(i+1)/n$. We also introduce the fractional
part of $ns$,
$$
  \calF_n(s)=ns-\lfloor ns\rfloor \, .
$$
Since $k_n$ is obtained by linear interpolation, we have
$$
  k_n=\okn+\calF_n\Delta \okn \, .
$$

Since the nonnegative function $\calF_n$ is at most $1$, the triangle
inequality yields
$$
  |k_n-k|_2
  \leq |\okn-k|_2+|\Delta \okn|_2 \, .
$$
Assumptions \HypKnCvLTwo\ and \HypKnQuadr\ respectively ensure that 
the first and second 
summand in this upper bound tend to $0$ as $n$ tends to infinity. This proves
that $k_n$ converges to $k$ in $\LTwo$.

Our proof of the second assertion is unfortunately onerous. Consider 
the bilinear form
$$
  \langle f,g\rangle_{(s,t)}
  = {1\over n} \sum_{1\leq i\leq n} 
  f\Bigl(t-{i\over n}\Bigr)g\Bigl(s-{i\over n}\Bigr) \, ,
$$
and set $|f|_{(s)}^2=\langle f,f\rangle_{(s,s)}$. Using the bilinearity,
the triangle and Cauchy-Schwartz inequalities,
$$\displaylines{
  |\langle k_n,k_n\rangle_{(s,t)}-\langle \okn,\okn\rangle_{(s,t)}|
  \hfill\cr\noalign{\vskip 3pt}\quad
   {}= \bigl|\langle k_n-\okn,k_n-\okn\rangle_{(s,t)} 
   + \langle k_n-\okn,\okn\rangle_{(s,t)}
   + \langle \okn,k_n-\okn\rangle_{(s,t)}\bigr| 
  \hfill\cr\noalign{\vskip 3pt}\quad
   {}\leq |k_n-\okn|_{(s)} |k_n-\okn|_{(t)} + |k_n-\okn|_{(s)}|\okn|_{(t)}
   +|\okn|_{(s)} |k_n-\okn|_{(t)} \, .
  \hfill\cr}
$$
Since $|\calF_n|\leq 1$, this upper bound is at most
$$
  |\Delta\okn|_{(s)} |\Delta\okn|_{(t)}
  + |\Delta\okn|_{(s)} |\okn|_{(t)}
  + |\Delta\okn|_{(t)} |\okn|_{(s)} \, .
$$
Note that for any nonnegative $s$ less than $1$,
$$
  |\Delta\okn|_{(s)}
  = \Bigl({1\over n} \sum_{1\leq i\leq n} 
    \Delta\okn\Bigl(s-{i\over n}\Bigr)^2\Bigr)^{1/2}
  \leq \Bigl( {1\over n} \sum_{1\leq i\leq n} 
    \Delta\okn\Bigl({i\over n}\Bigr)^2\Bigr)^{1/2} \, ,
$$
while
$$
  |\okn|_{(s)} 
  =\Bigl({1\over n}\sum_{1\leq i\leq n} 
    \okn\Bigl( {\lfloor ns\rfloor-i\over n}\Bigr)^2\Bigr)^{1/2}
  \leq |\okn|_2 \, .
$$
Hence, combining \HypKnCvLTwo\ and \HypKnQuadr, we obtain that
$$
  \limn \langle k_n,k_n\rangle_{(s,t)}-\langle\okn,\okn\rangle_{(s,t)}
  = 0 \, .
  \eqno{\equa{LTwoCvAa}}
$$
Note that, writing $\langle\cdot,\cdot\rangle_2$ for the inner product 
in $\LTwo$,
$$
  \langle\okn,\okn\rangle_{(s,t)}
  = \Bigl\langle \okn\Bigl({\lfloor ns\rfloor\over n} -\cdot\Bigr) ,
                 \okn\Bigl({\lfloor nt\rfloor\over n} -\cdot\Bigr) 
  \Bigr\rangle_2 \, .
$$
Hence, to conclude the proof it suffices to show that for any $s$,
$$
  \limn \Bigl| \okn\Bigl({\lfloor ns\rfloor\over n} -\cdot\Bigr) -\okn(s-\cdot)
        \Bigr|_2 = 0
  \eqno{\equa{LTwoCvA}}
$$
and use \HypKnCvLTwo\ and \LTwoCvAa. Thus, it suffices to show that
$$
  \int_0^1\Bigl(\okn\Bigl({\lfloor ns\rfloor\over n}-u\Bigr)-\okn(s-u)\Bigr)^2
  \d u
$$
tends to $0$ uniformly in $s$ as $n$ tends to infinity. This integral is
$$
  \int_0^1\biggl(\okn\Bigl({\big\lfloor\lfloor ns\rfloor-nu\bigr\rfloor\over n}
                    \Bigr)
              -\okn\Bigl({\lfloor n(s-u)\rfloor\over n}\Bigr)\biggr)^2
  \d u \, .
  \eqno{\equa{LTwoCvB}}
$$
The equality
$$
  \lfloor n(s-u)\rfloor -1
  =\lfloor n(s-u)-1\rfloor
  \leq \bigl\lfloor \lfloor ns\rfloor -nu\bigr\rfloor
  \leq \lfloor n(s-u)\rfloor
$$
shows that \LTwoCvB\ is at most $n^{-1} \sum_{1\leq i\leq n} \Delta\okn(i/n)^2$
and hence tends to $0$ as $n$ tends to infinity under \HypKnQuadr.\hfill\qed

\bigskip

We can now prove that the finite dimensional distributions of the processes
$(\MM_n)_{n\geq 1}$ converge.

\Proposition{\label{fidi}
  Let $(k_n)$ be a sequence of functions which satisfies \HypKnCvLTwo\ 
  and \HypKnQuadr. If $\limn n\oF(b_n)=\infty$, the finite dimensional 
  distributions of the processes $\MM_n$ converge to that of $k\star\dot W$.
}

\bigskip

\Proof We will use the Cram\'er-Wold device. Let $t_1,\ldots, t_m$ be in
$[\,0,1\,]$ and let $\lambda_1,\ldots , \lambda_m$ be some real numbers.
Consider the functions 
$$
  K_n=\sum_{1\leq j\leq m} \lambda_j k_n(t_j-\cdot) \, .
$$
and
$$
  K=\sum_{1\leq j\leq m} \lambda_j k(t_j-\cdot) \, .
$$
We claim that the sequence $(K_n)_{n\geq 1}$ converges to 
$K$ in $\LTwo$.  Indeed, since both $k_n$ and $k$ are supported
in $[\,0,1\,]$,
$$
  |(k_n-k)(t-\cdot)|_2^2
  =\int_0^t(k_n-k)^2(u)\d u
  \leq |k_n-k|_2^2 \, .
$$
Combined with Lemma \LTwoCv, this shows that $k_n(t-\cdot)$ converges 
to $k(t-\cdot)$ in $\LTwo$ under \HypKnCvLTwo\ and \HypKnQuadr, and the
claim follows from the Banach space structure of $\LTwo$.

Consider the random variable
$$
  T_n=n^{-1/2}\sum_{1\leq i\leq n} K_n\Bigl({i\over n}\Bigr)Z_{i,n} \, .
$$
Its variance is
\hfuzz=1pt
$$\displaylines{\qquad
  n^{-1}\sum_{1\leq i\leq n} K_n^2\Bigl({i\over n}\Bigr)
  \hfill\cr\hfill
  {}= \sum_{1\leq j_1,j_2\leq m}\lambda_{j_1}\lambda_{j_2}\, {1\over n} 
  \sum_{1\leq i\leq n} k_n\Bigl(t_{j_1}-{i\over n}\Bigr)
  k_n\Bigl(t_{j_2}-{i\over n}\Bigr) \, .
  \qquad\cr}
$$
\hfuzz=0pt
Applying Lemma \LTwoCv, this variance converges to $\int_0^1 K^2(u)
\d u=\Var(K\star \dot W)$. 

To prove the asymptotic normality of $T_n$ it then suffices to use Lemma 
\Lindeberg\ with  $A_{i,n}=K_n(i/n)/\sqrt n$ ---~note that the sum of
the squares of $A_{i,n}$ converges to $\int_0^1 K^2(u)\d u$ by Lemma \LTwoCv, 
and though it is not equal to $1$, Lemma
\Lindeberg\ still applies for $\lambda$ is arbitrary in that lemma. Condition
\LFA\ requires that
$$
  \limn n^{-1}\sum_{1\leq i\leq n} K_n^2(i/n)
  \One\bigl\{\,|K_n|(i/n) >\lambda\sqrt{n\oF(b_n)}\,\bigr\} = 0 \, .
  \eqno{\equa{fidiA}}
$$
Since $n\oF(b_n)$ tends to infinity with $n$, \fidiA\ follows from 
Lemma \LTwoZero\ and the convergence of $(K_n)$ in $\LTwo$ and
\HypKnQuadr.\hfill\qed

\bigskip


\subsection{Proof of tightness}
We will use Kolmogorov's criterion (see e.g.\ Stroock, 1994, Theorem 3.4.16) 
to prove the tightness. This requires us to calculate higher order moments
of the increments of our processes, and our first lemma will be useful to 
do so. Recall that \abDiverge\ and \abEquiv\ hold.

\Lemma{\label{ZMoments}
  There exists some positive constants $z_k$ such that
  for any integer $k$ at least $2$,
  $$
    \E |Z_{1,n}|^k \sim z_k\oF(b_n)^{1-k/2} \, .
  $$
  as $n$ tends to infinity.
}

\bigskip

\noindent
The value of $z_k$ can be made explicit but is irrelevant in what follows.

\bigskip

\Proof Lemma \meanVar\ and the tail balance condition show that
$\sigma_n^2\sim c b_n^2\oF(b_n)$.
Karamata's theorem for Stieltjes integral (Bingham, Goldie and Teugels, 1987,
Theorem 1.6.4) implies $\E X^k\One_{[0,b]}(X) \sim c b^k\oF(b)$
as $b$ tends to infinity. Using the binomial formula,
$$
  \E Z_{i,n}^k\One_{[0,b_n]}(X_i)
  \sim c {b_n^k\oF(b_n)\over\sigma_n^k}
  \sim c \oF(b_n)^{1-k/2} \, .
$$
Similarly, since \abEquiv\ holds, $\E Z_{i,n}^k\One_{[a_n,0]}
\sim c\oF(b_n)^{1-k/2}$, and the result follows.\hfill\qed

\bigskip

Our next lemma gives a simple bound for the moment of some key increments of 
the process $\MM_n$ in terms of the modulus of continuity 
$\overline\omega_{k_n,r}$ of $k_n$ as defined in \HypKnModulus.

\Lemma{\label{ModulusUnKnA}
  For any integer $r$ there exists a constant $c_r$ such that for any $n$ 
  large enough, for any integers $0\leq p\leq q\leq n$,
  $$
    \Bigl|\E \Bigl(\MM_n\Bigl({q\over n}\Bigr)-\MM_n\Bigl({p\over n}\Bigr)
             \Bigr)^r\Bigr|
    \leq c_r \overline\omega_{k_n,r}\Bigl({q-p\over n}\Bigr)\, .
  $$
  }

\noindent Note that how large $n$ has to be for the inequality to hold depends
on $r$. Thus, this lemma does not imply that the moment generating functions
of the increments of $\MM_n$ are uniformly (in $n$) bounded in a neighborhood 
of the origin. 

\bigskip

\Proof We introduce the notation
$$
  \Delta_i=k_n\Bigl({q-i\over n}\Bigr)-k_n\Bigl({p-i\over n}\Bigr) \, .
$$
We set $T_0=0$ and for any positive integer $m$ at most $n$, we define
$$
  T_m=\sum_{1\leq i\leq m} \Delta_i Z_{i,n}
  \, .
$$
This ensures that $\MM_n(q/n)-\MM_n(p/n)=n^{-1/2}T_n$. The proof of the lemma 
is by a double induction. More precisely, defining
$$
  \omega_{m,n,r}={1\over n}\sum_{1\leq i\leq m} |\Delta_i|^r \, ,
$$
we will prove that for any integer $r$ there exists a constant $c_r$ such
that whenever $n$ is large
enough, for any integer $m$ at most $n$,
$$
  |\E (n^{-1/2}T_m)^r| \leq c_r \omega_{m,n,r} \, .
  \eqno{\equa{ModulusA}}
$$

Since $\overline\omega_{k_n,0}=1$, the conclusion of the lemma holds 
when $r$ vanishes.
It also holds when $r$ is $1$ for $\MM_n$ is centered. Similarly, since
$Z_{i,n}$ has unit variance,
$$
  \E (n^{-1/2}T_m)^2
  = n^{-1}\sum_{1\leq i\leq m} \Delta_i^2
  = \omega_{m,n,2} \, ,
$$
so that \ModulusA\ holds for $r=2$.

We now consider $r$ to be at least $3$. When $m$ is $1$, Lemma \ZMoments\
implies, for any $n$ large enough,
$$\eqalign{
  |\E (n^{-1/2}T_1)^r|
  &{}\leq n^{-r/2} |\Delta_1|^r \E |Z_{1,n}|^r \cr
  &{}\leq {|\Delta_1|^r\over n} 2z_r\bigl(n\oF(b_n)\bigr)^{1-r/2} \cr
  &{}=\omega_{1,n,r} 2 z_r\bigl( n\oF(b_n)\bigr)^{1-r/2} 
   \, . \cr
  }
$$
Since $\bigl(n\oF (b_n)\bigr)^{1-r/2}$ tends to $0$ as $n$ tends 
to infinity, this proves \ModulusA\ when $m$ is $1$.

Let $\rho$ be an integer at least $3$. We now assume that \ModulusA\ holds 
for all $r$ less than $\rho$. Then, since
$Z_{i,n}$ is centered, the binomial formula yields
$$\displaylines{
  \E(n^{-1/2}T_{k+1})^\rho
  = \E (n^{-1/2}T_k)^\rho
  \hfill\cr\hfill 
  +\sum_{0\leq j\leq \rho-2} {\rho\choose j} 
  \E (n^{-1/2}T_k)^j \Delta_{k+1}^{\rho-j} 
  n^{(j-\rho)/2}\E Z_{k+1,n}^{\rho-j}\, .
  \cr}
$$
Using Lemma \ZMoments\ and the induction hypothesis, for $n$ large enough
and $k$ between $1$ and $m-1$,
$$\displaylines{
  |\E (n^{-1/2}T_{k+1})^\rho|
  \leq |\E(n^{-1/2}T_k)^\rho|
  \hfill\cr\hfill
  {}+{1\over n}\sum_{0\leq j\leq \rho-2} {\rho\choose j} c_j \omega_{k,n,j}
  |\Delta_{k+1}|^{\rho-j} 2 z_j \bigl( n\oF(b_n)\bigr)^{1-(\rho-j)/2} \, .
  \cr}
$$
Bounding $c_jz_j$ by its maximum value over $0\leq j\leq \rho-2$, summing
these inequalities over $k$ between $1$ and $m-1$, and using 
that $\limn \bigl(n\oF(b_n)\bigr)^{1-(\rho-j)/2}\leq 1$ for 
any $0\leq j\leq \rho-2$, we obtain
$$\displaylines{\quad
  |\E (n^{-1/2}T_m)^\rho|
  \leq |\E (n^{-1/2}T_1)^\rho|
  \hfill\cr\hfill
  {}+c {1\over n}\sum_{1\leq k\leq m-1} \sum_{0\leq j\leq \rho-2}
  {\rho\choose j} \omega_{k,n,j}|\Delta_{k+1}|^{\rho-j} \, .
  \quad\cr}
$$
Considering the inner summation in this upper bound, we have
\finetune{\hfuzz=2pt}
$$\eqalign{
  \sum_{0\leq j\leq \rho-2} {\rho\choose j} \omega_{k,n,j}
  |\Delta_{k+1}|^{\rho-j}
    &{}={1\over n} \sum_{1\leq i\leq k} \sum_{0\leq j\leq \rho-2} 
             {\rho\choose j}|\Delta_i|^j |\Delta_{k+1}|^{\rho-j} \cr
    &{}\leq {1\over n} \sum_{1\leq i\leq k} (|\Delta_i|+|\Delta_{k+1}|)^\rho
             \cr
    &{}\leq {2^\rho\over n} \sum_{1\leq i\leq k}
             (|\Delta_i|^\rho+|\Delta_{k+1}|^\rho) \, . \cr}
$$
\finetune{\hfuzz=0pt}
Since $k$ is between $1$ and $m-1$, this upper bound is at most
$2^\rho(\omega_{m,n,\rho}+|\Delta_{k+1}|^\rho)$ whenever $n$ exceeds $m-1$.
Thus, since $\rho$ is at least $3$,
$$\eqalign{
  |\E (n^{-1/2}T_m)^\rho|
  &{}\leq c {|\Delta_1|^\rho\over n} + {c \over n}2^\rho\sum_{1\leq k\leq m-1} 
   (\omega_{m,n,\rho}+|\Delta_{k+1}|^\rho) \cr
  &{}\leq c 2^{\rho+1} \omega_{m,n,\rho} \, , \cr
  }
$$
proving \ModulusA\ when $r$ is $\rho$ and the lemma.\hfill\qed

\bigskip

Combining Lemma \ModulusUnKnA\ and assumption \HypKnModulus, we can bound some
moments of the increments of $\MM_n$. 

\Lemma{\label{ModulusUnKnC}
  Under \HypKnModulus, there exists a positive integer $r$ and a positive
  $\epsilon$ such that for any $0\leq s\leq t\leq 1$,
  $$
    \E|\MM_n(t)-\MM_n(s)|^r\leq c(t-s)^{1+\epsilon} \, .
  $$
}

\Proof We take $r$ to be the even integer such that \HypKnModulus\ holds.
Let $p=\lceil ns\rceil$ and $q=\lfloor nt\rfloor$. The identity
$$
  \MM_n\Bigl({p\over n}\Bigr)-\MM_n(s)
  =(p-ns)\Bigl(\MM_n\Bigl({p\over n}\Bigr)-\MM_n\Bigl({p-1\over n}\Bigr)\Bigr)
  \, ,
$$
Lemma \ModulusUnKnA\ and assumption \HypKnModulus\ imply
$$
  \E\Bigl|\MM_n\Bigl({p\over n}\Bigr)-\MM_n(s)\Bigr|^r
  \leq (p-ns)^r {c\over n^{1+\epsilon}}
  = c n^{r-1-\epsilon} \Bigl( {p\over n}-s\Bigr)^r \, .
  \eqno{\equa{ModulusUnKnCEqA}}
$$
Since $0\leq (p/n)-s\leq 1/n$, we have $n\leq 1/\bigl((p/n)-s\bigr)$. Since
$r$ is at least $2$,
taking $\epsilon$ less than $1$, we obtain that the right hand side of
\ModulusUnKnCEqA\ is at most $c\bigl((p/n)-s\bigr)^{1+\epsilon}$. 

Similarly, we have
$$
  \E\Bigl|\MM_n(t)-\MM_n\Bigl({q\over n}\Bigr)\Bigr|^r
  \leq c\Bigl(t-{q\over n}\Bigr)^{1+\epsilon} \, .
$$
Combined with Lemma \ModulusUnKnA, this yields, if $p\leq q$,
$$\eqalign{
  \E|\MM_n(t)-\MM_n(s)|^r
  &{}\leq c\Bigl( \Bigl( t-{q\over n}\Bigr)^{1+\epsilon} 
    +\Bigl({q-p\over n}\Bigr)^{1+\epsilon}
    +\Bigl({p\over n}-s\Bigr)^{1+\epsilon}\Bigr) \cr
  &{}\leq 3 c (t-s)^{1+\epsilon} \, ; \cr
  }
$$
furthermore, if $q<p$ then
$$\eqalignno{
  \E|\MM_n(t)-\MM_n(s)|^r
  &{}=(t-s)^r \E\Bigl|\MM_n\Bigl({p\over n}\Bigr)-\MM_n\Bigl({p-1\over n}\Bigr)
   \Bigr|^r \cr
  &{}\leq |t-s|^r {c\over n^{1+\epsilon}} \, .\cr}
$$
Since $n$ is at least $1$, this upper bound is at most $c |t-s|^{1+\epsilon}$, 
provided that $r$ is at least $1+\epsilon$.\hfill\qed

\bigskip

Theorem \MiddleCv\ follows from Proposition \fidi, Lemma \ModulusUnKnC, 
and Kolmogorov's criterion (see Stroock, 1994, Theorem 3.4.16). 

\bigskip

\def\prevs{\the\sectionnumber .\the\snumber }
\def\preveq{(\the\sectionnumber .\the\equanumber)}

\section{L\'evy stable limit and putting things together}%
In this section, we study the part of the process $\oSS_n$ 
---~see \procXStep~--- left unanalyzed in the previous section, namely
$$
  t\mapsto \sum_{1\leq i\leq n} \okn\Bigl(t-{i\over n}\Bigr)
  X_i\One_{(-\infty,a_n)\cup (b_n\infty)}(X_i) \, .
  \eqno{\equa{procXExtr}}
$$
In the
first subsection, we obtain the asymptotic behavior of \procXExtr\ under
some suitable conditions. Combined with Theorem \MiddleCv, this will yield
the limiting behavior of \procXStep.

\bigskip

\def\prevs{\the\sectionnumber .\the\subsectionnumber .\the\snumber }
\def\preveq{(\the\sectionnumber .\the\subsectionnumber .\the\equanumber)}

\subsection{Convoluted L\'evy stable limit}%
We consider the process containing the upper tail of the sample,
$$
  \TT_n^+(t)
  ={1\over F^\leftarrow(1-1/n)}\sum_{1\leq i\leq n} 
  k_n\Bigl(t-{i\over n}\Bigr) X_i\One_{(b_n,\infty)}(X_i) \, .
$$
Similarly, we define the process involving the lower tail of the sample,
$$
  \TT_n^-(t)
  ={1\over F^\leftarrow(1/n)}\sum_{1\leq i\leq n} 
  k_n\Bigl(t-{i\over n}\Bigr) X_i\One_{(-\infty,a_n)}(X_i) \, .
$$
It is clear that the same analysis applies to $\TT_n^+$ and $\TT_n^-$, simply
switching the lower-\ and upper-tail of the distribution.

The following result shows that whenever
$\bigl(n\oF(b_n)\bigr)$ diverges to infinity, the distribution of
$\TT_n^+$ properly centered and rescaled converges. 

We define 
$$
  \mu_n^+=\E X_i\One_{(b_n,F^\leftarrow(1-1/n))}(X_i)
$$
and
$$
  \mu_n^-=\E X_i\One_{(F^\leftarrow(1/n),a_n)}(X_i)
$$
as well as the centerings
$$
  \ttT_n^+(t)
  ={1\over F^\leftarrow(1-1/n)}\sum_{1\leq i\leq n} k_n(t-i/n) \mu_n^+
  \, .
$$
and
$$
  \ttT_n^-(t)
  ={1\over F^\leftarrow(1/n)}\sum_{1\leq i\leq n} k_n(t-i/n) \mu_n^-
  \, .
$$
We write $L^+$ for a L\'evy stable process with L\'evy measure
$$
  {\d\nu^+\over \d\lambda} (x) = \alpha x^{-\alpha-1}\One_{(0,\infty)}(x) 
  \, ;
  \eqno{\equa{nuPlusDef}}
$$
the distribution of $L^+$ is specified by the characteristic 
function
$$
  \E e^{i\theta L^+(t)}
  = \exp\Bigl( t\int \bigl( e^{i\theta x}-1-i\theta x\One_{(0,1)}(x)\bigr)
  \d\nu^+(x)\Bigr) \, .
$$

\Theorem{\label{LevyLimit}%
  Assume that \AssumptionStable, \HypKnCvUnif, \HypKnModulus\ hold
  and that $k(0)=0$. For any nonnegative
  sequence $(b_n)$ such that $\bigl(n\oF(b_n)\bigr)$ diverges to infity, 
  we can find a 
  version $\tilde \TT_n^+$ of $\TT_n^+$ as
  well as a L\'evy stable process $L^+$ such that 
  $$
    \limn |\tilde \TT_n^+-\ttT_n^+-k\star \dot L^+|_{[0,1]}=0
  $$ 
  almost surely.
}


\Remark 
Without assuming \HypKnModulus, we will prove 
that there exists a sequence $(b_n)$
with $\bigl( n\oF(b_n)\bigr)$ diverging to infinity such that the conclusion
holds. Having constructed one such sequence, we will see in the next subsection
that Theorem \MiddleCv\ ---~which requires \HypKnModulus~--- implies 
that the result is true for all sequences
with $\bigl(n\oF(b_n)\bigr)$ diverging to infinity. 

\Remark
The nonnegativity of 
$(b_n)$ is in fact irrelevant and could be replaced by the condition
that $\bigl(nF(b_n)\bigr)$ diverges to infinity; this amounts to keeping in 
$\TT_n^+$ all the innovations but the most negative ones.

\bigskip

The proof of a partial form of Theorem \LevyLimit\ makes up the remainder 
of this subsection. 
We will construct the sequence $(b_n)$ as quantiles of $F$. To do so we will
first define a sequence $(m_n)$ of positive real numbers, and 
set $b_n=F^\leftarrow(1-m_n/n)$. At an
intuitive level, the sequence $(m_n)$ gives the order of magnitude of the
number of extreme order statistics $X_{i,n}$ which we will retain in our
analysis. 

To define the sequence $(m_n)$, recall that the uniform convergence 
theorem for regularly varying functions
(Bingham, Goldie and Teugels, 1989, Theorem 1.2.1) asserts
that the convergence defining regularly varying functions is in
fact locally uniform. Our next lemma shows that for a given
regularly varying function, we can take increasing compactas in
the uniform convergence theorem.

\Lemma{\label{UCT}%
  Let $f$ be a regularly varying function of index $\beta$.
  There exists a nondecreasing
  function $m(\cdot)$ which tends to infinity at infinity and such that
  $$
    \lim_{x\to\infty}m(x)\sup_{1/m(x)\leq\lambda\leq m(x)}
    \Bigl| {f(\lambda x)\over f(x)}-\lambda^\beta\Bigr| = 0 \, .
  $$
}

\Proof For any positive integer $n$ let
$$
  h_n(x)=\sup_{y\geq x}\sup_{1/n\leq\lambda\leq n}
  \Bigl|{f(\lambda y)\over f(y)}-\lambda^\beta\Bigr| \, .
$$
For any fixed $n$ we have $\lim_{x\to\infty}h_n(x)=0$. We define an
increasing sequence $(x_n)_{n\geq 0}$ by $x_0=1$ and by induction
$$
  x_n=\inf\{\, x \, :\, h_n(x)\leq 1/n^2 \hbox{\ and \ }x\geq x_{n-1}+1\,\} 
  \, .
$$
Since $x_n\geq x_{n-1}+1$, this sequence tends to infinity. We set
$$
  m(x)=\sum_{n\geq 1} n\One_{(x_n,x_{n+1}]}(x) \, .
$$
We see that if $x_n<x\leq x_{n+1}$, then $m(x)=n$
and $h_n(x)\leq 1/n^2$.\hfill\qed

\bigskip

If $m(\cdot)$ is a function which satisfies the conclusions of Lemma \UCT, 
then any function with a slower growth does as well. Thus, given 
Lemma \UCT, we consider a divergent sequence $(m_n)$ growing slowly 
enough so that
$$
  \limn m_n\sup_{1/m_n\leq \lambda\leq m_n}
  \Bigl|{F^\leftarrow(1-\lambda/n)\over F^\leftarrow(1-1/n)}
  -\lambda^{-1/\alpha}\Bigr|
  = 0 \, .
  \eqno{\equa{mnDef}}
$$
Karamata's theorem implies that
$$
  \limt {1\over t\oF(t)} \int_{\lambda t}^t x\d F(x)
  = \alpha {\lambda^{1-\alpha}-1\over\alpha-1}
$$
locally uniformly in $\lambda$ ---~the limit should be read $-\log\lambda$ if
$\alpha$ is $1$, and what follows should be read accordingly. Thus,
$$
  \limn {n\over F^\leftarrow(1-1/n)}
  \int_{F^\leftarrow(1-\lambda/n)}^{F^\leftarrow(1-1/n)} x\d F(x)
  = \alpha {\lambda^{1-1/\alpha}-1\over \alpha-1} \, .
$$
Hence, up to making $(m_n)$ growing to infinity at an even slower rate, we
can assume that
$$
  {n\over F^\leftarrow(1-1/n)}
  \int_{F^\leftarrow(1-m_n/n)}^{F^\leftarrow(1-1/n)} x\d F(x)
  = \alpha {m_n^{1-1/\alpha}-1\over \alpha-1} +o(1)\, .
  \eqno{\equa{mnDefB}}
$$
Still making $(m_n)$ growing to infinity at an even slower rate if needed, 
we can assume that
$$
  \limn m_n\sqrt{\log\log n\over n}=0 \, .
  \eqno{\equa{mnGrowthA}}
$$
Furthermore, still making $(m_n)$ growing slower if needed, given \HypKnCvUnif,
we assume that
$$\eqalignno{
  &\limn m_n^{1-1/\alpha}|k_n-k|_{[0,1]}=0 \hbox{ if } \alpha>1\, ,
  &\equa{mnGrowthBa} \cr
  &\limn \log m_n\, |k_n-k|_{[0,1]}=0 \hbox{ if } \alpha=1 \, .
  &\equa{mnGrowthBb} \cr
  }
$$
Consider the modulus of continuity of $k$,
$$
  \omega_k(\delta)
  =\sup_{0\leq s\leq t\leq s+\delta\leq 1}|k(t)-k(s)|
  \, .
$$
Since $k$ is continuous on $[\,0,1\,]$, this modulus tends to $0$ 
with $\delta$. Hence, possibly restricting its growth, we can assume 
that $(m_n)$ satisfies
$$\eqalignno{
  &\limn m_n^{1-1/\alpha} \omega_k(n^{-1/2}\log n) = 0 
   \quad\hbox{ if } \alpha>1\, ,
  &\equa{mnGrowthCa}\cr
  &\limn \log m_n\ \omega_k(n^{-1/2}\log n) = 0
   \quad\hbox{ if } \alpha=1\ \, .
  &\equa{mnGrowthCb} \cr
  }
$$

Finally, if a sequence $(m_n)$ satisfies \mnDef--\mnGrowthCb, so is any 
sequence diverging at a slower rate, and so is any other
sequence asymptotically equivalent to it. Thus, we can assume that $1-(m_n/n)$
is in the range of $F$ for any $n$ large enough. In this case, since the 
quantile function is c\`agl\`ad, the inequality
$F^\leftarrow (1-u)>F^\leftarrow(1-m_n/n)$ is equivalent 
to $u<m_n/n$ (see e.g.\ Shorack and Wellner, 1986, \S I.1, pp.\ 5--7). We take 
$$
  b_n=F^\leftarrow(1-m_n/n) \, .
$$

Let $\tau=\tau_n$ be the permutation defining the `anti-reversed-ranks' of
the innovations, that is, writing $X_{i,n}$ for the $i$-th
order statistic of $(X_1,\ldots, X_n)$, we have
$X_{\tau(i)}=X_{n-i+1,n}$. We have
$$\eqalign{
  \TT_n^+(t)
  &{}=\sum_{1\leq i\leq n} k_n\Bigl(t-{i\over n}\Bigr) 
   {X_i\over F^\leftarrow(1-1/n)}  \ltrunc{X_i} \cr
  &{}=\sum_{1\leq i\leq n} k_n\Bigl(t-{\tau(i)\over n}\Bigr) 
   {X_{n-i+1,n}\over F^\leftarrow(1-1/n)}  \ltrunc{X_{n-i+1,n}} \, .\cr
  }
$$ 
\hfuzz=0pt
Recall that the distribution of $\tau$
is uniform on the set of permutations of $n$ elements, and, as
a random variable taking values in the symmetric group, it is
independent of the order statistics $(X_{i,n})_{1\leq i\leq
n}$.

Let $(w_i)_{i\geq 1}$ be a sequence of independent
exponential random variables with mean $1$. We set
$W_k=w_1+\cdots+w_k$ to be their partial sum. The
distributional equality
$$
  (X_{n-i+1,n})_{1\leq i\leq n}
  \eqd \Bigl(F^\leftarrow\Bigl(1-{W_i\over W_{n+1}}\Bigr)\Bigr)_{1\leq i\leq n}
  \eqno{\equa{repOrderStat}}
$$
holds ---~see e.g.\ Shorack and Wellner (1986, Chapter 8, \S 2, Proposition 1,
p.\ 335).

Define
$$A_{0,n}(t)
  =\sum_{1\leq i\leq n} k_n\Bigl(t-{\tau(i)\over n}\Bigr)
  {F^\leftarrow (1-W_i/W_{n+1})\over F^\leftarrow(1-1/n)}
  \One\Bigl\{\, {W_i\over W_{n+1}}< {m_n\over n} \,\Bigr\} \, .
$$
We then have the distributional equality of processes,
$$
  \TT_n^+\eqd A_{0,n} \, .
  \eqno{\equa{repVnAZero}}
$$
Thus, studying the limiting distribution of $\TT_n^+$ amounts to studying
that of $A_{0,n}$. For this purpose, we will approximate this latter process
by a simpler one; this will be done through four successive
approximations.

\bigskip

\noindent{\it First approximation.} Define the process
$$
  A_{1,n}(t)=\sum_{1\leq i\leq n} k_n\Bigl(t-{\tau(i)\over n}\Bigr)
  \Bigl(W_i{n\over W_{n+1}}\Bigr)^{-1/\alpha} \Utruncn \, .
$$

\Lemma{\label{AOne}%
  $\limn|A_{0,n}-A_{1,n}|_{[0,1]}=0$ almost surely.
}

\bigskip

\Proof  The strong law of large numbers ensures that
almost surely for $n$ larger than some random integer,
$$
  {1\over m_n}
  \leq {W_1\over 2}
  \leq W_1 {n\over W_{n+1}} 
  \leq W_i {n\over W_{n+1}} \, .
$$
Thus, given \mnDef, whenever $W_i/W_{n+1}\leq m_n/n$ and $n$ is 
large enough, we have
$$
\Bigl| {F^\leftarrow(1-W_i/W_{n+1})\over F^\leftarrow(1-1/n)}
   -\Bigl(W_i{n\over W_{n+1}}\Bigr)^{-1/\alpha}\Bigr|
  =o(m_n^{-1}) \, .
$$
Therefore, for any nonnegative $t$ at most $1$, and since \HypKnCvUnif\ ensures
that $(|k_n|_{[0,1]})$ is a bounded sequence
$$\eqalign{
  |A_{0,n}-A_{1,n}|(t)
  &\leq \sum_{1\leq i\leq n} |k_n|_{[0,1]} o(m_n^{-1}) \Utruncn \cr
  &=o(m_n^{-1})
   \sharp\Bigl\{\, i\geq 1\, : \, {W_i\over W_{n+1}}\leq {m_n\over n}\,\Bigr\}
   \, .\cr
}
$$
The strong law of large numbers implies that the above
cardinality is of order $m_n$.\hfill\qed

\bigskip

\noindent{\it Second approximation.} In the process $A_{1,n}$,
we would like to approximate the term
$(W_in/W_{n+1})^{-1/\alpha}$ by $W_i^{-1/\alpha}$. For this
purpose, let
$$
  A_{2,n}(t)=\sum_{1\leq i\leq n} k_n\Bigl( t-{\tau(i)\over n}\Bigr)
  W_i^{-1/\alpha} \Utruncn \, .
$$

\Lemma{\label{ATwo}%
  $\limn |A_{1,n}-A_{2,n}|_{[0,1]}=0$ almost surely.
}

\bigskip

\Proof Taylor's formula and the law of the iterated logarithm imply
$$
  \Bigl| \Bigl({n\over W_{n+1}}\Bigr)^{-1/\alpha}-1\Bigr|
  =O\Bigl(\sqrt{\log\log n\over n}\Bigr)
$$
almost surely as $n$ tends to infinity. Therefore,
$$\displaylines{\quad
  |A_{1,n}-A_{2,n}|_{[0,1]}
  \hfill\cr\hfill
  \leq |k_n|_{[0,1]} O\Bigl(\sqrt{\log\log n\over n}\Bigr)
  \sum_{1\leq i\leq n} W_i^{-1/\alpha}
  \One\{\, W_i\leq m_n W_{n+1}/n\,\} \, .
  \quad\cr}
$$
By the strong law of large numbers, $W_i\sim i$ as $i$ tends to infinity. 
Since $W_i^{-1/\alpha}\leq W_1^{-1/\alpha}$, the above upper bound is
$O(n^{-1/2}\sqrt{\log\log n}m_n)$ regardless the value of $\alpha$ in $(0,2)$.
Therefore, it converges to $0$ almost surely thanks to \mnGrowthA.\hfill\qed

\bigskip

\noindent{\it Third approximation.} In the indicator function involved 
in $A_{2,n}$, we would like to replace
the inequality $W_i/W_{n+1}\leq m_n/n$ by $W_i\leq m_n$. Define
$$
  A_{3,n}(t)
  =\sum_{1\leq i\leq n} k_n\Bigl(t-{\tau(i)\over n}\Bigr)
   W_i^{-1/\alpha}\One\{\, W_i\leq m_n\,\}\, . 
$$

\Lemma{\label{AThree}%
  $\limn |A_{2,n}-A_{3,n}|_{[0,1]}=0$ almost surely.
}

\bigskip

\Proof
The law of the iterated
logarithm implies that $|W_i- i|\leq 2\sqrt{i\log\log i}$ whenever $i$ is at
least $m_n$ and $n$ is large enough. 
In particular, $W_{n+1}/n=1+O(n^{-1}\log\log n)^{1/2}$ almost surely 
as $n$ tends to infinity.

Comparing the expressions of $A_{2,n}$ and $A_{3,n}$, we consider
the difference
$$
  \Utruncn-\One\{\, W_i\leq m_n\,\} \, .
  \eqno{\equa{AThreeA}}
$$
This is $1$ if $m_n<W_i<m_n W_{n+1}/n$ and this is $-1$ 
if $m_nW_{n+1}/n\leq W_i\leq m_n$. 
Then, the law of the iterated logarithm applied to $W_{n+1}$ and assumption
\mnGrowthA\ imply that for \AThreeA\ to be $1$, we must have $m_n<W_i<m_n+o(1)$
while for \AThreeA\ to be $-1$, we must have $m_n+o(1)<W_i<m_n$. The law
of the iterated logarithm applied to $W_i$ then shows that for \AThreeA\ not to
vanish we should have $m_n$ between $i-2\sqrt{i\log\log i}$ 
and $i+2\sqrt{i\log\log i}$. Such $i$ is between $m_n-3\sqrt{m_n\log\log m_n}$
and $m_n+3\sqrt{m_n\log\log m_n}$ whenever $n$ is large enough. Thus,
$$\eqalign{
  |A_{2,n}-A_{3,n}|(t)
  &{}\leq |k_n|_{[0,1]} \sum_{|i-m_n|\leq 3\sqrt{m_n\log\log m_n}} 
    W_i^{-1/\alpha} \cr
  &{}\leq |k_n|_{[0,1]} m_n^{-1/\alpha} \sqrt{m_n\log\log m_n}O(1) \cr
  }
$$
almost surely as $n$ tends to infinity. This upper bound tends to $0$ since
$\alpha$ is less than $2$.\hfill\qed

\bigskip

\noindent{\it Fourth approximation.} In the process $A_{3,n}$, we would like
to replace $k_n$ by $k$ and $\tau(i)/n$ by a uniform random variable. 
For this, we introduce
a sequence $(V_i)_{i\geq 1}$ of independent uniform random variables 
on $[\,0,1\,]$, independent of all random variables introduced so far. 
We write $G_n$ for the empirical distribution function of
$(V_i)_{1\leq i\leq n}$, that is
$$
  G_n(x)=n^{-1}\sum_{1\leq i\leq n} \One_{[0,x]}(V_i) \, .
$$
The vector $\bigl(nG_n(V_i)\bigr)_{1\leq i\leq n}$ is the vector of the ranks
of the $V_i$, $1\leq i\leq n$. It is uniformly distributed over the set of all
permutations of $n$ elements. Therefore, without any loss of generality, 
we can assume that
$$
  \bigl(\tau(i)\bigr)_{1\leq i\leq n}=\bigl( nG_n(V_i)\bigr)_{1\leq i\leq n}
  \, .
$$
For this version of $\tau$,
$$
  A_{3,n}(t)=\sum_{1\leq i\leq n} k_n\bigl(t-G_n(V_i)\bigr) W_i^{-1/\alpha}
  \One\{\, W_i\leq m_n\,\}
  \, .
$$
Define
$$
  A_{4,n}(t)
  =\sum_{1\leq i\leq n} k(t-V_i) W_i^{-1/\alpha}\One\{\, W_i\leq m_n\,\}
  \, .
$$

\Lemma{\label{AFour}%
  $\limn|A_{3,n}-A_{4,n}|_{[0,1]}=0$ almost surely.
}

\bigskip

\Proof Finkelstein's (1971) law of the iterated logarithm for the empirical 
process ensures that $|\sqrt n(G_n-\Id)|_{[0,1]}=O(\log\log n)^{1/2}$
almost surely. Therefore, since $k$ vanishes at $0$, almost surely 
as $n$ tends to infinity,
$$
  \max_{1\leq i\leq n}
  \bigl|k\bigl(\cdot -G_n(V_i)\bigr)-k(\cdot-V_i)\bigr|_{[0,1]}
  \leq \omega_k(n^{-1/2}\log n) \, .
  \eqno{\equa{AFourEqA}}
$$
Using again the law of large numbers, we obtain that, almost surely for $n$ 
large enough, $|A_{3,n}-A_{4,n}|_{[0,1]}$ is at most
$$
  \bigl(|k_n-k|_{[0,1]}+\omega_k(n^{-1/2}\log n)\bigr)
            \sum_{1\leq i\leq 2m_n} W_i^{-1/\alpha} \, .
$$
The strong law of large numbers ensures that
$$
  \sum_{1\leq i\leq 2m_n} W_i^{-1/\alpha}
  = \cases{ O(m_n^{1-1/\alpha}) & if $\alpha>1$,\cr
            O(\log m_n)         & if $\alpha=1$,\cr
            O(1)                & otherwise,\cr}
$$
and the result follows from conditions \mnGrowthBa--\mnGrowthCb.\hfill\qed

\bigskip

\Remark If $k$ does not vanish at $0$, the inequality \AFourEqA\ does
not hold; one needs replace $k$ by $k-k(0)$ and handle the part of the process
driven by $k(0)$ separately.  This does not introduce any further difficulty
because this part is a usual partial sum. However, the uniform topology
needs to be replaced by the Skorokhod (1956) ${\rm J}_1$ one.

\bigskip

Combining Lemmas \AOne--\AFour, we see that, as $n$ tends to infinity,
$$
  \TT_n^+\eqd A_{4,n}+o(1)
  \eqno{\equa{SnFinalApprox}}
$$
where the $o(1)$ term is uniform on $[\,0,1\,]$ and almost sure.

We now introduce the processes
$$\eqalign{
  \Delta_{1,n}(t)
    &{}=\sum_{1\leq i\leq n}(W_i^{-1/\alpha}-i^{-1/\alpha}) k(t-V_i)
     \One\{\, W_i\leq m_n\,\} \, , \cr
  \Delta_{2,n}(t)
    &{}=\sum_{1\leq i\leq n} i^{-1/\alpha}\bigl( k(t-V_i)-\E k(t-V_i)\bigr) 
     \One\{\, W_i\leq m_n\,\} \, .\cr
  }
$$
Using the law of the iterated logarithm as in the proof of Lemma \AThree,
$$\eqalignno{
  A_{4,n}
  &{}=\Delta_{1,n}+\Delta_{2,n} 
    +\sum_{1\leq i\leq n}i^{-1/\alpha}\E k(\cdot-V_i)
    \One\{\, W_i\leq m_n\,\} \, .\cr
  &{}=\Delta_{1,n}+\Delta_{2,n}+\sum_{1\leq i\leq m_n}i^{-1/\alpha}
    \E k(\cdot -V_i) + o(1) \, , \qquad\qquad
  &\equa{AFourDecomposition} \cr}
$$
where the $o(1)$ is uniform over $[\,0,1\,]$.

Not surprisingly, our next results show that, as $n$ tends to infinity,
$\Delta_{1,n}$ and $\Delta_{2,n}$ converge respectively to
$$\eqalign{
  \Delta_1(t)
    &=\sum_{i\geq 1} (W_i^{-1/\alpha}-i^{-1/\alpha})k(t-V_i) \, , \cr
  \Delta_2(t)
    &=\sum_{i\geq 1} i^{-1/\alpha} \bigl( k(t-V_i)-\E k(t-V_i) \bigr)\, .\cr
  }
$$

\Lemma{\label{DeltaConverge}
  The sequences of processes $(\Delta_{1,n})$ and $(\Delta_{2,n})$
  converge almost surely
  to respectively $\Delta_1$ and $\Delta_2$, and those limiting processes
  are continuous.
}

\bigskip

\Proof Since $m_n=o(n)$ as $n$ tends to infinity by \mnGrowthA, the 
strong law of large numbers implies that almost surely for $n$ large enough,
$$
  \Delta_{1,n}(t)
  =\sum_{i\geq 1} (W_i^{-1/\alpha}-i^{-1/\alpha}) k(t-V_i)
   \One\{\, W_i\leq m_n\,\} 
$$
and
$$
  \Delta_{2,n}(t)
  =\sum_{i\geq 1} i^{-1/\alpha}  \bigl( k(t-V_i)-\E k(t-V_i)\bigr)
   \One\{\, W_i\leq m_n\,\} \, .
$$
The strong law of large numbers and the
law of the iterated logarithm yield, uniformly on $[\,0,1\,]$,
$$\displaylines{\qquad
    \limn \Bigl|\sum_{i\geq 1} (W_i^{-1/\alpha}-i^{-1/\alpha})k(t-V_i)
    \One\{\, W_i>m_n\,\}\Bigr|
  \hfill\cr\hfill
    {}\leq 2|k|_{[0,1]} \limn \sum_{i\geq 1} 
    {\sqrt{i\log\log i}\over i^{1+(1/\alpha)}} \One\{\, i>m_n/2\,\} \, .
  \qquad\cr
  }
$$
Since $\alpha<2$, this proves the convergence of $\Delta_{1,n}$ to $\Delta_1$.

To prove the continuity of the limiting process $\Delta_1$, we have for
$0\leq s\leq t\leq 1$,
$$
  |\Delta_1(t)-\Delta_1(s)|
  \leq \sum_{i\geq 1} |W_i^{-1/\alpha}-i^{-1/\alpha}| \omega_k(t-s)
$$
and the series in this upper bound converges almost surely.

To prove the convergence of $\Delta_{2,n}$ to $\Delta_2$ requires three steps.

\noindent{\it Step 1.} Convergence of the finite dimensional distributions.
The remainder
$$
  \sum_{i\geq 1} i^{-1/\alpha} \bigl( k(t-V_i)-\E k(t-V_i)\bigr)
  \One\{\, W_i>m_n\,\}
$$
has variance at most $c\sum_{i\geq 1} i^{-2/\alpha}\Prob\{\, W_i> m_n\,\}$.
This variance tends to $0$ as $n$ tends to infinity, and hence, the finite
dimensional distributions of $(\Delta_{2,n})$ converge to those of $\Delta_2$.

\noindent{\it Step 2.} Tightness. Recall that the modulus of continuity
$\omega_{k,r}$ was defined before \HypKModulus. We first prove that 
for any integer $r$ at least $2$,
$$
  \E|\Delta_{2,n}(t)-\Delta_{2,n}(s)|^r \leq c\omega_{k,r}(t-s) \, .
  \eqno{\equa{DeltaConvergenceEqA}}
$$
We consider first the case $r=2$. Since the summands in $\Delta_{2,n}$ are
centered and independent given $(W_i)_{i\geq 1}$,
$$\displaylines{\qquad
  \E\Bigl(\bigl(\Delta_{2,n}(t)-\Delta_{2,n}(s)\bigr)^2\Bigm| 
          (W_i)_{i\geq 1}\Bigr)
  \hfill\cr\noalign{\vskip 4pt}\hfill
  {}=\sum_{i\geq 1} i^{-2/\alpha} 
              \Var\bigl( k(t-V_i)-k(s-V_i)\bigr) 
              \One\{\, W_i\leq m_n\,\} \, .
  \qquad\cr}
$$
Hence,
$$\displaylines{\quad
  \E\bigl( \Delta_{2,n}(t)-\Delta_{2,n}(s)\bigr)^2
  \hfill\cr\noalign{\vskip 5pt}\hfill
  \eqalign{
  {}\leq{}& \sum_{i\geq 1} i^{-2/\alpha} \E\bigl( k(t-V_i)-k(s-V_i)\bigr)^2
              \Prob\{\, W_i\leq m_n\,\} \cr
  {}\leq{}& \omega_{k,2}(t-s)\sum_{i\geq 1} i^{-2/\alpha} \, .\cr
  }
  \quad\cr}
$$

Assume now that $r$ is an integer greater than $2$. Given $(W_i)_{i\geq 1}$, 
the summands in $\Delta_{2,n}$ are independent and centered. Writing 
$$
  A_i=k(t-V_i)-k(s-V_i)-\E\bigl( k(t-V_i)-k(s-V_i)\bigr) \, ,
$$
the Marcinkiewicz-Zygmund inequality (see e.g.\ Chow and Teicher, 1988, 
\S 10.3, Theorem 2) implies
$$\displaylines{\quad
  \E\bigl( |\Delta_{2,n}(t)-\Delta_{2,n}(s)|^r \bigm| (W_i)_{i\geq 1}\bigr)
  \hfill\cr\noalign{\vskip 3pt}\hfill
  \leq c \E \Bigl( \Bigl|\sum_{i\geq 1} A_i^2 i^{-2/\alpha}
  \One\{\, W_i\leq m_n\,\}\Bigr|^{r/2}\Bigm| (W_i)_{i\geq 1}\Bigr) \, .
  \quad
  \equa{DeltaConvergenceEqB}\cr}
$$
Since $\alpha$ is less than $2$, the inequality $\alpha/r<1-\alpha(r-2)/(2r)$
holds. Thus, one can find a $\theta$ between $\alpha/r$ and 
$1-\alpha(r-2)/(2r)$. Such a $\theta$ is clearly between $0$ and $1$ as well.
Thus, writing $i^{-2/\alpha}$ as $i^{-2\theta/\alpha}i^{-2(1-\theta)/\alpha}$
and using H\"older's inequality with conjugate exponents $r/2$ and $r/(r-2)$,
we have
$$\displaylines{
  \sum_{i\geq 1} A_i^2 i^{-2/\alpha} \One\{\, W_i\leq m_n\,\}
  \hfill\cr\hfill
  \leq \Bigl( \sum_{i\geq 1} |A_i|^r i^{-\theta r/\alpha}\Bigr)^{2/r}
  \Bigl( \sum_{i\geq 1} i^{-2(1-\theta)r/\alpha(r-2)}\One\{\, W_i\leq m_n\,\}
  \Bigr)^{(r-2)/r} \, .
  \cr}
$$
Using this bound and taking expectation with respect to $(W_i)_{i\geq 1}$
in \DeltaConvergenceEqB, we obtain that the expectation over $(W_i)_{i\geq 1}$
of the left hand side of \DeltaConvergenceEqB\ is at most
$$
  c \E \sum_{i\geq 1} |A_i|^r i^{-\theta r/\alpha} 
  \Bigl( \sum_{i\geq 1} i^{-2(1-\theta)r/\alpha(r-2)}\Bigr)^{(r-2)/2} \, ;
$$
our choice of $\theta$ ensures that both series in this upper bound are finite
for the exponents of $i$ are both less than $-1$. Hence, the expectation 
over $(W_i)_{i\geq 1}$ of the left hand side of \DeltaConvergenceEqB\ is at 
most $c\E |A_1|^r$.

We have
$$\eqalign{
  \E |A_i|^r
  &{}\leq 2^r\Bigl( \E |k(t-V_i)-k(s-V_i)|^r  \cr
  &\hskip 90pt    + \bigl|E\bigl(k(t-V_i)-k(s-V_i)\bigr)\bigr|^r\Bigr) \cr
  &{}\leq 2^{r+1} \E |k(t-V_i)-k(s-V_i)|^r \, . \cr
  }
$$
This proves \DeltaConvergenceEqA.

Assumption \HypKModulus\ and Kolmogorov's criterion (see e.g.\ Stroock, 1994, 
Theorem 3.4.16) then yield the convergence of $\Delta_{2,m}$ to $\Delta_2$
in distribution. 

\noindent{\it Step 3.} The Ito-Nisio theorem (see e.g.\ Kwapie\'n and 
Woyczy\'nski, 1992, \S 2.1) implies that this convergence in 
distribution is in fact an almost sure one.

The continuity of $\Delta_2$ comes from assumption \HypKModulus\ and 
Kolmogorov's criterion.\hfill\qed

\bigskip

Combining \SnFinalApprox, \AFourDecomposition\ and Lemma \DeltaConverge\ 
we obtain that
$$
  \TT_n^+ -\sum_{1\leq i\leq m_n} i^{-1/\alpha} \E k(\,\cdot-V_i) 
  \eqd \Delta_1+\Delta_2 + o_P(1)
$$
where the $o_P(1)$ term is uniform.

\bigskip

It remains to identify the limiting process as indeed an integral with respect
to a L\'evy stable process. 

\bigskip

Note that
$$
  \Delta_1+\Delta_2
  =\sum_{i\geq 1}W_i^{-1/\alpha}k(\cdot-V_i)-i^{-1/\alpha}\E k(\cdot-V_i) \, .
$$
We consider the random measure
$$
  N=\sum_{i\geq 1} \delta_{(V_i,W_i^{-1/\alpha})} \, .
$$
Exercise 5.4 in Resnick (2007, \S 5.7) shows that $N$ is a Poisson random 
measure. Its mean measure is given by
$$\eqalign{
  EN\bigl([\,0,t\,]\times [\,x,\infty)\bigr)
  &=t\E\sharp\{\, i\,:\, W_i\leq x^{-\alpha}\,\}  \cr
  &=t x^{-\alpha} \, . \cr
}
$$

Consider the family of stochastic processes, for $0<\epsilon<1$,
$$
  L_\epsilon^+(s)=\int_{[0,s]\times (1,\infty)}x \d N(v,x)
  + \int_{[0,s]\times(\epsilon^{1/\alpha},1]} x \d \bigl(N-\E N\bigr)(v,x) \, .
$$
It follows from the Ito representation of L\'evy processes (see e.g.\ Resnick,
2007, \S 5.5.1) that $L^+=\lim_{\epsilon\to 0} L_\epsilon^+$ exists almost 
surely and is a L\'evy stable process, and this convergence is in fact 
locally uniform on the nonnegative real line (Ito, 1969, \S 1.7, 
Proposition 3). We rewrite $L^+_\epsilon(s)$ as
$$
  \sum_{i\geq 1} \One_{[0,s]}(V_i)\One_{[0,1/\epsilon]}(W_i) W_i^{-1/\alpha}
  -\int_{[0,s]\times (\epsilon^{1/\alpha},1)} \alpha x^{-\alpha}\d v \d x 
  \, .
$$
Agreeing to read $(\epsilon^{1/\alpha-1}-1)/(\alpha-1)$ as $-\log\epsilon$
when $\alpha$ is $1$, this yields
$$
  \d L_\epsilon^+ (s)=\sum_{i\geq 1} W_i^{-1/\alpha}\d \delta_{V_i}(s)
  \One\{\, W_i\leq 1/\epsilon\,\}-{\alpha\over \alpha-1}(\epsilon^{1/\alpha-1}-1)
  \d s \, .
$$
Therefore, whenever $0\leq t\leq 1$,
\hfuzz=2pt
$$\displaylines{
  (k\star \dot L_\epsilon^+)(t)
  \hfill\cr\ 
  {}= \int k(t-s)\d L_\epsilon^+ (s) 
  \hfill\cr\ 
  {}=\sum_{i\geq 1} k(t-V_i) W_i^{-1/\alpha} \One\{\, W_i\leq 1/\epsilon\,\}
    -{\alpha\over \alpha-1} (\epsilon^{1/\alpha-1}-1) \int k(t-s)\d s \, . 
  \hfill\cr
  }
$$
\hfuzz=0pt
Hence, we see that almost surely whenever $n$ is large enough,
$$
  A_{4,n}(t)
  =k\star \dot L^+_{1/m_n}(t)+{\alpha\over\alpha-1}(m_n^{1-1/\alpha}-1)
  \int k(t-s)\d s 
  \, .
  \eqno{\equa{AFourRepL}}
$$

In our next lemma, we show that the term added to $k\star\dot L_{1/m_n}^+$ to
obtain $A_{4,n}$ is essentially the centering  $\ttT_n^+$.

\Lemma{\label{lemmaCentering} 
  For any $\alpha$ in $(0,2)$, we have
  $$
    \limn\sup_{0\leq t\leq 1} \Bigl| \ttT_n^+(t)-{\alpha\over\alpha-1} 
    (m_n^{1-1/\alpha}-1)\int_0^1 k(t-s)\d s\Bigr| = 0 \, .
  $$
}

\Proof We give the proof when $\alpha\not=1$. With obvious 
changes, the same proof applies when $\alpha$ is $1$.

Using \mnDefB, the expectation involved in $\ttT_n^+$ satisfies
$$
  {n\over F^\leftarrow(1-1/n)}\mu_n^+
  ={\alpha\over\alpha-1}(m_n^{1-1/\alpha}-1)+o(1) \, .
  \eqno{\equa{lemmaCenteringB}}
$$

The summation involved in $\ttT_n^+$ is $\sum_{1\leq i\leq n} k_n(t-i/n)$.
Since
$$
  \Bigl|k_n\Bigl(t-{i\over n}\Bigr)-n\int_{(i-1)/n}^{i/n} k_n(t-u)\d u\Bigr|
  \leq \omega_{k_n}(1/n)
$$
and $\omega_{k_n}\leq \omega_k+2|k_n-k|_{[0,1]}$, we have
$$\displaylines{\qquad
  \sup_{0\leq t\leq 1} \Bigl|\sum_{1\leq i\leq n} k_n\Bigl(t-{i\over n}\Bigr)
  -n\int_0^1 k(t-u)\d u\Bigr|
  \hfill\cr\hfill
  \leq n\bigl(\omega_k(1/n)+3|k_n-k|_{[0,1]}\bigr) \, .
  \qquad
  \equa{lemmaCenteringC}\cr
  }
$$
Combining \lemmaCenteringB\ and \lemmaCenteringC, we obtain that, $\ttT_n^+(t)$
is equal to
$$\displaylines{
        \Bigl(\int_0^1 k(t-u)\d u
        +O\bigl(\omega_k(1/n)+|k_n-k|_{[0,1]}\bigr)\Bigr)
  \hfill\cr\hfill
        {}\times\Bigl({\alpha\over\alpha-1}(m_n^{1-1/\alpha}-1)+o(1)\Bigr) 
  \cr\hfill\cr
        {}={\alpha\over\alpha-1} (m_n^{1-1/\alpha}-1)\int_0^1 k(t-u)\d u
        +O\bigl( \omega_k(1/n)(m_n^{1-1/\alpha}+1)\bigr) 
  \qquad\cr\hfill
        {}+ O\bigl((m_n^{1-1/\alpha}+1)|k_n-k|_{[0,1]}\bigr)+o(1)\, .\cr
  }
$$
The result follows if $\alpha$ is less than $1$, and follows from 
\mnGrowthBa--\mnGrowthBb\ and \mnGrowthCa--\mnGrowthCb\ if $\alpha$ is at
least $1$.\hfill\qed

\bigskip

It remains for us to show that $L^+$ has the proper L\'evy measure $\nu^+$
defined in \nuPlusDef.

\Lemma{\label{LPlusLMeasure}
  The following equality holds
  $$
    \E \exp\bigl(itL^+(1)\bigr)
    =\exp\Bigl(\int_0^\infty \bigl(e^{itx}-1-itx\One_{(0,1)}(x)\bigr)\d\nu^+(x)
    \Bigr) \, .
  $$ 
}%

\Proof Set 
$$
  f(v,x)=x\One_{[0,1]\times [\epsilon^{1/\alpha},\infty)}(v,x)
  \hbox{ and }
  g(v,x)=x\One_{[0,1]\times[\epsilon^{1/\alpha},1]}(v,x) \, .
$$
We see that $L^+_\epsilon(1)=\int f\d N -\int g\d\E N$. Thus,
$$\displaylines{\qquad
  \E \exp\bigl(itL^+_\epsilon(1)\bigr)
  \hfill\cr\noalign{\vskip 3pt}\hfill
  \eqalign{
  {}={}&\exp\Bigl(\int e^{itf}-1\d \E N -it\int g \d\E N\Bigr) \cr
  {}={}&\exp\Bigl(\int_{\epsilon^{1/\alpha}}^\infty (e^{itx}-1)\alpha 
    x^{-\alpha-1}-it\int_{\epsilon^{1/\alpha}}^1\alpha x^{-\alpha}
    \d x \Bigr) \, . \cr}
  \qquad\cr}
$$
This characteristic function converges to the one given in the right hand
side of the equality stated in the lemma. Since the distribution
of $L^+_\epsilon$ converges to that of $L^+$ (see e.g.\ Resnick, 2007, 
\S 5.5.1), this proves the lemma.\hfill\qed

\bigskip


\subsection{Putting things together}
After completing the proof of Theorem \LevyLimit, we will combine Theorems
\MiddleCv\ and \LevyLimit\ to prove Theorem \LimitCD. Of essential importance
will be the identity
$$
  \SS_n-s_n= F^\leftarrow(1-1/n)(\TT_n^+-\ttT_n^+)
  + F^\leftarrow(1/n) (\TT_n^--\ttT_n^-)+n^{1/2}\sigma_n\MM_n \, .
  \eqno{\equa{repSn}}
$$

\bigskip

\noindent{\bf End of the proof of Theorem \LevyLimit.} In the previous 
subsection, we proved that there exists a sequence $(b_n)$ such that
$\bigl( n\oF(b_n)\bigr)$ tends to infinity and the conclusion of Theorem
\LevyLimit\ holds. To go from this existence statement to a universal one, we
keep writing $(b_n)$ for the specific sequence which we constructed, and we
write $(c_n)$ for a positive sequence such that $\bigl( n\oF(c_n)\bigr)$ tends
to infinity. We chose $a_n$ so that both $\bigl(b_n/(-a_n)\bigr)$ 
and $\bigl(c_n/(-a_n)\bigr)$ diverge to infinity, so that \abEquiv\ 
holds. We write $\TT_n^+[b_n]$ for the process $\TT_n^+$ based on 
a truncation at $(b_n)$ and $\TT_n^+[c_n]$ when the truncation is at $c_n$.
Similarly, we write $\ttT_n^+[b_n]$ and $\ttT_n^+[c_n]$ for the 
corresponding expectations, as well as $\MM_n[a_n,b_n]$ and $\MM_n[a_n,c_n]$. 
Subtracting equality \repSn\ evaluated with the sequences $(a_n)$ and $(c_n)$ 
from equality \repSn\ evaluated with the sequences $(a_n)$ and $(b_n)$,
dividing by $F^\leftarrow(1-1/n)$, we obtain
$$\displaylines{
  0 = \TT_n^+[b_n]-\ttT_n^+[b_n]-(\TT_n^+[c_n]-\ttT_n^+[c_n])
  + {n^{1/2}\sigma_n[a_n,b_n]\over F^\leftarrow(1-1/n)}\MM_n[a_n,b_n]
  \hfill\cr\hfill
  {}- {n^{1/2}\sigma_n[a_n,c_n]\over F^\leftarrow(1-1/n)}\MM_n[a_n,c_n] \, .
  \qquad\equa{EndProofA}
  \cr}
$$
Since the function $f=F^\leftarrow(1-1/\Id)/\sqrt{\Id}$ is regularly varying
of positive index $(1/\alpha)-(1/2)$, we have
$$
  {n^{1/2}\sigma_n[a_n,b_n]\over F^\leftarrow(1-1/n)}
  \sim c {n^{1/2} b_n\sqrt{\oF(b_n)}\over F^\leftarrow(1-1/n)}
  \sim c {f(n/m_n)\over f(n)}
  = o(1)
$$
as $n$ tends to infinity. Thus, Theorem \MiddleCv\ implies
$$
  {n^{1/2}\sigma_n[a_n,b_n]\over F^\leftarrow(1-1/n)}|\MM_n[a_n,b_n]|_{[0,1]}
  = o_P(1)
  \eqno{\equa{EndProofAa}}
$$
as $n$ tends to infinity, and similarly with $c_n$ substituted for $b_n$.
Thus, equality \EndProofA\ implies
$$
  \TT_n^+[b_n]-\ttT_n^+[b_n]=\TT_n^+[c_n]-\ttT_n^+[c_n] +o_P(1)
$$
where the $o_P(1)$ term is uniform in $[\,0,1\,]$, proving Theorem \LevyLimit.

\bigskip

\subsection{Proof of Theorem \LimitCD}
Our proof of Theorem \LimitCD\ requires us to establish a joint convergence of 
$(\TT_n^--\ttT_n^-,\TT_n^+-\ttT_n^+)$.

The limiting behavior of $\TT_n^--\ttT_n^-$ can be obtained from that 
of $\TT_n^+-\ttT_n^+$ by replacing the
distribution $F$ by $\M_{-1}F$. This is straightforward when $F$ is continuous,
but requires some care when it is not; this is the purpose of the following
lemma. When the distribution function of 
$X_i$ is $F$, we write $\ttT_n^+[a_n,F]$ and $\ttT_n[b_n,F]$ for what we wrote 
previously as $\ttT_n^-[a_n]$ and $\ttT_n^+[b_n]$. Similarly, 
we write $\TT_n^+[b_n,F]$, $\TT_n^-[a_n,F]$,
$\mu_n^-[a_n,F]$ and $\mu_n^+[b_n,F]$, thus
adding the underlying distribution function as a bracketted argument.

\Lemma{\label{TnPlusTnMinus}
  If $p$ and $q$ do not vanish then
  $$
    (\TT_n^--\ttT_n^-)[a_n,F]
    \eqd (\TT_n^+-\ttT_n^+)[-a_n,\M_{-1}F]\bigl(1+o(1)\bigr) +o(1) \, ,
  $$
  as $n$ tends to infinity with the $o(1)$ terms being deterministic and 
  uniform over $[0,1]$.
  }

\bigskip

\Proof Considering $\TT_n^-$, we have
$$\eqalignno{
  \TT_n^-[a_n,F](t)
  &{}={-1\over F^\leftarrow(1/n)}\sum_{1\leq i\leq n} k_n(t-i/n)(-X_i)
    \One_{(-a_n,\infty)}(-X_i) \cr
  &{}= {(\M_{-1}F)^\leftarrow(1-1/n)\over -F^\leftarrow(1/n)}
    \TT_n^+[-a_n,\M_{-1}F](t) \, .
  &\equa{TTnPlusMinus}\cr}
$$
In order to compare $\ttT_n^-[a_n,F]$ and $\ttT_n^+[-a_n,\M_{-1}F]$, we first
need to compare $\mu_n^-[a_n,F]$ and $\mu_n^+[-a_n,\M_{-1}F]$, and this requires
some discussion on $(\M_{-1}F)^\leftarrow(1-1/n)$ and $F^\leftarrow(1/n)$. Note
that since $(\M_{-1}F)(x)=\oF(-x-)$,
$$\eqalign{
  -(\M_{-1}F)^\leftarrow(1-s)
  &{}=-\inf\{\, x\,:\, F(-x-)\leq s\,\} \cr
  &{}=\sup\{\, y\,:\, F(y-)\leq s\,\} \, .\cr
  }
$$
Thus, if a real number $x$ is greater than $-(\M_{-1}F)^\leftarrow(1-s)$,
then $F(x)>s$. Given the definition of $F^\leftarrow$, this implies that
$x$ is at least $F^\leftarrow(s)$. Consequently, for any positive $s$ less than
$1$,
$$
  -(\M_{-1}F)^\leftarrow(1-s) \geq F^\leftarrow(s) \, ,
$$
and, furthermore, for any $1>t>s>1$,
$$
  F^\leftarrow(t)\geq -(\M_{-1}F)^\leftarrow(1-s) \, .
  \eqno{\equa{FMF}}
$$
Both inequalities may be strict if $F$ has a jump 
at $-(\M_{-1}F)^\leftarrow(1-s)$.

Considering now $\mu_n^-$ and a random variable $X$ having distribution 
function $F$, 
$$\eqalign{
  \mu_n^-[a_n,F]
  &{}=-\E \bigl(-X\One_{(-a_n,-F^\leftarrow(1/n))}(-X)\bigr) \cr
  &{}=-\mu_n^+[-a_n,\M_{-1}F]  \cr
  &\hskip 60pt {}
    - \E\bigl( -X\One_{[(\M_{-1}F)^\leftarrow(1-1/n),-F^\leftarrow(1/n))}(-X)
        \bigr) \cr
  &{}=-\mu_n^+[-a_n,\M_{-1}F] \cr
  &\hskip 60pt {}
      +\E\bigl( X\One_{(F^\leftarrow(1/n),-(\M_{-1}F)^\leftarrow(1-1/n)]}(X)
         \bigr) \, .\cr
  }
$$
Since $F$ is regularly varying, \FMF\ implies that for any positive $u<v<1$,
$$\displaylines{\qquad
  \bigl| \E X\One_{(F^\leftarrow(1/n),-(\M_{-1}F)^\leftarrow(1-1/n)]}(X)\bigr|
  \hfill\cr\noalign{\vskip 5pt}\hfill
  \eqalign{
    {}\leq{}& \bigl| \E X\One_{(F^\leftarrow(u/n),F^\leftarrow(v/n))}(X)\bigr|
              \cr
    {}\sim{}& \alpha\, {|F^\leftarrow(1/n)|\over n}
                {v^{1-1/\alpha}-u^{1-1/\alpha}\over \alpha-1}\, , \cr
  }
  \qquad\cr}
$$
the last fraction being interpreted as $\log (v/u)$ if $\alpha$ is $1$.
Since $u$ and $v$ are arbitrary close to $1$,
$$
  \E X\One_{(F^\leftarrow(1/n),-(\M_{-1}F)^\leftarrow(1-1/n)]}(X)
  = o(F^\leftarrow(1/n)/n)
$$
as $n$ tends to infinity. It follows that
$$
  \mu_n^-[a_n,F]=-\mu_n^+[-a_n,\M_{-1}F]+o\Bigl({F^\leftarrow(1/n)\over n}\Bigr)
  \, .
$$
Then, as $n$ tends to infinity and uniformly in $t$ in $[\,0,1\,]$,
$$\eqalign{
  \ttT_n^-[a_n,F](t)
  &{}= {1\over -F^\leftarrow(1/n)} \sum_{1\leq i\leq n} \Bigl(k_n(t-i/n)
         \mu_n^+[-a_n,\M_{-1}F] \cr
  &\hskip 160pt {}+o\Bigl({F^\leftarrow(1/n)\over n}\Bigr)\Bigr) 
         \cr
  &{}={(\M_{-1}F)^\leftarrow(1-1/n)\over -F^\leftarrow(1/n)}
         T_n^+[-a_n,\M_{-1}F]+o(1) \, . \cr
  }
$$
This, \TTnPlusMinus\ and the asymptotic 
equivalence $(\M_{-1}F)^\leftarrow(1-1/n)\sim -F^\leftarrow(1/n)$ 
yield the lemma.\hfill\qed

\bigskip

Write $\TT_n^+[b_n,F]$ and
$\TT_n^-[a_n,F]$ when the distribution of the $X_i$ is $F$. Assume until
further notice that both $p$ and $q$ do not vanish. Lemma \TnPlusTnMinus\ and
Theorem \LevyLimit\ imply that the distribution of $(\TT_n^--\ttT_n^-)[a_n,F]$
converges to that of $k\star \dot L^+$.

To obtain a joint limiting distribution
for $(\TT_n^--\ttT_n^-,\TT_n^+-\ttT_n^+)$, still assuming that $q$ does not
vanish, note that the limiting distribution
of $\TT_n^+-\ttT_n^+$ was obtained through the representation \repOrderStat\
of the order statistics. This representation allowed us to approximate
$\TT_n^+$ by a process based on $\omega_1,\ldots, \omega_i$ with $i$ of order 
$m_n$. Using the same representation for $\TT_n^-$ leads to a representation
in terms of $\omega_{n-i+1},\ldots, \omega_n$ with $i$ of order $n$. Since the
partial sums $\omega_{n-i+1}+\omega_{n-i+2}+\cdots+\omega_{n-i+k}$, 
$1\leq k\leq m_n$ are not amenable to the usual strong law of large numbers
or the law of the iterated logarithm, our proof for $\TT_n^+$ does not yields
a joint limiting result for $(\TT_n^-,\TT_n^+)$. The following modified 
construction will do. Let $(\omega_i^{(1)})_{i\geq 1}$ 
and $(\omega_i^{(2)})_{i\geq 1}$ be two independent sequences of 
random variables having a standard exponential distribution.
Set $W_i^{(1)}=\omega_1^{(1)}+\cdots+\omega_i^{(1)}$; similarly $W_i^{(2)}$
is the $i$-th partial sum of the $(\omega_i^{(2)})$. We set $W_0^{(1)}$ and
$W_0^{(2)}$ to be $0$. Define $N=\lfloor n/2\rfloor$ and define
$$
  W_i^{(3)}=\cases{ W_i^{(1)}                             & if $i\leq N$,\cr
                    \noalign{\vskip 3pt}
                    W_N^{(1)}+W_{n-N}^{(2)}-W_{n-i}^{(2)} & if $N<i\leq n$.\cr}
$$
By construction, $(W_i^{(3)})_{1\leq i\leq n}$ and $(W_i)_{1\leq i\leq n}$
involved in \repOrderStat\
have the same distribution for each fixed $n$. Note that 
$(W_i^{(3)})_{1\leq i\leq n}$ 
and $(W_n^{(3)}-W_{n-i+1}^{(3)})_{1\leq i\leq n}$ 
obey the strong law of large numbers and the
law of the iterated logarithm. Let $\omega$ be a standard exponential random
variable independent of all random variables introduced to far.
Substituting $(W_i^{(3)})_{1\leq i\leq n}$
for $(W_i)_{1\leq i\leq n}$ and $W^{(3)}_n+\omega$ for $W_{n+1}$ 
in \repOrderStat, the proof of Theorem \LevyLimit, and in 
particular the representation of $\TT_n^+$ with $A_{4,n}$, shows that 
$\TT_n^+$ and $\TT_n^-$ are asymptotically independent, and that
we can find $(a_n)$ and $(b_n)$ such that
$(\TT_n^- -\ttT_n^-,\TT_n^+-\ttT_n^+)$ converges in distribution 
to $(k\star\dot L_1^+,k\star\dot L_2^+)$, where $L_1^+$ and $L_2^+$ are 
independent copies of $L^+$.
The meaning of this convergence is that we can find versions of those processes
so that convergence holds almost surely in uniform norm. Representation 
\repSn\ and Theorem \MiddleCv ---~see also \EndProofAa~--- imply
$$\displaylines{
  {\SS_n-s_n\over F_*^\leftarrow(1-1/n)}
  = {F^\leftarrow (1-1/n)\over F_*^\leftarrow(1-1/n)} (\TT_n^+-\ttT_n^+)
  + {F^\leftarrow (1/n)\over F_*^\leftarrow(1-1/n)} (\TT_n^--\ttT_n^-)
  \hfill\cr\noalign{\vskip 5pt}\hfill
  +o_P(1) \, ,
  \qquad \equa{SnVn}
  \cr}
$$
where the $o_P(1)$ term is uniform over $[\,0,1\,]$. Since
$$
  \limn {F^\leftarrow (1-1/n)\over F_*^\leftarrow(1-1/n)} = p^{1/\alpha}
  \quad\hbox{and}\quad
  \limn {F^\leftarrow (1/n)\over F_*^\leftarrow(1-1/n)} = -q^{1/\alpha} \, ,
  \eqno{\equa{FFStarQuantile}}
$$
we obtain, using Theorem \LevyLimit, that $(\SS_n-s_n)/F_*^\leftarrow(1-1/n)$ 
converges in distribution 
to $p^{1/\alpha}k\star\dot L_2^+-q^{1/\alpha}k\star \dot L_1^+$. 

Let $L^-$ be a L\'evy process whose L\'evy measure $\nu^-$ is given by
$$
  {\d\nu^-\over \d\lambda} (x) = \alpha |x|^{-\alpha-1}\One_{(-\infty,0)}(x) 
  \, ;
$$
specifically, the characteristic function of $L^-(1)$ is
$$
  \E e^{i\theta L^-(1)}
  = \exp\Bigl( \int\bigl(e^{i\theta x}-1-i\theta x\One_{(-1,0)}(x)\bigr)
  \d\nu^-(x)\Bigr) \, .
$$
One can check that $-L^+$ and $L^-$ have the same distribution, simply by 
comparing their characteristic functions at $1$. Therefore, considering two
independent processes $L^+$ and $L^-$, we have
$$
  p^{1/\alpha}k\star\dot L_2^+-q^{1/\alpha}k\star L_1^+
  \eqd p^{1/\alpha}k\star\dot L^++q^{1/\alpha}k\star\dot L^- \, .
$$
Since $p^{1/\alpha}L^++q^{1/\alpha}L^-$ has the same distribution as 
the process $L$, it follows that the limiting process may be interpreted 
as $k\star\dot L$. This proves that, when $q$ does not vanish, the 
process $(\SS_n-s_n)/F_*^\leftarrow(1-1/n)$ converges in distribution to
$k\star\dot L$.

We assume now that $q$ vanishes. Given \SnVn, it suffices to prove that there
exists a sequence $(a_n)$ diverging to infinity, such that \abEquiv\ holds and
$$
  \limn \Bigl| {F^\leftarrow(1/n)\over F_*^\leftarrow(1-1/n)}(\TT_n^--\ttT_n^-)
  \Bigr|_{[0,1]}=0 
  \eqno{\equa{TnMinusNeglect}}
$$
in probability. Indeed, if such sequence exists, Theorem \MiddleCv\ implies 
\SnVn, and, combined with \TnMinusNeglect, Theorem  \LimitCD\ will follow.

Since $F_*^\leftarrow(1-1/n)$ is regularly varying and 
$F^\leftarrow(1/n)/F_*^\leftarrow(1-1/n)$ tends to $0$ as $n$ tends to 
infinity, we can find a sequence $(m_n)$ diverging to infinity and such that
for any positive $M$,
$$
  \limn m_n F^\leftarrow\bigl(1/(Mn)\bigr)/F_*^\leftarrow(1-1/n)= 0 \, .
  \eqno{\equa{mnTMinus}}
$$
If a sequence $(m_n)$ satisfies \mnTMinus, so does any sequence diverging at
a slower rate. Thus, we can assume that $(m_n)$ is such that $m_n/n$ is in the
range of $F$. We set $a_n=F^\leftarrow(m_n/n)$. We can then find a 
sequence $(m_n')$ diverging to infinity such
that, setting $b_n=F^\leftarrow(1-m_n'/n)$, condition \abEquiv\ holds. With
these choices of $(a_n)$ and $(b_n)$, we apply Theorems \MiddleCv\ and
\LevyLimit, so that, once we will have proved \TnMinusNeglect,  identity \SnVn\
yields Theorem \LimitCD\ when $q$ vanishes. 

For our choice of $(a_n)$ we have
$$
  \Bigl| {F^\leftarrow(1/n)\over F_*^\leftarrow(1-1/n)}\TT_n^-(t)\Bigr|_{[0,1]}
  \leq {|k_n|_{[0,1]}|X_{1,n}|\over F_*^\leftarrow(1-1/n)} 
  \sharp\{\, i\,:\, X_i\leq a_n\,\} \, .
  \eqno{\equa{TnMinusBoundA}}
$$
Let $(U_i)_{i\geq 1}$ be a sequence of independent random variables, all 
uniformly distributed over $[\,0,1\,]$. We write $(U_{i,n})_{1\leq i\leq n}$
for the order statistics of $(U_i)_{1\leq i\leq n}$.
Without loss of generality, we assume that $X_i=F^\leftarrow(U_i)$.

Let $\epsilon$ be a positive real number. Since the distribution of $nU_{1,n}$ 
converges to an exponential one, we can find a positive $M$ such that
for any $n$ large enough
$$
  \Prob\{\, U_{1,n}\leq 1/(nM)\,\} \leq \epsilon \, .
$$
Writing $G_n$ for the empirical distribution function 
of $(U_i)_{1\leq i\leq n}$ and $F_n$ for that of $(X_i)_{1\leq i\leq n}$, we
have, following Shorack and Wellner (1986, Theorem 2, p.4, and Proposition 1, 
p.5), 
$$
  \sharp\{\, i\,:\, X_i\leq a_n\,\}
  = n F_n(a_n)
  \eqd n G_n\bigl(F(a_n)\bigr)
  = n G_n(m_n/n) \, .
$$
Using Daniels' linear bound (see e.g.\ Shorack and Wellner, 1986,
Theorem 2, p.345), $G_n(m_n/n)\leq M m_n/n$ with probability at 
least $1-\epsilon$, provided $M$ is large enough. Thus, with probability at 
least $1-2\epsilon$, the right hand side of \TnMinusBoundA\ is at most
$$
  |k_n|_{[0,1]} {\bigl|F^\leftarrow\bigl(1/(nM)\bigr)\bigr|\over 
                 F_*^\leftarrow(1-1/n)} Mm_n \, .
$$
This tends to $0$ under \mnTMinus\ and \HypKnCvUnif. This proves that the left
hand side of \TnMinusBoundA\ tends to $0$ in probability. Finally,
$$\eqalign{
  \Bigl| {F^\leftarrow(1/n)\over F_*^\leftarrow(1-1/n)} \ttT_n^-(t)
  \Bigr|_{[0,1]}
  &{}\leq {n|k|_{[0,1]}\over F_*^\leftarrow(1-1/n)} 
    \E |X_i|\One_{(F^\leftarrow(1/n),a_n)}(X_i) \cr
  &{}\leq n |k|_{[0,1]} 
    {F(a_n)|F^\leftarrow(1/n)|\over F_*^\leftarrow(1-1/n)} \cr
  &{}=|k|_{[0,1]} {m_n|F^\leftarrow(1/n)|\over F_*^\leftarrow(1-1/n)} \cr
  }
$$
tends to $0$ as $n$ tends to infinity since \mnTMinus\ holds. 
This proves \TnMinusNeglect\ as well as the convergence 
in distribution of the 
process $(\SS_n-s_n)/F_*^\leftarrow(1-1/n)$ in Theorem \LimitCD, regardless
whether $q$ vanishes or not.

It remains to prove that the distribution of the piecewise constant
process $(\oSS_n-\osn)/F_*^\leftarrow(1-1/n)$ converges to the distribution
$k\star\dot L$. Since $\oSS_n-\osn$ and $\SS_n-s_n$ coincide on the points
$(i/n)_{0\leq i\leq n}$ and the latter process is a linear interpolation of the
former one on those points, we see that
$$\displaylines{\qquad
  |\SS_n-s_n-(\oSS_n-\osn)|_{[0,1]}
  \hfill\cr\hskip 60pt
  {}=\max_{0\leq i\leq n} \Bigl| (\oSS_n-\osn)\Bigl({i+1\over n}\Bigr)
    -(\oSS_n-\osn)\Bigl({i\over n}\Bigr)\Bigr| 
  \hfill\cr\hskip 60pt
  {}=\max_{0\leq i\leq n} \Bigl| (\SS_n-s_n)\Bigl({i+1\over n}\Bigr)
    -(\SS_n-s_n)\Bigl({i\over n}\Bigr)\Bigr| \, .
  \hfill\equa{EndProofB}\cr}
$$
Since the limiting process $k\star\dot L$ is continuous, \EndProofB\ implies 
that 
$$
  {|\SS_n-s_n-(\oSS_n-\osn)|_{[0,1]}\over F_*^\leftarrow(1-1/n)}=o_P(1)
  \eqno{\equa{EquivSnOSn}}
$$
as $n$ tends to infinity. Thus the weak convergence of the distribution 
of $(\SS_n-s_n)/F_*^\leftarrow(1-1/n)$
implies that of $(\oSS_n-\osn)/F_*^\leftarrow(1-1/n)$ to the same limit, 
and this proves Theorem \LimitCD.

\bigskip

\subsection{Proof of Corollary \LimitCDCorollary}
Considering $A_{i,n}=1/\sqrt n$ in Lemma \Lindeberg, we see 
that $(Z_{i,n}/\sqrt{n})_{1\leq i\leq n}$ obeys the Lindeberg condition. 
Prokhorov's theorem (1956, Theorem 3.1) implies the convergence in
distribution of the partial sum
process $t\mapsto n^{-1/2}\sum_{1\leq i\leq nt}Z_{i,n}$ to a Wiener
process. In particular, this partial sum process is stochastically bounded 
over $[\,0,1\,]$. Then, \repSn\ and \EndProofAa\ imply
$$\displaylines{\qquad
  \SS_n-s_n 
  = F^\leftarrow(1-1/n)(\TT_n^+-\ttT_n^+)+F^\leftarrow(1/n)(\TT_n^--\ttT_n^-) 
  \hfill\cr\hfill
  {}+ o_P\bigl( F_*^\leftarrow(1-1/n)\bigr)
  \qquad\cr}
$$
as $n$ tends to infinity, and the $o_P$-term is uniform over $[\,0,1\,]$. 
Define
$$
  \LL_n^+(t)
  ={1\over F^\leftarrow(1-1/n)}\sum_{1\leq i\leq nt} X_i
  \One_{(b_n,\infty)}(X_i)
$$
and the corresponding centering
$$
  \ttL_n^+(t)={\lfloor nt\rfloor \mu_n^+\over F^\leftarrow(1-1/n)} \, .
$$
We write $\TT_n^+(t)$ as the sum of $k_n(0)\LL_n^+(t)$ and
$$
  \TT_{0,n}(t)
  ={1\over F^\leftarrow(1-1/n)}\sum_{1\leq i\leq nt} 
  \Bigl(k_n\Bigl(t-{i\over n}\Bigr)-k_n(0)\Bigr) X_i\One_{(b_n,\infty)}(X_i)
  \, .
$$
We also define
$$
  \ttT_{0,n}^+(t)
  ={1\over F^\leftarrow(1-1/n)}\sum_{1\leq i\leq nt} 
  \Bigl(k_n\Bigl(t-{i\over n}\Bigr)-k_n(0)\Bigr) \mu_n^+ \, .
$$
The process $\LL_n^+$ is not amenable to 
Theorem \LevyLimit\ because assumption \HypKnModulus\ fails
if $k_n=\One_{[0,1]}$, and therefore, the following proposition requires 
a proof. 

\Proposition{\label{JointCv}
  The process $(\LL_n^+-\ttL_n^+,\TT_{0,n}^+-\ttT_{0,n}^+)$ converges in 
  distribution
  to $\bigl(L^+,(k-k(0))\star\dot L^+\bigr)$ 
  in ${\rm D}[\,0,1\,]\times{\rm C}[\,0,1\,]$ endowed with the product of the
  Skorokhod and uniform topologies.
}

\bigskip

To prove Proposition \JointCv, we will use notations copied from those of the
proof of Theorem \LevyLimit; in particular, we consider the version of the 
innovations $(X_i)_{1\leq i\leq n}$ given by \repOrderStat. For the 
corresponding version of $\TT_{0,n}^+$, the proof of Theorem \LevyLimit\ shows 
that
$$
  \limn \bigl|\TT_{0,n}^+-\ttT_{0,n}^+-\bigl(k-k(0)\bigr)\star\dot L^+
  \bigr|_{[0,1]} = 0 
$$
in probability. Thus, to prove Proposition \JointCv, we concentrate 
on $\LL_n^+$. 

With the same version of the innovations, considering $k_n=\One_{[0,1]}$
in the proof of Theorem \LevyLimit, we see that
$$
  \limn|\LL_n^+-A_{3,n}|_{[0,1]}=0
  \eqno{\equa{LnPlusAThree}}
$$
in probability. Note that for $t$ in $[\,0,1\,]$,
$$
  A_{3,n}(t)=\sum_{1\leq i\leq n} \One\{\, G_n(V_i)\leq t\,\} W_i^{-1/\alpha}
  \One\{\, W_i\leq m_n\,\} \, , 
$$
while
$$
  A_{4,n}(t)=\sum_{1\leq i\leq n} \One\{\, V_i\leq t\,\} W_i^{-1/\alpha}
  \One\{\, W_i\leq m_n\,\} \, .
$$
Recall that ${\rm D}[\,0,1\,]$ equipped with the Skorokhod topology can be 
metrized with the following distance (see Billingsley, 1968, \S 14): if $f$
and $g$ are in ${\rm D}[\,0,1\,]$, define their distance $\dsk(f,g)$ as the
infimum of all positive $\epsilon$ such that there exists a continuous and
increasing map $\lambda$ from $[\,0,1\,]$ into itself such that
$|\lambda-\Id|_{[0,1]}\leq\epsilon$ and 
$|f-g\circ\lambda|_{[0,1]}\leq\epsilon$.

\Lemma{\label{dAThreeAFour} The following limit holds in probability,
  $$
    \limn \dsk\Bigl(A_{3,n}-\ttL_n^+,A_{4,n}
    -\Id\alpha{\ds m_n^{1-1/\alpha}-1\over\ds\alpha-1}\Bigr)=0 \, .
  $$
}

\Proof Consider the set 
$$
  I_n=\{\, i \, : \, 1\leq i\leq n \, ;\, W_i\leq m_n\,\}
$$ 
and the event 
$$
  \Omega_n=\{\, \max_{i\in I_n}G_n(V_i)\leq 1-1/n\,\} \, . 
$$
We claim that $\limn \Prob(\Omega_n)=1$. Indeed, $\Omega_n$ does not occur
if and only if $\max_{i\in I_n}V_i=\max_{1\leq i\leq n} V_i$. 
Since $(V_i)$ and $(W_i)$ are independent,
$$
  \Prob\{\, \max_{i\in I_n}V_i=\max_{1\leq i\leq n} V_i\,\mid \, I_n\,\}
  =\sharp I_n/n \, .
$$
The strong law of large number ensures that $\sharp I_n\sim m_n$ almost
surely as $n$ tends to infinity. The claim follows from the 
dominated convergence theorem applied to the sequence
$\sharp I_n/n$, upon noting that this sequence is less than $1$ and converges
almost surely to $0$.

Consider now a function $\lambda_n$ such that $\lambda_n(V_i)=G_n(V_i)$
whenever $i$ is in $I_n$. Since all the $V_i$ are almost surely distinct, this
function is increasing on $\{\, V_i \,:\, i\in I_n\,\}$. We extend it 
to $[\,0,1\,]$ by defining $\lambda(0)=0$ and $\lambda(1)=1$, and, furthermore,
requiring that $\lambda_n$ is linearly interpolated between the points where
we defined it so far, that is, $0$, $\{\, V_i \, :\, i\in I_n\,\}$ and $1$. 
On $\Omega_n$ the function $\lambda_n$ is increasing. 

We see that $A_{3,n}\circ\lambda_n$ and $A_{4,n}$ are step functions whose
jumps are at the points $(V_i)_{i\in I_n}$, and that they coincide on their
jumps. Hence, they are equal.

From \mnDefB\ and, say, \mnGrowthA, we deduce that
$$\eqalign{
  \ttL_n^+(t)
  &{}={\lfloor nt\rfloor\over n} \Bigl( \alpha{m_n^{1-1/\alpha}-1\over\alpha-1}
    +o(1)\Bigr) \cr
  &{}= t\alpha {m_n^{1-1/\alpha}-1\over\alpha-1} +o(1) \cr
  }
$$
as $n$ tends to infinity, and with the $o(1)$-term being uniform over
$[\,0,1\,]$. Thus, we have
$$\displaylines{\qquad
    \Bigl| (A_{3,n}-\ttL_n^+)\circ \lambda_n-
           \Bigl( A_{4,n}-\Id\alpha{m_n^{1-1/\alpha}-1\over \alpha-1}\Bigr)
    \Bigr|_{[0,1]}
  \hfill\cr\hfill
  \eqalign{
  {}\leq{}& |\ttL_n^+\circ\lambda_n-\ttL_n^+|_{[0,1]}
           +\Bigl|\ttL_n^+-\Id\alpha{m_n^{1-1/\alpha}-1\over \alpha-1}
            \Bigr|_{[0,1]} \cr
  {}\leq{}&c|\lambda_n-\Id|_{[0,1]}(m_n^{1-1/\alpha}\vee 1)+o(1) \, . \cr
  }
  \qquad\cr
  }
$$
Finkelstein's (1971) theorem and \mnGrowthA\ imply that this upper bound tends 
to $0$
almost surely as $n$ tends to infinity. Since $|\lambda_n-\Id|_{[0,1]}$
tends to $0$ almost surely as $n$ tends to infinity, the lemma 
follows.\hfill\qed

\bigskip

\noindent{\bf Proof of Proposition \JointCv.}
Given \AFourRepL, Lemma \dAThreeAFour\ asserts
$$
  \limn \dsk(A_{3,n}-\ttL_n^+,L_{1/m_n}^+) = 0
$$
in probability. Recall (Ito, 1969, Proposition 3, \S 1.7) that
$$
  \limn|L_{1/m_n}^+-L^+|_{[0,1]}=0
  \eqno{\equa{Ito}}
$$
almost surely. In particular, $\dsk(L_{1/m_n}^+,L^+)$ converges to $0$ almost
surely as $n$ tends to infinity. 

Combining \LnPlusAThree, Lemma \dAThreeAFour\ and \Ito, we obtain that
$\limn\dsk(\LL_n^+-\ttL_n^+,L^+)=0$ in probability, and this proves Proposition
\JointCv.\hfill\qed

\bigskip

While addition of function in ${\rm D}[\,0,1\,]$ endowed with the Skorokhod 
topology is not continuous, one can check from the definition of the distance
$\dsk$ that if a sequence $(f_n)$ converges to some limit $f$ in the
Skorokhod topology and a sequence $(g_n)$ converges to a limit $g$ uniformly,
then $(f_n+g_n)$ converges to $f+g$ in the Skorokhod topology.
Then, it follows from Proposition \JointCv\ that, for the particular 
representation \repOrderStat\ of the innovations, $\TT_n^+-\ttT_n^+$ 
converges to 
$\bigl(k-k(0)\bigr)\star\dot L^++k(0)L^+$, that is, to $k\star\dot L^+$, in the
Skorokhod topology. We then argue as in the proof of Theorem \LevyLimit\ to
show that the undue restrictions on the growth of $(m_n)$ can be removed,
that $\TT_n^+$ and $\TT_n^-$ are asymptotically independent, and, finally,
that indeed the conclusion of Corollary \LimitCDCorollary\ is valid.

\bigskip


\section{Point process limit and putting things together}%
The purpose of this section is to prove Theorem \limitPP. Its structure is
very similar to that of the previous section and shows how the technique of
a separate consideration of a middle and two extreme parts can be 
applied for the slightly
different problem of proving the convergence of the distributions of 
point processes.

\bigskip


\subsection{Point process limit}%
In the context of Theorem \limitPP, we consider the weak$*$ convergence 
of the distributions of the point processes 
$\sum_{1\leq i\leq n} \delta_{(i/n,(\TT_n^+-\ttT_n^+)(i/n))}$. For notational convenience,
it is easier to think of the sequence $\kappa_i$ as a function
$$
  \kappa(x)=\kappa_{\lfloor x\rfloor} \, .
$$
The analogue of Theorem \LevyLimit\ is then the following.

\Theorem{\label{VnPP}
  Assume that \AssumptionStable, \HypKnKappa\ and \HypKappaDecay\ hold. 
  For any positive sequence $(b_n)$ such that $\limn n\oF(b_n)=\infty$, 
  the distribution of the point process 
  $$
    \sum_{1\leq i\leq n} \delta_{(i/n,(\TT_n^+- \ttT_n^+)(i/n))}
  $$ 
  converges to that
  of $\sum_{j\geq 1}\sum_{i\geq 0} \delta_{(V_j,\kappa_i W_j^{-1/\alpha})}$.
}

\bigskip

The remainder of this subsection is devoted to the proof of
this theorem for a particular sequence $(b_n)$ which allows one not to use
the centering $\ttT_n^+$. The extension to all proper sequences $(b_n)$
will be done in the next subsection. The partial result is as follows.

\Proposition{\label{VnPPPartial}
  Under \AssumptionStable\ and \HypKnKappa, there exists a sequence $(b_n)$
  such that $\limn n\oF(b_n)=\infty$ and the distribution of the point
  process
  $$
    \sum_{1\leq i\leq n} 
    \delta_{(i/n,\TT_n^+(i/n))}
  $$ 
  converges to that
  of $\sum_{j\geq 1}\sum_{i\geq 0} \delta_{(V_j,\kappa_i W_j^{-1/\alpha})}$.
}

\bigskip

We now prove this partial result.

As in the proof of Theorem \LevyLimit, we construct $b_n$ as a quantile
$F^\leftarrow(1-m_n/n)$ and impose restrictions on the growth of $(m_n)$. Our
first requirement is that $(m_n)$ satisfies \mnDef. We 
strengthen \mnGrowthA\ into
$$
  \limn {m_n\log m_n\over n^{1/4}}(\log\log n)^{1/4} = 0 \, . 
  \eqno{\equa{mnGrowthD}}
$$
This implies $\limn m_n/n^{1/4}=0$.
Furthermore, since $(\kappa_n)$ converges to $0$, so does the function
$$
  \kappa^\downarrow(x)=\sup\{\, |\kappa_n| \, :\, n\geq x\,\} \, .
$$
We then require $(m_n)$ to diverge slowly enough so that
$$
  \limn m_n\kappa^\downarrow\Bigl({n\over m_n^3}\Bigr)=0 \, .
  \eqno{\equa{mnGrowthE}}
$$

We agree that substituting $\okn$ for $k_n$ in a quantity $A_{j,n}$ defined in
section \fixedref{3.1} results in $\oA_{j,n}$.

Similarly to equality \repVnAZero, we see that $\TT_n^+$ has the same
distribution as $\oA_{0,n}$. Substituting $\okn$ for $k_n$ in the 
definition of $A_{3,n}$ in the proof of Theorem \LevyLimit, let
$$
  \oA_{3,n}(t)
  =\sum_{1\leq i\leq n}\okn\Bigl( t-{\tau(i)\over n}\Bigr) W_i^{-1/\alpha}  
  \One\{\, W_i\leq m_n\,\}
  \, .
$$
The proofs of Lemmas \AOne, \ATwo, and \AThree, show that under \mnDef\
and \mnGrowthA,
$$
  \limn |\oA_{0,n}-\oA_{3,n}|_{[0,1]}=0
$$
almost surely.

It is convenient to introduce an extra approximation,
$$
  \oA_{5,n}(t)
  =\sum_{1\leq i\leq m_n} \okn\Bigl( t-{\tau(i)\over n}\Bigr) W_i^{-1/\alpha}
  \, .
$$

\Lemma{\label{AFive}
  $\limn |\oA_{3,n}-\oA_{5,n}|_{[0,1]} = 0$ almost
  surely.
}

\bigskip

\Proof Comparing the expressions for $\overline A_{3,n}$ and 
$\overline A_{5,n}$, we consider the difference
$$
  \One\{\, W_i\leq m_n\,\}-\One\{\, i\leq m_n\,\}
$$
and argue as in the proof of Lemma \AThree.\hfill\qed

\bigskip

As in the proof of Theorem \LevyLimit, $G_n$ denotes the empirical distribution
function of $n$ independent random variables $(V_i)_{1\leq i\leq n}$
uniformly distributed over $[\,0,1\,]$ and independent from all the previous
random variables. We write $\tau(i)=nG_n(V_i)$. Extending
$\kappa$ by setting $\kappa(i)=0$ if $i$ is negative, we
have
$$
  \oA_{5,n}(t)
  =\sum_{1\leq i\leq m_n} \kappa\Bigl(n\bigl(t-G_n(V_i)\bigr)\Bigr) 
  W_i^{-1/\alpha} \, .
$$

Let $(V_{i,m_n})_{1\leq i\leq m_n}$ be the order statistics for 
$V_1,\ldots ,V_{m_n}$. Let $\theta$ be the permutation 
on $\{\, 1,2,\ldots,m_n\,\}$ such
that $V_{i,m_n}=V_{\theta(i)}$. For notational simplicity, 
define $\nu_{0,n}=0$,
$$
  \nu_{i,n}=nG_n(V_{i,m_n}) \, , \qquad 1\leq i\leq m_n  \, ,
$$
and $\nu_{m_n+1,n}=n+1$. We see that
$$
  \oA_{5,n}(i/n)=\sum_{1\leq j\leq m_n} \kappa(i-\nu_{j,n})
  W_{\theta(j)}^{-1/\alpha} \, .
$$
The intuition for what follows is that $i-\nu_{j,n}$ is either negative or 
very large for all but one $j$, the one such 
that $\nu_{j,n}\leq i<\nu_{j+1,n}$. So, except for one $j$, the 
term $\kappa(i-\nu_{j,n})$ is small and should not contribute significantly
to $\oA_{5,n}(i/n)$. To set up this approximation by a single term, for any
$i=1,2,\ldots ,n$, let $j^*=j^*(i)$ be such that
$$
  \nu_{j^*,n}\leq i<\nu_{j^*+1,n} \, .
$$
Set
$$
  B_n(i)=\cases{\kappa(i-\nu_{j^*,n}) W_{\theta(j^*)}^{-1/\alpha} 
                  & if $i\geq \nu_{1,n}$, \cr
                \noalign{\vskip 2pt}
                0 & if $i<\nu_{1,n}$. \cr}
$$
Since $\kappa(\cdot)$ vanishes on the negative half-line, we have
$$
  \oA_{5,n}(i/n)-B_n(i)
  =\sum_{1\leq j<j^*} \kappa(i-\nu_{j,n})W_{\theta(j)}^{-1/\alpha} \, .
  \eqno{\equa{approxB}}
$$
To bound $i-\nu_{j,n}$, from below, we will use the following result.

\Lemma{\label{nuSpacings}
  Assume that \mnGrowthD\ holds. Then,

  \medskip

  \noindent(i) the distribution of 
  $n^{-1}m_n^2\min_{0\leq j\leq m_n}\nu_{j+1,n}-\nu_{j,n}$ converges to 
  a standard exponential one;

  \medskip

  \noindent(ii) $\ds\limn\max_{0\leq j\leq m_n}
  \Bigl|{\ds\nu_{j+1,n}-\nu_{j,n}\over \ds n(V_{j+1,m_n}-V_{j,m_n})}-1\Bigr|=0$
  in probability;

  \medskip

  \noindent(iii) $\ds\limn\min_{0\leq j\leq m_n} \nu_{j+1,n}-\nu_{j,n}
  =+\infty$ almost surely.
}

\bigskip

\Proof (i) The central limit theorem for the empirical process implies
$$
  \nu_{j+1,n}-\nu_{j,n}=n(V_{j+1,m_n}-V_{j,m_n})+O_P(\sqrt n) \, ,
  \eqno{\equa{nuSpacingsA}}
$$
the $O_P(1)$-term being uniform in $0\leq j\leq m_n$. The result is then 
a restatement of the asymptotic behavior of the smallest spacing 
$\min_{0\leq j\leq m_n}V_{j+1,m_n}-V_{j,m_n}$ (see Shorack and Wellner,
1986, chapter 21, \S 2, Theorem 1, p.\ 726).

\noindent(ii) Given \nuSpacingsA, assertion (i) of the current
lemma and the convergence of the distribution of the normalized smallest 
spacing to an exponential one,
$$\eqalign{
  \Bigl|{\nu_{j+1,n}-\nu_{j,n}\over n(V_{j+1,m_n}-V_{j,m_n})}-1\Bigr|
  &{}=O_P\Bigl({\sqrt n\over n\min_{0\leq j\leq m_n} (V_{j+1,m_n}-V_{j,m_n})}
   \Bigr)\cr
  &{}=O_P\Bigl({m_n^2\over \sqrt n}\Bigr) \, .\cr
  }
$$
The result follows from \mnGrowthD.

\noindent (iii) Finkelstein's (1971) law of the iterated logarithm implies
$$
  \min_{0\leq j\leq m_n}\nu_{j+1,n}-\nu_{j,n}
  = n\min_{0\leq j\leq m_n} \Bigl( V_{j+1,m_n}-V_{j,m_n}
  +O\Bigl({\log\log n\over n}\Bigr)^{1/2}\Bigr)
  \eqno{\equa{nuSpacingsEqB}}
$$
almost surely. Theorem 3.1 in Devroye (1982) implies
$$
  \limn m_n^2 \log^2 m_n \min_{0\leq j\leq m_n} V_{j+1,m_n}-V_{j,m_n}=+\infty
$$
almost surely. Hence, using \mnGrowthD, the right hand side of \nuSpacingsEqB\
is at least
$$
  n\Bigl( {1\over m_n^2\log^2 m_n}+O\Bigl({\log\log n\over n}\Bigr)^{1/2}\Bigr)
  \sim {n\over m_n^2\log^2 m_n} \, .
$$
Under \mnGrowthD, this diverges to infinity.\hfill\qed

\bigskip

Our next result shows that $\oA_{5,n}(i/n)$ is well approximated by $B_n(i)$
uniformly in $i$.

\Lemma{\label{ABApprox}
  $\limn \max_{1\leq i\leq n} |\oA_{5,n}(i/n)-B_n(i)|=0$ in probability.
}

\bigskip

\Proof Given \approxB\ and that $j^*$ is at most $m_n$,
$$\eqalignno{
  |\oA_{5,n}(i/n)-B_n(i)|
  &{}\leq \sum_{1\leq j< j^*} \kappa^\downarrow(\nu_{j^*,n}-\nu_{j,n})
    W_{\theta(j)}^{-1/\alpha}\cr
  &{}\leq m_n \kappa^\downarrow\bigl(\min_{1\leq j< j^*} 
    (\nu_{j^*,n}-\nu_{j,n})\bigr)W_1^{-1/\alpha} \, .\qquad
  &\equa{ABApproxA}\cr
  }
$$
Thus, given Lemma \nuSpacings, for any positive $\epsilon$ we can find a 
positive $M$ such that with probability at least $1-\epsilon$, the upper bound
in \ABApproxA\ is at most $m_n\kappa^\downarrow(n/Mm_n^2)W_1^{-1/\alpha}$.
Assumption \mnGrowthE\ implies that this upper bound tends to $0$ 
with $n$.\hfill\qed

\bigskip

Our next lemma proves the convergence of the point process based on $B_n$.

\Lemma{\label{ppBnCv}
  The distribution of the point process
  $$
    \sum_{1\leq i\leq n} \delta_{(i/n,B_n(i))}
  $$
  converges to that of\quad
  $\ds\sum_{j\geq 1} \sum_{i\geq 0} \delta_{(V_j,\kappa_i W_j^{-1/\alpha})}$.
}

\bigskip

\Proof Write $N_n$ for the point process 
$\sum_{1\leq i\leq n} \delta_{(i/n,B_n(i))}$, and let $f$ be a function with
compact support in $[\,0,1\,]\times (\RR\setminus\{\,0\,\})$. Writing $f$ as
$f\One_{[0,1]\times (-\infty,0)}+f\One_{[0,1]\times (0,\infty)}$, there is
no loss of generality in what follows to assume that the support of $f$ is
a compact subset of $[\,0,1\,]\times (0,\infty)$. Thus, this support lies in a
strip $[\,0,1\,]\times (\epsilon,\infty)$ for some positive $\epsilon$.
The strong law of large numbers guarantees the existence of a random 
integer $j_0$ such that
$W_j^{-1/\alpha}\leq \epsilon /(2|\kappa|_\infty)$ for any $j$ at least $j_0$.
Since $\kappa$ tends to $0$ at infinity and the sequence $(W_j^{-1/\alpha})$ is
almost surely bounded, there exists a random $i_0$ such that the
function $\kappa$ is at most $\epsilon/(2\max_{j\geq 1} W_j^{-1/\alpha})$ 
on $[\,i_0,\infty)$.
Since $\min_{1\leq j\leq m_n} \nu_{j,n}-\nu_{j-1,n}$ tends to infinity almost
surely, as shown in Lemma \nuSpacings.(iii), 
this implies that there exists a random $n_0$ such that for any $n$
at least $n_0$, for a point $\bigl(i/n,B_n(i)\bigr)$ to be in the support
of $f$ we must have $i-\nu_{j^*,n}\leq i_0$ and $\theta(j)\leq j_0$. Thus,
for $n$ at least $n_0$, and since $B_n(i)$ vanishes for $i<\nu_{1,n}$,
$$\eqalignno{
  &\hskip -23pt\sum_{1\leq i\leq n} \delta_{(i/n,B_n(i))}f
  \cr
  {}={}&\sum_{1\leq j\leq m_n}\sum_{\nu_{j,n}\leq i<\nu_{j+1,n}} 
        f\bigl( i/n,\kappa(i-\nu_{j,n})W_{\theta(j)}^{-1/\alpha}\bigr) \cr
  {}={}&\sum_{1\leq j\leq m_n} \One\{\, \theta(j)\leq j_0\,\} 
        \sum_{0\leq i\leq i_0} 
        f\Bigl({\nu_{j,n}+i\over n} ,\kappa(i)W_{\theta(j)}^{-1/\alpha}
         \Bigr)\, .\qquad
       &\equa{ppBnCvA}\cr
  }%
$$
Since $f$ is uniformly continuous and only $j_0$ integers $j$ are such that
$\theta(j)\leq j_0$, \ppBnCvA\ implies
$$
  N_n f=\sum_{1\leq j\leq m_n} \One\{\, \theta(j)\leq j_0\,\}
  \sum_{0\leq i\leq i_0} 
  f\Bigl( {\nu_{j,n}\over n},\kappa_i W_{\theta(j)}^{-1/\alpha}\Bigr)
  + o(1)
$$
almost surely. The 
set $\{\, \theta(j)\,:\,1\leq j\leq m_n\,\}$ covers $\{\, 1,\ldots ,j_0\,\}$
whenever $m_n$ exceeds $j_0$. 
Moreover, by the Glivenko-Cantelli theorem,
$$
  |n^{-1}\nu_{j,n}-V_{\theta(j)}|=|n^{-1}\nu_{j,n}-V_{j,m_n}|=o(1)
$$
almost surely. Thus,
$$
  N_nf=\sum_{1\leq j\leq j_0} \sum_{0\leq i\leq i_0} 
  f(V_j,\kappa_i W_j^{-1/\alpha}) + o(1) \, .
  \eqno{\equa{ppBnCvB}}
$$
Given how $i_0$ and $j_0$ were defined, this implies
$$
  N_nf=\sum_{j\geq 1} \sum_{i\geq 0} f(V_j,\kappa_i W_j^{-1/\alpha})
  + o(1)
$$
almost surely. This proves a stronger result than stated in the 
lemma.\hfill\qed

\bigskip

We can now conclude the proof of Proposition \VnPPPartial. 
Since $A_{0,n}$ and $\oA_{0,n}$ coincide on the points $(i/n)_{1\leq i\leq n}$,
equality \repVnAZero\ 
shows that, for each $n$, the tuple $\bigl(\TT_n^+(i/n)\bigr)_{1\leq i\leq n}$ 
has the same distribution
as $\bigl( \oA_{0,n}(i/n)\bigr)_{1\leq i\leq n}$. 
Lemmas \AFive\ and \ABApprox\ imply that
$\max_{1\leq i\leq n} |\oA_{0,n}(i/n)-B_n(i)|$ tends to $0$ in probability
as $n$ tends to infinity. Then, it follows from \ppBnCvB\ that
$$
  \sum_{1\leq i\leq n} \delta_{(i/n,A_{0,n}(i/n))}f-N_nf = o_P(1)
$$
as $n$ tends to infinity. Combined with Lemma \ppBnCv, this 
implies Proposition \VnPPPartial.

\bigskip


\subsection{Putting things together}%
Combining Theorem \MiddleCv\ and Proposition \VnPPPartial, we complete the 
proofs of Theorems \VnPP.  

\bigskip

\noindent{\bf Proof of Theorem \VnPP}. As in section \fixedref{3.2}, we write
$\TT_n^+[b_n]$ and so on for the process $\TT_n^+$ based on a truncation
$(b_n)$ constructed in the previous section, which guarantees that the 
conclusion of Proposition \VnPPPartial\ holds.

Our first lemma shows that the process $\MM_n$, properly rescaled, can be 
ignored.

\Lemma{\label{UnNeglect}
  $\limn {\ds\sqrt n\sigma_n\over \ds F_*^\leftarrow(1-1/n)}|\MM_n|_{[0,1]}
  =0$ in probability.
}

\bigskip

\Proof From the proof of Lemma \ModulusUnKnA, specifically 
inequality \ModulusA, we deduce that for any $0\leq p\leq q\leq n$,
$$\eqalign{
  \Bigl|\E\Bigl(\MM_n\Bigl({q\over n}\Bigr)-\MM_n\Bigl({p\over n}\Bigr)\Bigr)^r
  \Bigr|
  &{}\leq {c_r\over n} \sum_{1\leq i\leq n} \Bigl| k_n\Bigl({q-i\over n}\Bigr)
    - k_n\Bigl({p-i\over n}\Bigr)\Bigr|^r \cr
  &{}\leq {c_r\over n} \sum_{1\leq i\leq n} |\kappa_{q-i}-\kappa_{p-i}|^r \cr
  &{}\leq 2^{r+1} {c_r\over n} \sum_{1\leq i\leq n} |\kappa_i|^r \, .\cr
  }
$$
Hence, using \HypKappaDecay, taking $r$ even and large enough so 
that $\sum_{i\geq 1} |\kappa_i|^r$ is finite, we obtain
$$
  \E\Bigl|\MM_n\Bigl({q\over n}\Bigr)-\MM_n\Bigl({p\over n}\Bigr)\Bigr|^r
  \leq {c\over n} \, .
$$
It then follows from Chebychef's inequality that
$$
  \Prob\Bigl\{\, \Bigl|\MM_n\Bigl({q\over n}\Bigr)-\MM_n\Bigl({p\over n}\Bigr)
                 \Bigr| >t\,\Bigr\} \leq {c\over n t^r} \, .
$$
Thus, Bonferroni's inequality implies
$$
  \Prob\Bigl\{\, \max_{0\leq q\leq n} \Bigl|\MM_n\Bigl({q\over n}\Bigr)
  -\MM_n(0)\Bigr|>t\,\Bigr\} \leq {c\over t^r} \, .
$$
Since $\MM_n(0)=0$ and the process $\MM_n$ is linearly interpolated between
the points $(i/n)_{0\leq i\leq n}$, this implies that $|\MM_n|_{[0,1]}=O_P(1)$
as $n$ tends to infinity. Since $\sqrt n\sigma_n/F^\leftarrow(1-1/n)$ tends
to $0$ as $n$ tends to infinity, the result follows.\hfill\qed

\bigskip

Our next lemma shows that it is possible to ignore the centering.

\Lemma{\label{centeringNeglect}
  We can choose $(b_n)$ so that Proposition \VnPPPartial\ holds 
  and $|\ttT_n^+|_{[0,1]}=o(1)$ as $n$ tends to infinity.
}

\bigskip

\Proof We have
$$\eqalign{
  \ttT_n^+(t)
  &{}={1\over F^\leftarrow(1-1/n)} \sum_{1\leq i\leq n} \kappa(nt-i)\mu_n^+ \cr
  &{}\leq c {n\mu_n^+\over F^\leftarrow (1-1/n)}
    {1\over n}\sum_{1\leq i\leq n} |\kappa(i)| \, . \cr
  }
  \eqno{\equa{togetherA}}
$$
Since the sequence $(\kappa_i)$ tends to $0$ 
and $n\mu_n^+/F^\leftarrow(1-1/n)$ is the left hand side of \mnDefB, the 
result follows from \mnDefB\ by imposing $(m_n)$ to grow slowly enough
when $\alpha$ is at least $1$.\hfill\qed

\bigskip

Representation \repSn\ and Lemma \UnNeglect\ imply
$$
  {\oSS_n-\osn\over F_*^\leftarrow(1-1/n)}
  =  p^{1/\alpha}(\TT_n^+-\ttT_n^+)- q^{1/\alpha}(\TT_n^- -\ttT_n^-)
  + o_P\bigl(1)
  \eqno{\equa{togetherB}}
$$
where the $o_P$ term is uniform over $[\,0,1\,]$. Subtracting this equality
with the truncations $(a_n)$ and $(b_n)$ from the same equality with the
truncations $(a_n)$ and $(c_n)$, we obtain
$$
  \bigl|\TT_n^+[c_n]-\ttT_n^+[c_n]-(\TT_n^+[b_n]-\ttT_n^+[b_n])\bigr|_{[0,1]}
  = o_P(1)
$$
as $n$ tends to infinity. Thus, Lemma \centeringNeglect\ yields
$$
  \TT_n^+[c_n]-\ttT_n^+[c_n]=\TT_n^+[b_n]+o_P(1) \, .
$$
Then, the same arguments as in the end of the proof of Proposition 
\VnPPPartial\ ---~see the proof of Lemma \ppBnCv\ and after~--- imply that
the distribution of 
$\sum_{1\leq i\leq n}\delta_{(i/n,(\TT_n^+[c_n]-\ttT_n^+[c_n])(i/n))}$
and that of $\sum_{1\leq i\leq n}\delta_{(i/n,\TT_n^+[b_n](i/n))}$ have the
same limit. This proves Theorem \VnPP.

\bigskip

\subsection{\bf Proof of Theorem \limitPP}%
The theorem follows from 
Lemma \UnNeglect, \togetherB\ and Theorem \VnPP, with the 
same arguments as at the end of proof of Theorem \LimitCD\ to show that
the limiting behavior of $\TT_n^-$ can be deduced from that of $\TT_n^+$
and that $\TT_n^-$ and $\TT_n^+$ are asympotically independent.

\bigskip


\def\prevs{\the\sectionnumber .\the\snumber }
\def\preveq{(\the\sectionnumber .\the\equanumber)}

\section{Application to \poorBold{$(g,F)$}-processes}
A $(g,F)$-process is both an extension and an abstraction of a FARIMA one.
Beyond being more general, the main motivation for introducing 
$(g,F)$-processes is that they are better suited for calculations than 
FARIMA ones, because it substitutes the asymptotics related to the 
gamma function involved in the analysis of the coefficients of FARIMA processes
by simpler relations involving regularly varying functions. Such a process is
given by a function $g$, analytic on $(-1,1)$, and a distribution function $F$,
hence its name. It is defined as follows. Let $(X_i)_{i\geq 1}$ be a sequence
of independent random variables, equidistributed according to $F$. We set
$X_i=0$ if $i$ is nonpositive. Let $B$ be the backward shift on sequences,
defined by $BX_i=X_{i-1}$. The analytic function $g$ has a series 
representation $g(x)=\sum_{i\geq 0} g_i x^i$. A $(g,F)$-process is a 
discrete time stochastic process $(S_n)_{n\geq 0}$ given by
$$
  S_n
  =g(B) X_n
  =\sum_{0\leq i<n} g_i X_{n-i}
  =\sum_{1\leq i\leq n} g_{n-i}X_i \, .
$$
Taking the function $g$ to be $1/(1-\Id)$ yields $(S_n)$ to be the random
walk with increments $(X_n)$. Taking $g$ to be a rational function defined
on the complex unit disk yields a causal ARIMA process. Finally, FARIMA 
processes are obtained by taking $g$ to be $(1-\Id)^{-d}$ times a rational 
function defined on the complex unit disk.

The usefulness of this class of processes in terms of calculation comes
from Karamata's theorem for power series (Bingham, Goldie and Teugels, 1989, 
Corollary 1.7.3), which allows us to relate the 
asymptotic behavior of the coefficients $(g_n)$ to that of the 
function $g(1-1/\Id)$. We first restate this result under the form given in
Lemma 5.1.1 in Barbe and McCormick (2008), which is tailored to our
analysis. For any positive real number $x$, we define the partial sum of the
coefficients,
$$
  \g{x}=\sum_{0\leq i<x} g_i
$$
and set $\g{0}=0$.

\Lemma{\label{BS}%
  Assume that $(g_n)$ is ultimately monotone. Let $\gamma$ be a real number
  greater than $-1$. The sequence $(g_n)$ is regularly varying of index
  $\gamma-1$ if and only if the function $g(1-1/\Id)$ is regularly
  varying of index $\gamma$ at infinity. In this case, as $n$ tends to 
  infinity and uniformly in $x$ in any compact subset of the positive 
  half-line,

  \smallskip

  \noindent (i) $\ds g_{\lfloor nx\rfloor}
  \sim {\ds x^{\gamma-1}\over\ds \Gamma(\gamma)}\,{\ds g(1-1/n)\over\ds n}$,

  \smallskip

  \noindent (ii) $\ds \g{nx}
  \sim {\ds x^\gamma\over\ds\Gamma(1+\gamma)}\, g(1-1/n)$.
}

\bigskip

We will need to assume in fact slightly more than the regular variation 
involved in Lemma \BS, namely that
$$
  (g_n)_{n\geq 0} \hbox{ is normalized regularly varying of index $\gamma-1$;}
$$
this means (see e.g.\ Bingham, Goldie and Teugels, 1989, Theorem 1.9.8) that
$$
  {g_{n+1}\over g_n} = 1+{\gamma-1\over n} \bigl( 1+o(1)\bigr)
  \eqno{\equa{gNRV}}
$$
as $n$ tends to infinity, interpreted as $g_{n+1}/g_n=1+o(1/n)$ when $\gamma$ 
is $1$. Equivalently, it means that the Karamata 
representation of the sequence is, for $n$ at least $1$,
$$
  g_n=Cn^{\gamma-1}\exp\Bigl(\sum_{1\leq i\leq n}\delta_i/i\Bigr)
$$
where $(\delta_i)$ is a sequence converging to $0$ (see Bingham, Goldie and
Teugels, 1989, \S 1.9).

Except when $\gamma=1$, assumption \gNRV\ implies that $(g_n)$ is ultimately
monotone, and therefore, Lemma \BS\ applies.

We define the centering for our $(g,F)$-process,
$$
  c_{n,k}=\g{k}\int_{F^\leftarrow(1/n)}^{F^\leftarrow(1-1/n)} x\d F(x) \, ,
  \qquad k\geq 1\, .
  \eqno{\equa{cnkDef}}
$$

Paralleling Theorem \LimitCD, we have the following.

\Theorem{\label{GammaGTOne}
  Assume that $\gamma$ is greater than $1$ and that \AssumptionStable\
  and \gNRV\ hold. Let $L$ be a L\'evy process satistfying \cfL. The
  distribution of the process 
  $$
    t\in [\,0,\infty)\mapsto 
    {n( S_{\lfloor nt\rfloor}-c_{n,\lfloor nt\rfloor})\over 
     \g{n}F^\leftarrow_*(1-1/n)}
  $$
  converges to that of the fractional L\'evy process
  $$
    t\in [\,0,\infty)\mapsto
    \int_0^t \gamma(t-s)^{\gamma-1} \d L(s)
  $$
  in ${\rm D}[\,0,\infty)$ equipped with the topology of the uniform 
  convergence on compactas. 
  The limiting process is continuous.
}

\bigskip

In contrast, when $\gamma$ is less than $1$, we have the following result. It 
does not require $(g_n)$ to be a regularly
varying sequence, but it assumes that
$$
  \sum_{i\geq 0} |g_i|^r<\infty \quad\hbox{ for some positive $r$.}
  \eqno{\equa{gSummability}}
$$

\Theorem{\label{GammaLTOne}
  Assume that $\gamma$ is less than $1$ and that \AssumptionStable\ 
  and \gSummability\
  hold. Let $\sum_{i\geq 1} \delta_{(t_i,x_i)}$ be a Poisson random measure
  with mean intensity $\d\lambda\otimes\d\nu$.
  The distribution of the point process
  $$
    \sum_{i\geq 1} 
    \delta_{\bigl(i/n\, ,\,
                  {\ss S_i-c_{n,i}\over\ss F_*^\leftarrow(1-1/n)}
            \bigr)}
  $$
  converges to that of
  $$
    \sum_{i\geq 1}\sum_{j\geq 0} \delta_{(t_i,g_jx_i)} \, .
  $$
}

\bigskip

If $(g_n)$ is regularly varying with index $\gamma-1$,
Theorems \GammaGTOne\ and \GammaLTOne\ leave open the case $\gamma=1$. This
case is more difficult though not impossible to settle. First, if all the
$g_i$'s are equal, the standard result on 
partial sum process applies. Second, if the $(g_i)$ take 
distinct values, Avram and Taqqu's (1992) results
suggest that one should consider convergence in the ${\rm M}_1$-topology. One
can reasonably conjecture that when $\gamma=1$, if $\liminf g_n>0$ then
the convergence stated in Theorem \GammaGTOne\ holds in 
the ${\rm M}_1$-topology, while if $\limn g_n=0$, then the conclusion of 
Theorem \GammaLTOne\ remains. 

\bigskip

In order to prove Theorems \GammaGTOne\ and \GammaLTOne, we will use a scaling
argument. The following lemma is stated for convenience. Define the
function
$$
  h(M)=
  \cases{ {\ds\alpha\over\ds\alpha-1}(p^{1/\alpha}-q^{1/\alpha})
            (1-M^{(1/\alpha)-1})
          & if $\alpha\not=1$,\cr\noalign{\vskip 3pt}
          (p-q)\log M & if $\alpha=1$.\cr}
$$

\Lemma{\label{CenteringAdjust}
  For any $M$ greater than $1$,
  $$
    c_{nM,k}-c_{n,k}
    \sim \g{k} h(M) {F_*^\leftarrow(1-1/n)\over n}
  $$
  uniformly in nonnegative integers $k$ and as $n$ tends to infinity.
}

\bigskip

\Proof We prove the lemma when $p$ does not vanish. The proof when $p$ vanishes
follows by analogous arguments. The proof of Karamata's theorem shows that,
as $n$ tends to infinity,
$$\displaylines{\qquad
  \int_{F^\leftarrow(1-1/n)}^{F^\leftarrow(1-1/nM)} x \d F(x)
  \hfill\cr\noalign{\vskip 3pt}\hfill
  {}\sim
  \cases{\noalign{\vskip -3pt} {\ds\alpha\over\ds\alpha-1} (1-M^{(1/\alpha)-1}) 
          {\ds F^\leftarrow(1-1/n)\over \ds n}   & if $\alpha\not=1$,\cr
          \noalign{\vskip 3pt}
          \log M  {\ds F^\leftarrow(1-1/n)\over\ds n} & if $\alpha=1$.\cr}
  \qquad\cr}
$$
We then use \FFStarQuantile\ to substitute $p^{1/\alpha}F_*^\leftarrow(1-1/n)$
for $F^\leftarrow(1-1/n)$ in this asymptotic equivalent.

Similarly, and even if $q$ vanishes, the same arguments shows that
$$\displaylines{\qquad
  \int_{F^\leftarrow(1/nM)}^{F^\leftarrow(1/n)} x\d F(x)
  \hfill\cr\noalign{\vskip 3pt}\hfill
  {}\sim
  \cases{\noalign{\vskip -3pt}
          -q^{1/\alpha}{\ds\alpha\over\ds\alpha-1}(1-M^{(1/\alpha)-1}) 
          {\ds F_*^\leftarrow(1-1/n)\over\ds n}
          & if $\alpha\not=1$, \cr
          \noalign{\vskip 3pt}
          -q\log M {\ds F_*^\leftarrow(1-1/n)\over\ds n} & if $\alpha=1$,\cr}
  \qquad\cr}
$$
provided that, when $q$ vanishes this asymptotic equivalence is understood as
meaning that the left hand side is $o\bigl(F_*^\leftarrow(1-1/n)/n\bigr)$
as $n$ tends to infinity. Given \cnkDef, this implies the lemma.\hfill\qed

\bigskip

\noindent{\bf Proof of Theorem \GammaLTOne.} To relate our $(g,F)$-processes
to the processes introduced in the \fixedref{first} section, referring to 
\HypKnKappa, set $\kappa_i=g_i$ if $i$ is nonnegative, and $\kappa_i=0$
otherwise. Then, considering \procXStep,
$$\eqalign{
  \oSS_n(t)
  &{}=\sum_{1\leq i\leq \lfloor nt\rfloor} \kappa_{\lfloor nt\rfloor-i}X_i \cr
  &{}=\sum_{0\leq i<\lfloor nt\rfloor} g_i X_{\lfloor nt\rfloor-i} \cr
  &{}=S_{\lfloor nt\rfloor}  \cr
  }
$$
for any $t$ nonnegative. Thus,
$$
  \osn(t)
  = \g{nt}\int_{F^\leftarrow(1/n)}^{F^\leftarrow(1-1/n)} x\d F(x)
  = c_{n,\lfloor nt\rfloor} \, .
$$

Let $M$ be a positive integer. Lemma \CenteringAdjust\ implies that
$$
  \max_{0\leq i\leq nM} {|c_{nM,i}-c_{n,i}|\over F_*^\leftarrow(1-1/n)}
  \leq c {\g{nM}\over n} \, .
$$
Since $\gamma$ is less than $1$, this upper bound tends to $0$ as $n$ 
tends to infinity. 

Set
$$
  R_n(t)
  ={S_{\lfloor nt\rfloor}-c_{n,\lfloor nt\rfloor}\over F_*^\leftarrow(1-1/n)}
  \, .
$$
The identity
$$
  R_n(i/n)
  = {F_*^\leftarrow(1-1/nM)\over F_*^\leftarrow(1-1/n)}R_{nM}(i/nM)
  + {c_{nM,i}-c_{n,i}\over F_*^\leftarrow(1-1/n)}
$$
and Lemma \CenteringAdjust\ then imply that, as $n$ tends to infinity,
$$
  R_n(i/n)=M^{1/\alpha}R_{nM}(i/nM)\bigl( 1+o(1)\bigr) +o(1)
  \eqno{\equa{PfGammaLTOneB}}
$$
where the $o(1)$-terms are uniform in $i$ between $0$ and $nM$.

Theorem \limitPP\ with $\kappa_i=g_i$ yields that the distribution of 
the random measure $N_n=\sum_{1\leq i\leq n} \delta_{(i/n\, ,\,R_n(i/n))}$
converges to that of
$$
  \sum_{j\geq 1}\sum_{i\geq 0} 
  \delta_{(V_j\, ,\,p^{1/\alpha} g_i W_j^{-1/\alpha})}
  +\delta_{(V_j'\, ,\, -q^{1/\alpha}g_i{W_j'}^{-1/\alpha})} \, .
  \eqno{\equa{PfGammaLTOneA}}
$$
Exercise 5.4, p.163, in Resnick (2007) shows that
$$\Pi_M=
  \sum_{j\geq 1} 
  \delta_{(V_jM\, ,\,p^{1/\alpha} W_j^{-1/\alpha}M^{1/\alpha})}
  +\delta_{(V_j'M\, ,\, -q^{1/\alpha}{W_j'}^{-1/\alpha}M^{1/\alpha})}
$$
is a Poisson random measure with mean intensity $(\One_{[0,M]}\d\lambda)\otimes
\d\nu$. Moreover, the random measure \PfGammaLTOneA\ is $\sum_{i\geq 0}
\int\delta_{(x,g_iy)}\d\Pi_1(x,y)$. Combined with \PfGammaLTOneB,
this implies that for any integer $M$, the distribution of the
random measure
$$\displaylines{\qquad
  \sum_{1\leq i\leq nM} \delta_{(i/n,R_n(i/n))}
  \hfill\cr\hfill
  = \sum_{1\leq i\leq nM}\delta_{\bigl(\ss(i/nM)M,M^{1/\alpha}R_{nM}(i/nM)
  (1+o(1))+o(1)\bigr)}
  \qquad\cr}
$$
converges to that of
$$
  \sum_{i\geq 0} \int\delta_{(xM,g_iyM^{1/\alpha})}\d\Pi_1(x,y)
  = \sum_{i\geq 0}\int\delta_{(x,g_iy)} \d\Pi_M(x,y) \, .
$$
The result follows since $M$ is arbitrary.\hfill\qed

\bigskip

\noindent{\bf Proof of Theorem \GammaGTOne.}
Define $\okn(t)=ng_{\lfloor nt\rfloor}/\g n$ for $t$ nonnegative, 
and $\okn(t)=0$ for $t$ negative. As in \HypKnLin, the function $k_n$ is
defined as a linear interpolation of $\okn$ between the points 
$(i/n)_{0\leq i\leq n}$. We have
$$
  {nS_{\lfloor nt\rfloor}\over \g n}
  = \sum_{1\leq i\leq nt} {ng_{\lfloor nt\rfloor -i}\over \g n} X_i
  = \sum_{1\leq i\leq n} \okn(t-i/n) X_i \, ,
$$
and
$$
  \osn(t)
  ={n\g{nt}\over\g{n}}\int_{F^\leftarrow(1/n)}^{F^\leftarrow(1-1/n)}x\d F(x)
  = {nc_{n,\lfloor nt\rfloor}\over\g{n}}
  \, .
$$
Theorem \GammaGTOne\ is then obtained by an application of 
Corollary \LimitCDCorollary.

Lemma \BS\
implies that the sequence $(\okn)$ converges to $k=\gamma\Id^{\gamma-1}$ 
locally uniformly in $(0,\infty)$. Since $\gamma$ is greater than $1$,
this convergence is in fact locally uniformly in $[\,0,\infty)$, showing 
that \HypKnCvUnif\ holds. We are left to check \HypKModulus\ and \HypKnModulus.

\bigskip

\noindent{\it Checking \HypKModulus.} We take $r=2$. 
Let $0\leq s\leq t\leq 1$ and set $\delta=t-s$. Define
$$\eqalign{
  \phi(\delta,s)
  &{}=\gamma^{-2}\int_0^1\bigl( k(t-u)-k(s-u)\bigr)^2\d u\cr
  &{}=\int_0^t (t-u)^{2(\gamma-1)}\d u +\int_0^s(s-u)^{2(\gamma-1)}\d u \cr
  &\phantom{{}=\int_0^t (t-u)^{2(\gamma-1)}\d u}
   {}- 2\int_0^s (t-u)^{\gamma-1}(s-u)^{\gamma-1} \d u \cr
  &{}={(s+\delta)^{2\gamma-1}+s^{2\gamma-1}\over 2\gamma-1}
    -2\int_0^s (\delta+v)^{\gamma-1} v^{\gamma-1} \d v \, . \cr
  }
$$
Since
$$\eqalign{
  {\partial\phi\over\partial s} (\delta,s)
  &{}=(s+\delta)^{2(\gamma-1)} + s^{2(\gamma-1)} 
    -2 (s+\delta)^{\gamma-1}s^{\gamma-1} \cr
  &{}=\bigl( (s+\delta)^{\gamma-1}-s^{\gamma-1}\bigr)^2 \, , \cr
  }
$$
the function $s\mapsto \phi(\delta,s)$ is nondecreasing. Thus, if $s$ is
between $0$ and $1-\delta$, this function is at most
$$
  \phi(\delta,1-\delta)={1+(1-\delta)^{2\gamma-1}\over 2\gamma-1}
  - 2\int_0^{1-\delta} (\delta+v)^{\gamma-1} v^{\gamma-1} \d v \, .
  \eqno{\equa{checkKA}}
$$
Note that
$$
  2\int_0^{1-\delta} v^{2\gamma-2} \d v
  = {2\over 2\gamma-1} -2\delta + O(\delta^2) \, , 
$$
and 
$$
 2(\gamma-1)\delta\int_0^{1-\delta} v^{2\gamma-3}\d v = \delta+O(\delta^2) \, .
$$
Since
$$
  {1+(1-\delta)^{2\gamma-1}\over 2\gamma-1}
  ={2\over 2\gamma-1} -\delta +O(\delta^2) \, ,
$$
equality \checkKA\ yields
$$
  \phi(\delta,1-\delta)= 2\int_0^{1-\delta} v^{2\gamma-2}+(\gamma-1)\delta 
  v^{2\gamma-3}-(\delta+v)^{\gamma-1}v^{\gamma-1} \d v + O(\delta^2)
$$
as $\delta$ tends to $0$. The change of variable $v=w\delta$ gives
$$\displaylines{\quad
  \phi(\delta,1-\delta)= 2\delta^{2\gamma-1}\int_0^{(1/\delta)-1} 
  w^{2\gamma-2}+(\gamma-1) w^{2\gamma-3}
  \hfill\cr\hfill
  -(1+w)^{\gamma-1}w^{\gamma-1} \d w + O(\delta^2)\, .
  \quad\cr}
$$
The above integral is convergent at $0$ since $\gamma$ exceeds $1$. 
The integrand is
$$
  w^{2\gamma-2} \Bigl( 1+{\gamma-1\over w}-\Bigl(1+{1\over w}\Bigr)^{\gamma-1}
  \Bigr)
  \sim -{(\gamma-1)(\gamma-2)\over 2} w^{2(\gamma -2)}
$$
as $w$ tends to infinity. Thus, we obtain
$$
  \phi(\delta,1-\delta) 
  =\cases{ O(\delta^2)                          & if $\gamma>3/2$,  \cr
           \noalign{\vskip 3pt}
           O\bigl(\delta^2\log (1/\delta)\bigr) & if $\gamma=3/2$,  \cr
           \noalign{\vskip 3pt}
           O(\delta^{2\gamma-1})                & if $1<\gamma<3/2$\cr}
$$
as $\delta$ tends to $0$. This confirms \HypKModulus\ since $\gamma>1$. 

\bigskip

\noindent{\it Checking \HypKnModulus.} 
To evaluate the modulus of continuity $\omega_{k_n}$ of $k_n$, we first
relate it to the increments of the sequence $(g_i)$. For this, for any 
nonnegative integer $r$, define
$$
  \Omega_n(r)
  =\max_{0\leq j\leq n-r} |g_{j+r}-g_j| \, .
$$
Our next result relates the discrete modulus of continuity of $k_n$ 
to $\Omega_n$ and is stated for convenience.

\Lemma{\label{omegaOmega}
  For any integer $n$ and any nonnegative $\delta$,
  $$
    \overline\omega_{k_n,r}(\delta)
    \leq \Bigl({n\over \g{n}}\Omega_n(\lfloor n\delta\rfloor)\Bigr)^r
    \, .
  $$
}

\Proof This is immediate since
$$
  k_n\Bigl({q-i\over n}\Bigr)-k_n\Bigl({p-i\over n}\Bigr)
  = {n\over \g{n}}(g_{q-i}-g_{p-i}) \, .
  \eqno{\qed}
$$

\smallskip

We can now check the assumption of Lemma \ModulusUnKnC\ 
when $\gamma$ exceeds $2$ and in some cases when $\gamma=2$.

\Lemma{\label{omegaGammaLarge}
  If $\gamma$ exceeds $2$, or, more generaly, if $(g_n/n)$ is asymptotically
  equivalent to a diverging nondecreasing sequence, then there 
  exists a constant $c$ such 
  that $\overline\omega_{k_n,r}(\delta)\leq c\delta^r$ 
  for all $n$ large enough.
}

\bigskip

\Proof For any positive $r$ at most $n$, the inequality
$$
  |g_{j+r}-g_j|\leq \sum_{j\leq i<j+r} |g_{i+1}-g_i|
$$
implies $\Omega_n(r)\leq r\Omega_n(1)$. Since $(g_n/n)$ is ultimately bounded
from below by a positive constant, assumption \gNRV\ implies that
$\Omega_n(1)$ is of order $g_n/n$. Thus, using Lemma \BS, there exists a 
constant $c_1$ such that
$$
  {n\over \g n} \Omega_n(\lfloor n\delta\rfloor)
  \leq {n\over \g n} \lfloor n\delta\rfloor c {g_n\over n} 
  \leq c_1\delta \, .
$$
The result follows from Lemma \omegaOmega.\hfill\qed

\bigskip

In order to evaluate $\overline\omega_{k_n,r}(\delta)$ when $g_n/n$ is 
not ultimately
bounded from below by a positive number, we need the following extension lemma.

\Lemma{\label{smoothRV}
  A sequence $(c_n)_{n\geq 1}$ is normalized regularly varying if and only
  if there exists a normalized regularly varying function $c^*$ such
  that $c_n=c^*(n)$ ultimately.
}

\bigskip

\Proof Clearly the condition is sufficient. To prove that it is necessary, 
consider the Karamata representation of $(c_n)_{n\geq 1}$,
$$
  c_n=n^\rho C\exp\Bigl( \sum_{1\leq i\leq n} {\delta_i\over i}\Bigr) \, ,
$$
where $(\delta_n)$ converges to $0$.
Define the function $\delta(u)=u\delta_{i+1}/(i+1)$ if $i\leq u<i+1$, 
$i\geq 0$. This ensures that
$$
  \sum_{1\leq i\leq n} {\delta_i\over i} = \int_0^n {\delta (u)\over u} \d u 
  \, .
$$
The function $\delta(\cdot)$ tends to $0$ at infinity. The function
$$
  c^*(x)=x^\rho C\exp\Bigl(\int_0^x {\delta(u)\over u}\d u\Bigr)
$$
is normalized regularly varying and, restricted to the integers, coincides
with the sequence $(c_n)_{n\geq 1}$.\hfill\qed

\bigskip

We now evaluate $\overline\omega_{k_n,r}(\delta)$ when $\gamma$ is at most $2$.

\Lemma{\label{omegaGammaSmall}
  Let $\gamma$ be in $(1,2\,]$. Let $\theta$ be any positive number less than
  $\gamma-1$. There exists a constant $c$ such that for any $n$ large enough
  and any nonnegative integer $r$ at most $n$,
  $$
    {n\over \g n} \Omega_n(r) \leq c\Bigl( {r\over n}\Bigr)^\theta \, .
  $$
  In particular, there exists a constant $c$ such that for any positive $n$ 
  large enough,
  $$
    \overline\omega_{k_n,r}(\delta)\leq c \delta^{\theta r} \, .
  $$
}

\Proof Lemma \omegaOmega\ shows that the second assertion follows from the
first one. So this lemma is really about the first assertion.

\parshape=4 0 pt \hsize 
            0 pt \hsize 
            0 pt \hsize
            0 pt 2.3 in
To prove it, recalling the definition of $\Omega_n(r)$, we consider
increments $g_{i+r}-g_i$ in five different ranges of $i$ and $r$ which are
illustrated below. We let $i_0$ and $r_0$ be two positive integers and $M$
be a (large) real number. The real number $\theta$ is positive and less
than $\gamma-1$.
\setbox1=\hbox to 4.1cm{\eightpoints\vbox to 4cm{%
  \kern 3.5cm%
  \includegraphics{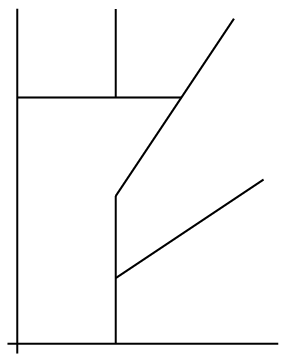}%
  \anote{.45}{-.05}{$i_1$}%
  \anote{-0.06}{1.05}{$r_1$}%
  \anote{-0.04}{1.4}{$r$}%
  \anote{1.1}{-.05}{$i$}%
  \anote{.95}{1.35}{$r=Mi$}%
  \anote{1.1}{.7}{$r=i/M$}%
  \anote{.1}{.5}{case 5}%
  \anote{.5}{1.25}{case 2}%
  \anote{.1}{1.25}{case 4}%
  \anote{.7}{.2}{case 3}%
  \anote{.6}{.7}{case 1}%
  \vfill%
  \vss%
}\hfill}%
\ht1=0pt\dp1=0pt\wd1=0pt
\setbox2=\hbox{\hskip 145pt\raise -95pt\box1\hss}
\wd2=0pt\ht2=0pt\dp2=0pt

\parshape=8 0 pt 2.3 in 
            0 pt 2.3 in 
            0 pt 2.3 in  
            0 pt 2.3 in 
            0 pt 2.3 in 
            0 pt 2.3 in 
            0 pt 2.3 in 
            0 pt \hsize
\noindent{\it Case 1.} \box2 We consider the range
$i\geq i_1$ and $i/M\leq r\leq Mi$. Since the sequence
$(g_i)$ is regularly varying of index $\gamma-1$ under \gNRV, 
the uniform convergence
theorem (Bingham, Goldie and Teugels, 1989, Theorem 1.2.1) implies that 
if $i_1$ is chosen large enough,\finetune{\break\vfill\eject\noindent}then for any $i$ at least $i_1$ and
any $\lambda$ in $[\, 1/M,M\,]$,
\finetune{\vskip -20pt}
$$\eqalignno{
  |g_{i+\lfloor i\lambda\rfloor}-g_i|
  &{}\leq 2 g_i\bigl( (\lambda+1)^{\gamma-1}-1\bigr)\cr
  &{}\leq 4 g_i(M+1)^{\gamma-1} \, .
  &\equa{omegaGammaSmallA}\cr
  }
$$
Thus, using Potter's bounds 
(Bingham, Goldie and Teugels, 1989, Theorem 1.5.6), since $\theta$ is less than
$\gamma-1$, provided that $i_1$ is large enough, \omegaGammaSmallA\ implies
$$
  {|g_{i+\lfloor i\lambda\rfloor}-g_i|\over g_n}
  \leq 4 (M+1)^{\gamma-1}\Bigl( {i\over n}\Bigr)^\theta \, .
$$
Writing $r$ for $\lfloor i\lambda\rfloor$, this yields, provided that $i_1$
is chosen large enough
$$\eqalign{
  {|g_{i+r}-g_i|\over g_n}
  &{}\leq 8(M+1)^{\gamma-1}{1\over\lambda^\theta}\Bigl({r\over n}\Bigr)^\theta
   \cr
  &{}\leq 8 (M+1)^{\theta+\gamma-1}\Bigl({r\over n}\Bigr)^\theta \, .\cr}
$$

\noindent{\it Case 2.} We consider the range $i\leq r/M$ and $r\geq r_1$ and
$i\geq i_1$. Since $\gamma$ exceeds $1$, assumption \gNRV\ implies that the
sequence $(g_n)$ is ultimately increasing. Thus, in the current range 
for $i$ and $r$, provided $i_1$, $r_1$ and $M$ are large enough,
$$
  0\leq g_{i+r}-g_i\leq g_{i+r}\leq 2 g_{\lfloor r ((1/M)+1)\rfloor} \, .
$$
Therefore, using Potter's bounds, provided that $r$ exceeds $r_1$ and $r_1$
is large enough,
$$\eqalign{
  0\leq {g_{i+r}-g_i\over g_n}
  &\leq 4 \Bigl( {1\over M}+1\Bigr)^{\gamma-1} {g_r\over g_n} \cr
  &\leq 8 \Bigl({r\over n}\Bigr)^\theta \, . \cr
  }
$$

\noindent{\it Case 3.} We consider the range $i\geq i_1$ and $i\geq Mr$.
Combining Lemma \smoothRV\ and Proposition 5.2.2 in Barbe and 
McCormick (2009), we see that
$$
  \Bigl| {g_{i(1+r/i)}-g_i\over {r\over i}g_i}-(\gamma-1)\Bigr|
  \leq {\gamma-1\over 2}
$$
provided $i_1$ and $M$ are large enough. Thus,
$$
  |g_{i+r}-g_i|\leq {3\over 2}(\gamma-1){r\over i} g_i \, .
$$
Then, Potter's bounds imply
$$\eqalignno{
  {|g_{i+r}-g_i|\over g_n}
  &{}\leq 4 (\gamma-1){r\over i} \Bigl({i\over n}\Bigr)^\theta \cr
  &{}\leq 4 (\gamma-1){r\over n} \Bigl({i\over n}\Bigr)^{\theta-1} \, . 
  &\equa{omegaGammaSmallB}\cr
  }
$$
Note that $\theta-1$ is nonpositive since $\gamma$ is at most $2$. 
Since $i\geq Mr$, the inequality $i/n\geq M r/n$ holds, and therefore, 
\omegaGammaSmallB\ yields
$$\eqalign{
  {|g_{i+r}-g_i|\over g_n}
  &\leq 4 (\gamma-1) {r\over n} \Bigl( {Mr\over n}\Bigr)^{\theta-1} \cr
  &= 4(\gamma-1) M^{\theta-1} \Bigl({r\over n}\Bigr)^\theta \, .\cr
  }
$$

\noindent{\it Case 4.} We consider the range $i\leq i_1$ and $r>r_1$. In this
range, povided $r_1$ is large enough, since $(g_n)$ tends to infinity,
$$
  |g_{i+r}-g_i|\leq 2 g_{i+r}\leq c g_r \, .
$$
Thus, using Potter's bounds, provided $r_1$ and $n$ are large enough,
$$
  {|g_{i+r}-g_i|\over g_n}
  \leq c {g_r\over g_n}
  \leq 2 c \Bigl({r\over n}\Bigr)^\theta \, .
$$

\noindent{\it Case 5.} We now cover a region where both $i$ and $r$ remain
bounded. In this region, there exists a constant $c$ such that
$$
  {|g_{i+r}-g_i|\over g_n} \leq {c\over g_n} \, .
  \eqno{\equa{omegaGammaSmallC}}
$$
Since $\theta$ is less than $\gamma-1$, it follows from Proposition 1.3.6 in
Bingham, Goldie and Teugels (1989) that for $n$ large enough, 
$g_n\geq n^\theta$. Thus, since $r$ is bounded, \omegaGammaSmallC\ implies
that for some constant $c$,
$$
  {|g_{i+r}-g_i|\over g_n}\leq c \Bigl({r\over n}\Bigr)^\theta \, .
$$

Combining all five cases, we proved that there exists a constant $c$ such that
for all $n$ large enough and all nonnegative integers $r$ at most $n$, the 
inequality $\Omega_n(r)/g_n \leq c(r/n)^\theta$ holds. Then, the first
assertion of the lemma follows from Lemma \BS.\hfill\qed

\bigskip

To conclude the proof of Theorem \GammaGTOne, define the fractional
L\'evy process
$$
   L^{(\gamma-1)}(t)=\int_0^t\gamma(t-u)^{\gamma-1}\d L(u) \, .
$$
Taking $r$ large enough, Lemmas \omegaGammaLarge\ and \omegaGammaSmall\ show
that \HypKnModulus\ holds. Hence, Theorem \LimitCD\ implies that the
distribution of the process
$$
  t\in [\,0,1\,]\mapsto 
  {n(S_{\lfloor nt\rfloor}-c_{n,\lfloor nt\rfloor})\over 
   \g{n}F_*^\leftarrow(1-1/n)}
$$
converges to that of $t\in [\,0,1\,]\mapsto L^{(\gamma-1)}(t)$ 
in ${\rm D}[\,0,1\,]$ equipped with the supremum norm. To extend
this con\-ver\-gence in $\D[\,0,1\,]$ to one in $\D[\,0,\infty)$, let $M$ be a
positive integer. For any $t$ in $[\,0,M\,]$,
$$\displaylines{\quad
  {n (S_{\lfloor nt\rfloor}-c_{n,\lfloor nt\rfloor})
    \over \g{n}F_*^\leftarrow(1-1/n)}
  \hfill\cr\hfill
  {}= {nM (S_{\lfloor nMt/M\rfloor}-c_{nM,\lfloor nMt/M\rfloor})\over \g{nM}
   F_*^\leftarrow(1-1/nM)} M^{\gamma+(1/\alpha)-1} \bigl( 1+o(1)\bigr)
  \qquad\cr\hfill
  {}+{n\over \g{n}F_*^\leftarrow(1-1/n)} 
  (c_{nM,\lfloor nt/\rfloor}-c_{n,\lfloor nt\rfloor}) \, .
  \quad\equa{GammaGTOneA}\cr
  }
$$
Hence, given \cnkDef\ and Lemma \CenteringAdjust, the distribution of the 
right hand side of \GammaGTOneA, as a process index by $t$ converges to that of 
$$
  L^{(\gamma-1)}(t/M)M^{\gamma+(1/\alpha)-1} +t^\gamma h(M) \, . 
  \eqno{\equa{GammaGTOneB}}
$$
The remainder of the proof shows that this limiting process has the same 
distribution as $L^{(\gamma-1)}$.

Referring to the process $L_0$ introduced in \DefLZero, we see that the 
fractional integral
$$
  L_0^{(\gamma-1)}(t)=\int_0^t\gamma (t-x)^{\gamma-1} \d L_0(x)
$$
satisfies, as function of $t$
$$\eqalignno{
  L_0^{(\gamma-1)}(\lambda t)
  &{}= \int_0^{\lambda t}\gamma(\lambda t-y)^{\gamma-1}\d L_0(y) \cr
  &{}\eqd\lambda^{\gamma-1+(1/\alpha)} L_0^{(\gamma-1)}(t) \, ;
  &\equa{LZeroFracSelfSimilar}\cr}
$$
in other words, the process $L_0^{(\gamma-1)}(\lambda\,\cdot\,)$ has the same 
distribution as $\lambda^{\gamma-1+(1/\alpha)}L_0^{(\gamma-1)}$.
From the definition of $L_0$ in \DefLZero, one can check that
\hfuzz=1pt
$$\displaylines{
  L^{(\gamma-1)}(t) = L_0^{(\gamma-1)}(t)
  \hfill\cr\noalign{\vskip 5pt}\hfill
  {}+\cases{{\ds\alpha\over\ds\alpha-1}(p^{1/\alpha}-q^{1/\alpha})t^\gamma
          &if $\alpha\not=1$,\cr
          (p-q)\ds\int_0^t\gamma(t-x)^{\gamma-1}(1+\log x)\d x
          &if $\alpha=1$.\cr}
  \quad
  \equa{LLZeroFrac}\cr}
$$
\hfuzz=0pt
Consider now the case where $\alpha\not=1$. Equalities \LZeroFracSelfSimilar\ 
and \LLZeroFrac\ yield
$$\displaylines{\qquad
 L^{(\gamma-1)}(t/M)
 \hfill\cr\noalign{\vskip 2pt}\hfill
 \eqalign{
  {}\eqd{}&(1/M)^{\gamma-1+1/\alpha}L_0^{(\gamma-1)}(t)
        + {\alpha\over\alpha-1} (p^{1/\alpha}-q^{1/\alpha})(t/M)^\gamma \cr
  {}={}&(1/M)^{\gamma-1+1/\alpha} 
        \bigl( L^{(\gamma-1)}(t) -{\alpha\over\alpha-1}
        (p^{1/\alpha}-q^{1/\alpha})t^\gamma\bigr)\cr
  &\hskip 120pt
    {}+{\alpha\over\alpha-1} (p^{1/\alpha}-q^{1/\alpha})(t/M)^\gamma \cr
  {}={}&(1/M)^{\gamma-1+1/\alpha} \bigl( L^{(\gamma-1)}(t)-t^\gamma h(M)\bigr) 
   \, . \cr
  }
  \cr}
$$
This proves that the process in \GammaGTOneB\ has the same distribution as
$L^{(\gamma-1)}$.

When $\alpha=1$, the same arguments show that the process
$L^{(\gamma-1)}(t/M)$ indexed by $t$ has the same distribution as
$$\displaylines{\qquad
  (1/M)^\gamma L_0^{(\gamma-1)}(t)
        + (p-q)\int_0^{t/M}\gamma\Bigl({t\over M}-x\Bigr)^{\gamma-1} 
        (1+\log x)\d x 
  \hfill\cr\hfill
  \eqalign{
  {}={}&(1/M)^\gamma\Bigl( L^{(\gamma-1)}(t) -
        (p-q)\int_0^t \gamma(t-x)^{\gamma-1}(1+\log x)\d x\Bigr) \cr
       &\hskip 30pt{} + {p-q\over M^\gamma}\int_0^t \gamma (t-y)^{\gamma-1} 
        (1+\log y-\log M)\d y \cr
  {}={}&{1\over M^\gamma} L^{(\gamma-1)}(t)-\Bigl({t\over M}\Bigr)^\gamma h(M) 
        \, . \cr
  }
  \cr}
$$
Again, this implies that the process in \GammaGTOneB\ has the same distribution
as $L^{(\gamma-1)}$ and this proves Theorem \GammaGTOne.\hfill\qed

\bigskip


\section{Application to FARIMA processes}
In this short section we consider a traditional FARIMA model. It is defined
by two polynomials, $\Theta$ and $\Phi$, with $\Theta(1)\not=0$, the roots 
of $\Phi$ being outside the complex unit disk, and 
a positive real number $\gamma$, setting
$$
  (1-B)^\gamma\Phi(B)S_n=\Theta(B)X_n \, .
$$
Thus, this is a $(g,F)$-process with $g=(1-\Id)^{-\gamma}\Theta/\Phi$. In 
order to apply the results of the previous section, it suffices to show that
the coefficients $(g_n)$ form a normalized regularly varying sequence of
nonvanishing index. This follows from the following slight refinement 
of Property 2 in Akonom and Gouri\'eroux (1987).

\Lemma{\label{FarimaNRV}
  If $\gamma$ is greater than $1$, then the coefficients $(g_n)$ 
  associated to the function 
  $g=(1-\Id)^{-\gamma}\Theta/\Phi$ form a normalized regularly varying sequence
  of index $\gamma-1$. In particular $g_n\sim n^{\gamma-1}\Theta(1)/\Phi(1)$
  as $n$ tends to infinity.
}

\bigskip

As a consequence of the second assertion of Lemma \FarimaNRV, the sign of
$g_n$ ultimately coincides with that of $\Theta(1)/\Phi(1)$. 

Equipped with this result whose proof is deferred to the end of this section,
Theorems \GammaGTOne\ and \GammaLTOne\ readily apply. 

In particular, we see that if $\gamma$ is greater than $1$, the process
$S_{\lfloor n\cdot\rfloor}$ properly centered and rescaled converges to a
fractional L\'evy stable process. This result parallels Akonom and 
Gouri\'eroux's (1987) who prove that when $\gamma>1/2$ and the innovations
have a finite moment of order $\max\bigl( 2,1/(\gamma-1/2)\bigr)$, then 
$S_{\lfloor n\,\cdot\,\rfloor}$ properly normalized converges to a fractional
Brownian motion. In contrast, when the innovations are in the domain of
attraction of a nonGaussian stable distribution, Theorem \GammaLTOne\ 
shows that it is not possible to make the distributions of the rescaled 
processes $S_{\lfloor n\cdot\rfloor}$ converge when $\gamma$ is less than $1$.

In a slightly different spirit, Theorem \GammaLTOne\ and Lemma \FarimaNRV\
imply that when $\gamma$
is less than $1$, and $\Theta\Phi(1)$ is positive say, the distribution of 
$\max_{0\leq i\leq n}(S_i-s_i)/F_*^\leftarrow(1-1/n)$ converges to 
that of $p^{1/\alpha}\max_{i\geq 0} g_i W_1^{-1/\alpha}$. In particular,
if the innovations are symmetric,
$$
  \limn \Prob\{\, \max_{0\leq i\leq n} S_i<xF_*^\leftarrow(1-1/n)\,\}
  = \exp\Bigl(-x^{-\alpha}{\max_{i\geq 0} g_i^\alpha\over 2}\Bigr) \, ,
$$
a result of interest for extreme value analysis. It is clear that, using
the result of this paper, this result
can be considerably generalized by studying the extremal process associated
to the $(g,F)$-process following the lines of Resnick (1987, chapter 4).

In contrast, when $\gamma$ is greater than $1$, Theorem \GammaGTOne\ and Lemma
\BS.ii show that, assuming still that the innovations are symmetric,
$$\displaylines{\qquad
  \limn \Prob\Bigl\{\, \max_{1\leq i\leq n} S_i<x 
  {n^{\gamma-1}\over\Gamma(1+\gamma)} {\Theta\over\Phi}(1) 
  F_*^\leftarrow(1-1/n)\,\Bigr\}
  \hfill\cr\hfill
  {}= \Prob\Bigl\{\, \sup_{0\leq t\leq 1} \int_0^t\gamma(t-s)^{\gamma-1}
  \d L(s) < x\,\Bigr\} \, .
  \qquad\cr}
$$
Hence, while in the $\gamma$ less than $1$ case, the partial maxima of the
process grows like $F_*^\leftarrow(1-1/n)$, it grows at a rate $n^{\gamma-1}$ 
times faster in the $\gamma$ greater than $1$ case.

In the special case where $\Phi$ and $\Theta$ are the constant polynomials
equal to $1$, we rewrite the process as
$$
  (1-B)^{\gamma-1}(1-B)S_n=X_n \, .
$$
This corresponds to a random walk being differentiated $\gamma-1$ times. 
Theorem \GammaGTOne\ and \GammaLTOne\ show that the usual convergence result
of the random walk to a L\'evy stable process is, as one would expect, sharp:
any slight differentiation ($\gamma<1$) prevents convergence of the process
in ${\rm D}[\,0,1\,]$, while any slight integration ($\gamma>1$) ensures 
convergence in uniform norm.

\bigskip

\noindent{\bf Proof of Lemma \FarimaNRV.} We will prove a more general result
which is of independent interest.

\Lemma{\label{FARIMAExpansion}
  Let $m$ be a positive integer and let $d$ be a real number greater than $1$. 
  Let $A$ be an analytic function 
  on $[\,-1,1\,]$ having expansion $A(x)=\sum_{i\geq 0} A_i x^i$ and such that
  $$
    \sum_{k\geq 0} k^m |A_k|<\infty \, ,
    \eqno{\equa{HypFARIMAExpansionA}}
  $$
  and, as $n$ tends to infinity,
  $$
    \sum_{k\geq n} |A_k|=o(n^{d-1-m}) \, .
    \eqno{\equa{HypFARIMAExpansionB}}
  $$
  Then the function $g=(1-\Id)^{-d} A$ has an 
  expansion $g(x)=\sum_{i\geq 0} g_i x^i$ and there exists a polynomial $P_m$
  of degree $m$ such that, as $n$ tends to infinity,
  $$
    g_n=n^{d-1}P_m(1/n) +o(n^{d-1-m}) \, ;
  $$
  in particular, $P_1(x)=A(1)+(d-1)\bigl( (d/2)A(1)-A'(1)\bigr)x$.
}

\bigskip

\noindent Note that under \HypFARIMAExpansionA, 
$$
  \sum_{k\geq n} |A_k|
  \leq \sum_{k\geq n} \Bigl({k\over n}\Bigr)^m |A_k|
  =o(n^{-m})
  \eqno{\equa{FARIMAExpansionAa}}
$$
as $n$ tends to infinity. Thus, \HypFARIMAExpansionB\ is meaningful only when
$d$ is less than $1$.

\bigskip

\Proof Our proof is a refinement of that of Property 2 of Akonom and Gourieroux
(1987) who give only an asymptotic equivalent for $g_n$ and an upper bound
for the error term. 

We will use the symbol $\sim$ between a function and a possibly divergent 
series to indicate that the series is an asymptotic expansion of the function; 
for instance, we write $f(x)\sim \sum f_n x^n$ if for every integer $m$,
$$
  f(x)=\sum_{0\leq n\leq m} f_n x^n+o(x^m)
$$
as $x$ tends to $0$; the same symbol will be used with an analogous meaning
when dealing with an asymptotic expansion at infinity instead of the origin 
---~we refer to Olver (1974, \S 1.7) for the usage of this symbol.

It is convenient to structure the proof into three steps.

{\noindent\it Step 1.} In this preliminary step, we recall some known result on
the expansion for the ratio of two gamma functions. Starting with the identity
$$
  (1-x)^{-d} = \sum_{k\geq 0} {\Gamma(k+d)\over \Gamma(k+1)\Gamma(d)} x^k \, ,
$$
we also have
$$
  (1-x)^{-d} A(x)
  = \sum_{k\geq 0} {\Gamma(k+d)\over \Gamma(k+1)\Gamma(d)}x^k 
  \sum_{k\geq 0}  A_k x^k 
$$
and
$$
  g_n=\sum_{0\leq k\leq n} {\Gamma(k+d)\over \Gamma(k+1)\Gamma(d)} A_{n-k} \, .
$$
The function $q(t)=(e^t-1)^{-d}$ admits an expansion at $0$,
$$
  q(t)\sim\sum_{j\geq 0} (-1)^j q_j t^{j-d} \, ,
$$
with $q_0=1$ and $q_1=d/2$. Setting $Q_0=1$,
$$
  Q_j=(d-1)(d-2)\ldots (d-j) q_j \, , \qquad j\geq 1\, ,
$$
and defining
$$
  h_m(x)=\sum_{0\leq j\leq m} Q_j x^{d-1-j} \, ,
$$
formula 5.02 in Olver (1974, \S 4.5) provides the error bound
$$
  \Bigl|{\Gamma(k+d)\over \Gamma(k+1)}-h_m(k)\Bigr|=o(k^{d-1-m})
  \eqno{\equa{FARIMAExpansionA}}
$$
as $k$ tends to infinity. In particular, $h_m$ is the $m+1$ term expansion
of $\Gamma(k+d)/\Gamma(k+1)$ on the asymptotic scale $(\Id^{d-1-i})_{i\geq 0}$.
We see that $Q_1=d(d-1)/2$. 

\noindent{\it Step 2.} In this step 
we will prove that, as $n$ tends to infinity,
$$
  g_n={1\over \Gamma(d)} \sum_{0\leq k\leq n} h_m(k) A_{n-k} +o(n^{d-1-m}) \, .
  \eqno{\equa{FARIMAExpansionD}}
$$

Given the previous step, we have
$$\displaylines{\qquad
  \Gamma(d)\Bigl| g_n-{1\over \Gamma(d)}\sum_{0\leq k\leq n}h_m(k)A_{n-k}\Bigr|
  \hfill\cr\hfill
  \leq \sum_{0\leq k\leq n} \Bigl|{\Gamma(k+d)\over \Gamma(k+1)}-h_m(k)\Bigr|
  |A_{n-k}| \, .
  \qquad\equa{FARIMAExpansionB}\cr
  }
$$
The bound \FARIMAExpansionA\ ensures that for $k$ larger than some $k_0$,
\FARIMAExpansionB\ is at most
$$
  c\sum_{0\leq k\leq k_0} |A_{n-k}|
  +\epsilon\sum_{k_0\leq k\leq n} k^{d-1-m} |A_{n-k}| \, .
  \eqno{\equa{FARIMAExpansionC}}
$$
If $d-1-m$ is nonnegative, this is at most
$$\displaylines{\qquad
  c k_0\max_{n-k_0\leq i\leq n} |A_i| 
  +\epsilon n^{d-1-m}\sum_{0\leq k\leq n-k_0} \Bigl( 1-{k\over n}\Bigr)^{d-1-m}
  |A_k| 
  \hfill\cr\hfill
  {}\leq c k_0 \max_{n-k_0\leq i\leq n} |A_i|+\epsilon n^{d-1-m}\sum_{k\geq 0}
  |A_k| \, . \qquad\cr
  }
$$
Thus, under assumption \HypFARIMAExpansionA, this is at 
most $2\epsilon n^{d-1-m}$. If $d-1-m$ is negative, \FARIMAExpansionC\ is at
most
$$\displaylines{\quad
    c k_0\max_{n-k_0\leq i\leq n} |A_i|
    +\epsilon\sum_{k_0\leq k\leq n/2} |A_{n-k}| 
    +\epsilon \sum_{n/2\leq k\leq n} \Bigl({n\over 2}\Bigr)^{d-1-m} |A_{n-k}| 
  \hfill\cr\hfill
    {}\leq c k_0\max_{n-k_0\leq i\leq n} |A_i|+\epsilon \sum_{k\geq n/2}|A_k|
    +\epsilon {2^{m+1-d}\over n^{m+1-d}}\sum_{k\geq 0} |A_k| \, . \cr
  }
$$
Using \HypFARIMAExpansionB, this upper bound is at most $c\epsilon
\sum_{k\geq 0} |A_k|/n^{m+1-d}$. Since $\epsilon$ is arbitrary, we see that
in both cases, \FARIMAExpansionC\ is $o(n^{d-1-m})$. In other words,
\FARIMAExpansionD\ holds.

\noindent{\it Step 3.} We now conclude the proof.
For any positive $\epsilon$ we also have, using \FARIMAExpansionAa,
$$\eqalign{
  \Bigl|\sum_{0\leq k\leq n(1-\epsilon)} h_m(k)A_{n-k}\Bigr|
  &{}\leq O(1) h_m(n)\sum_{k\geq \epsilon n} |A_k| \cr
  &{}=O(n^{d-1}) \sum_{k\geq \epsilon n} |A_k| \cr
  &{}=o(n^{d-1-m}) \, .\cr
  }
$$
By Taylor's formula, there exists $\theta_{k,n}$ between $n-k$ and $n$ such
that
$$
  h_m(n-k)=\sum_{0\leq p\leq m-1}{(-1)^p\over p!} k^p h_m^{(p)}(n)
  + {(-1)^m\over m!} k^m h_m^{(m)}(\theta_{n,k}) \, .
$$
We then have
$$\displaylines{
  \sum_{0\leq k\leq n\epsilon} h_m(n-k)A_k
    =\sum_{0\leq p\leq m-1} {(-1)^p\over p!} h_m^{(p)}(n)
    \sum_{0\leq k\leq n\epsilon} k^p A_k 
  \qquad\cr\hfill
    {}+{(-1)^m\over m!} h_m^{(m)}(n) \sum_{0\leq k\leq n\epsilon} 
    {h_m^{(m)}(\theta_{k,n})\over h_m^{(m)}(n)} k^m A_k \, ,
  \quad\equa{FARIMAExpansionE}
\cr}
$$
where, in this equality, it is agreed that if $h_m^{(m)}$ vanishes 
identically, then $h_m^{(m)}(\theta_{k,n})/h_m^{(m)}(n)$ is to be read as $1$.

Consider an arbitrary positive real number $\delta$. Since $h_m^{(m)}$ is
regularly varying, the uniform convergence theorem (Bingham, Goldie and 
Teugels, 1989, Theorem 1.2.1) implies that for $\epsilon$ small enough,
$$
  \max_{0\leq k\leq n\epsilon} 
  \Bigl| {h_m^{(m)}(\theta_{k,n})\over h_m^{(m)}(n)}-1 \Bigr| 
  \leq\delta \, .
$$
Thus, \FARIMAExpansionE\ yields
$$\displaylines{\quad
    \Bigl| \sum_{0\leq k\leq n\epsilon} h_m(n-k)A_k -\sum_{0\leq p\leq m} 
    {(-1)^p\over p!} h_m^{(p)}(n)\sum_{k\geq 0} k^p A_k\Bigr|
  \hfill\cr\hfill
    {}\leq \sum_{0\leq p\leq m} {|h_m^{(p)}(n)|\over p!} 
    \sum_{k\geq\epsilon n} k^p|A_k|
    +{|h_m^{(m)}(n)|\over m!} \delta\sum_{0\leq k\leq n\epsilon} k^m |A_k| \, .
  \quad\equa{FARIMAExpansionF}\cr
  }
$$
Using assumptions \HypFARIMAExpansionA, \HypFARIMAExpansionB, and 
since $\delta$ is arbitrary, we see
that \FARIMAExpansionF\ implies
$$
  \sum_{0\leq k\leq n\epsilon} h_m(n-k) A_k 
  = \sum_{0\leq p\leq m} {(-1)^p\over p!} h_m^{(p)}(n)\sum_{k\geq 0} k^p A_k
  + o\bigl( h^{(m)}_m(n)\bigr) \, .
  \eqno{\equa{FARIMAExpansionG}}
$$
Combined with \FARIMAExpansionD, this yields the result. Referring to 
\FARIMAExpansionG, we see that
$$
  \sum_{p\geq 0} {(-1)^p\over p!}k^p h_m^{(p)}
  = e^{-k\D} h_m(n)
$$
so that $P_m$ coincides with $\Id^{1-d}\Gamma(d)^{-1}\bigl(A(e^{-\D})h_m\bigr)
(1/\Id)$ modulo the ideal generated by $\Id^{m+1}$.\hfill\qed

\bigskip


\noindent {\bf References.}

\medskip

{\leftskip=\parindent \parindent=-\parindent
 \par

J.\ Akonom, Chr.\ Gouri\'eroux (1987). A functional central limit theorem
for fractional processes, discussion paper 8801, CEPREMAP, Paris.

F.\ Avram, M.S.\ Taqqu (1992). Weak convergence of sums of moving averages
in the $\alpha$-stable domain of attraction, {\sl Ann.\ Probab.}, 20, 483--503.


Ph.\ Barbe, W.P.\ McCormick (2008). Veraverbeke's theorem at large --
On the maximum of some processes with negative drift and heavy tail
innovations, {\tt arxiv:0802.3638}.

Ph.\ Barbe, W.P.\ McCormick (2009). {\sl Asymptotic Expansions for Infinite
Weighted Convolutions of Heavy Tail Distributions and Applications}, {\sl
Mem.\ Amer.\ Math.\ Soc.}, 922.


C.\ Bender, T.\ Marquardt (2008). Stochastic calculus for convoluted L\'evy
processes, {\sl Bernoulli}, 14, 499--518.

J.\ Bertoin (1996). {\sl L\'evy Processes}, Cambridge.

P.\ Billingsley (1968). {\sl Convergence of Probability Measures}, Wiley.

N.H.\ Bingham, C.M.\ Goldie, J.L.\ Teugels (1989). {\sl Regular Variation},
2nd ed. Cambridge University Press.

A.\ Chakrabarty, G.\ Samorodnitsky (2010). Understanding heavy tails in a
bounded world or, is a truncated heavy tail heavy or not?, 
{\tt arxiv:1001.3218v1}.

Y.S.\ Chow, H.\ Teicher (1988). {\sl Probability Theory. Independence, 
Interchangeability, Martingales}, Springer.



S.\ Cs\"org\H o, L.\ Horv\`ath, D.\ Mason (1986). What portion of the sample
makes a partial sum asymptotically stable or normal?, {\sl Probab.\ Theor.\
Rel.\ Fields}, 72, 1--16.

L.\ Devroye (1982). Upper and lower class sequences for minimal uniform 
spacings, {\sl Z.\ Wahrsch.\ verw.\ Geb.}, 61, 237--254.


H.\ Finkelstein (1971). The law of the iterated logarithm for empirical
distributions, {\sl Ann.\ Math.\ Statist.}, 42, 607--615.






K.\ Ito (1969). {\sl Stochastic Processes}, Lecture Notes Series, 16, Aarhus
Universitet.

S.\ Kwapie\'n, W.\ Woyczy\'nski (1992). {\sl Random Series and Stochastic
Integrals: Single and Multiple}, Birkhauser.

R.\ LePage, M. Woodroofe, J.\ Zinn (1981). Convergence to a stable distribution
via order statistics, {\sl Ann.\ Probab.}, 9, 624--632.


M.\ Lo\`eve (1960). {\sl Probability Theory}, 2nd ed., Van Nostrand.

M.\ Magdziarz (2008). Fractional Ornstein-Uhlenbeck processes, Joseph effect
in models with infinite variance, {\sl Physica A}, 387, 123--133.

T.\ Marquardt (2006). Fractional L\'evy processes with an application to long
memory moving average processes, {\sl Bernoulli}, 12, 1099--1126.

F.W.J.\ Olver (1974). {\sl Asymptotic Analysis and Special Functions}, 
Academic Press.


Yu.V.\ Prokhorov (1956). Convergence of random processes and limit theorems
in probability theory, {\sl Theor.\ Probab.\ Appl.}, 1, 157--214.

S.I.\ Resnick (1987). {\sl Extreme Values, Regular Variation and Point 
Processes}, Springer.

S.I.\ Resnick (2007). {\sl Heavy-Tail Phenomena, Probabilistic and Statistical
Modeling}, Springer.

G.\ Samorodnitsky, M.S.\ Taqqu (1994). {\sl Stable Non-Gaussian Random 
Processes}, Chapman \& Hall.

G.R.\ Shorack, J.A.\ Wellner (1986). {\sl Empirical Processes with Applications
to Statistics}, Wiley.

A.V.\ Skorohod (1956). Limit theorems for stochastic processes, {\sl Theor.\
Probab.\ Appl.}, 1, 261--290.

D.W.\ Stroock (1994). {\sl Probability Theory, an Analytic View} (revised
ed.), Cambridge.


}


\bigskip

\setbox1=\vbox{\halign{#\hfil&\hskip 40pt #\hfill\cr
  Ph.\ Barbe            & W.P.\ McCormick\cr
  90 rue de Vaugirard   & Dept.\ of Statistics \cr
  75006 PARIS           & University of Georgia \cr
  FRANCE                & Athens, GA 30602 \cr
  philippe.barbe@math.cnrs.fr                        & USA \cr
                        & bill@stat.uga.edu \cr}}%
\box1%

\bye

\vfill\eject

\pageno=1
\def\DATE{{November 10, 2009 (corrected March 26, 2010)}}

\snumber=0
\equanumber=0
\def\prevs{A.\the\snumber }
\def\preveq{({\rm A}.\the\equanumber)}

\centerline{\bf L\'evy measure and characteristic function of \poorBold{$L(1)$}}

\bigskip

Let us use the notation, for $x$ positive
$$
  \nu^+(x,\infty)=x^{-\alpha}
  \quad\hbox{ and }\quad
  \nu^-(-\infty,-x)=x^{-\alpha} \, ,
$$
and $\nu=p\nu^++q\nu^-$. So,
$$
  {\d\nu\over\d\lambda}(x)
  =p\alpha x^{-\alpha-1}\One_{(0,\infty)}(x)x^{-\alpha-1}
  + q\alpha |x|^{-\alpha-1}\One_{(-\infty,0)}(x) \, .
$$
So, here the $\nu$ is as in Resnick and not as in our introduction; the one
in our introduction should be changed.

Lemma \LPlusLMeasure\ (which I think is correct, but, as usual, check it) 
asserts that
$$\displaylines{\quad
    \E \exp\bigl(itL^+(1)\bigr)
    \hfill\cr\hfill
    {}=\exp\Bigl(\int_1^\infty (e^{itx}-1)\d\nu^+(x) 
    +\int_0^1(e^{itx}-1-itx)\d\nu^+(x)\Bigr) \, .
    \quad\cr}
$$ 
(note that there is a typo in that lemma in the second integral)

Similarly, with now $L^-$,
$$\displaylines{\quad
    \E \exp\bigl(itL^-(1)\bigr)
    \hfill\cr\hfill
    {}=\exp\Bigl(\int_{-\infty}^{-1} (e^{itx}-1)\d\nu^-(x) 
    +\int_{-1}^0(e^{itx}-1-itx)\d\nu^-(x)\Bigr) \, .
    \quad\cr}
$$
(To see that this is the correct expression, think of $\TT_n^-$ as 
$\TT_n^+[\M_1F]$ (I hope the notation is clear) which makes the book-keeping
easier).

Now, identity \SnVn\ shows that the limiting L\'evy process to consider
has the representation
$$
  L=p^{1/\alpha}L^++q^{1/\alpha}L^- \, .
$$

The characteristic function of $L(1)$ is then
$$
  \E \exp\bigl(itL(1)\bigr) 
  =\E\exp\bigl( itp^{1/\alpha}L^+(1)+itq^{1/\alpha}L^-(1)\bigr)
$$
Since $L^+$ and $L^-$ are independent, its logarithm is
$$\displaylines{
  \int_1^\infty (e^{itp^{1/\alpha}x}-1)\d\nu^+(x) 
    +\int_0^1(e^{itp^{1/\alpha}x}-1-itp^{1/\alpha}x)\d\nu^+(x)
  \hfill\cr\hfill
    +\int_{-\infty}^{-1} (e^{-itq^{1/\alpha}x}-1)\d\nu^-(x) 
    +\int_{-1}^0(e^{-itq^{1/\alpha}x}-1+itq^{1/\alpha}x)\d\nu^-(x)\, .
  \cr\hfill
  \equa{AAOne}\cr}
$$
Now make the change of variable $y=p^{1/\alpha}x$ in the integral with respect
to $\nu^+$ and the change of variable $y=q^{1/\alpha}x$ in the integral with
respect to $\nu^-$. We have
$$
  \d\nu^+(x)
  =\alpha x^{-\alpha-1} \d x
  = p \alpha y^{-\alpha-1}\d y
  = p\d\nu^+(y) \, ,
$$
and therefore \AAOne\ becomes
$$\displaylines{\qquad
  \int_{p^{1/\alpha}}^\infty (e^{ity}-1)p\d\nu^+(y) 
    +\int_0^{p^{1/\alpha}}(e^{ity}-1-ity)p\d\nu^+(y)
  \hfill\cr\hfill
    +\int_{-\infty}^{-q^{1/\alpha}} (e^{-ity}-1)q\d\nu^-(y) 
    +\int_{-q^{1/\alpha}}^0(e^{-ity}-1+ity)q\d\nu^-(y)\, .
  \qquad\cr}
$$
This is equal to
$$
  \int_{(-\infty,-q^{1/\alpha})\cup (p^{1/\alpha},^\infty)} 
   (e^{ity}-1)\d\nu(y) 
    +\int_{(-q^{1/\alpha},p^{1/\alpha})}(e^{ity}-1-ity)\d\nu(y)
$$
So, setting, with (as Resnick),
$$
  {\d\nu\over\d\lambda}(x)
  = p\alpha x^{-\alpha-1}\One_{(0,\infty)}(x) + q\alpha|x|\One_{(-\infty,0)}(x)
  \, ,
$$
we have
$$\displaylines{
  \E\exp\bigl(itL(1)\bigr)
  =\exp\Bigl(\int_{(-\infty,-q^{1/\alpha})\cup (p^{1/\alpha},^\infty)} 
   (e^{ity}-1)\d\nu(y) 
  \hfill\cr\hfill
    +\int_{(-q^{1/\alpha},p^{1/\alpha})}(e^{ity}-1-ity)\d\nu(y)
  \Bigr)
  \qquad\equa{CFLOne}\cr}
$$
(Note that the expression for $\nu$ as well as for this characteristic 
function given in the introduction (p.\ 3) are wrong and must be changed; 
to be honest, I had not done those calculations. I surely should have.).

\bigskip

Now, let us settle the discrepency with Resnick.

For this, consider, at least formally, the partial sums obtained for $k=\Id$,
and write $R(1)$ for the limiting process that Resnick considers. For us, we 
have that 
$$
  {\SS_n(1)-s_n(1)\over F_*^\leftarrow(1-1/n)}\to L(1)
$$
with
$$
  s_n(1)=n\E X\One\{\, F^\leftarrow(1/n)< X\leq F^\leftarrow(1-1/n)\,\}
$$
while Resnick uses a different centering,
$$
  r_n(1)=n\E X \One\{\, |X|\leq F_*^\leftarrow(1-1/n)\,\}
$$
to obtain
$$
  {\SS_n(1)-r_n(1)\over F_*^\leftarrow(1-1/n)} \to R(1) \, .
$$
The characteristic function of $R(1)$ is given (following Resnick; I did not
check his calculation, but what follows shows that either we are all right or
all wrong!)
$$\displaylines{
  \E \exp\bigl(itR(1)\bigr)
  =\exp\Bigl(\int_{(-\infty,-1)\cup (1,\infty)} 
   (e^{ity}-1)\d\nu(y) 
  \hfill\cr\hfill
    +\int_{(-1,1)}(e^{ity}-1-ity)\d\nu(y)
  \Bigr)
  \cr}
$$

Given how $R(1)$ and $L(1)$ are obtained as limit, assuming for the time
being that
$$
  \delta=\limn {r_n(1)-s_n(1)\over F_*^\leftarrow(1-1/n)}
$$
exists, we should have
$$
  \E\exp\bigr(itR(1)+it\delta\bigr)=\E\exp\bigl(itL(1)\bigr) \, .
$$
So, let us check this identity.

First,
$$\displaylines{\qquad
  {r_n(1)-s_n(1)\over n} 
  =\E X\Bigl(\One\{\, -F_*^\leftarrow(1-1/n)\leq X<F^\leftarrow(1-1/n)\,\}
  \quad\cr\hfill
  {}+ \One\{\, F^\leftarrow(1-1/n)<X\leq F_*^\leftarrow(1-1/n)\,\}\Bigr)
  \, .\cr}
$$

We have
$$\displaylines{
  \E X\One\{\, F^\leftarrow(1-1/n)<X\leq F_*^\leftarrow(1-1/n)\,\}
  \hfill\cr\noalign{\vskip 3pt}\hfill
  \eqalign{
  {}={}   &\int_{F^\leftarrow(1-1/n)}^{F_*^\leftarrow(1-1/n)}x\d F(x)\cr
  {}\sim{}& {\alpha\over 1-\alpha}\Bigl( (\Id\oF)\circ F_*^\leftarrow(1-1/n)
             -(\Id\oF)\circ F^\leftarrow(1-1/n)\Bigr) \cr
  }
  \cr}
$$
Since $\oF\sim p F_*^\leftarrow$, we have $\oF\circ F_*^\leftarrow(1-1/n)
\sim p/n$ and $F^\leftarrow(1-1/n)\sim p^{1/\alpha}F_*^\leftarrow(1-1/n)$.
So,
$$
  (\Id\oF)\circ F_*^\leftarrow(1-1/n)\sim F_*^\leftarrow(1-1/n){p\over n}
$$
and
$$
  (\Id\oF)\circ F^\leftarrow(1-1/n)\sim p^{1/\alpha} 
  {F_*^\leftarrow(1-1/n)\over n} \, .
$$
It follows that
$$\displaylines{\qquad
  \E X\One\{\, F^\leftarrow(1-1/n)<X\leq F_*^\leftarrow(1-1/n)\,\}
  \hfill\cr\noalign{\vskip 3pt}\hfill
  \eqalign{
  {}\sim{}& {\alpha\over 1-\alpha}
            {F_*^\leftarrow(1-1/n)\over n} (p-p^{1/\alpha}) \cr
  {}\sim{}& {F_*^\leftarrow(1-1/n)\over n} \alpha p
            {1-p^{(1/\alpha)-1}\over 1-\alpha}\cr
  }
  \qquad\cr}
$$
Thus,
$$
  {r_n(1)-s_n(1)\over F_*^\leftarrow(1-1/n)}
  \sim 
  \alpha \Bigl( p{1-p^{(1/\alpha)-1}\over 1-\alpha}
   -q{1-q^{(1/\alpha)-1}\over 1-\alpha}\Bigr)
  =\delta
  \, .
$$

Now, 
$$\displaylines{
  {\E\exp\bigr(itR(1)\bigr)\over \E\exp\bigr(itL(1)\bigr)}
    =\exp\Bigl(\int_{(-q^{1/\alpha},-1)\cup (1,p^{1/\alpha})} 
   (e^{ity}-1)\d\nu(y) 
  \hfill\cr\hfill
    +\int_{(-1,-q^{1/\alpha})\cup (p^{1/\alpha},1)}(e^{ity}-1-ity)\d\nu(y)
  \Bigr)
  \cr}
$$
Note that the parts involving $e^{ity}-1$ cancel. So this ratio of 
characteristic functions is simply
$$
\exp\Bigl(-\int_{(-1,-q^{1/\alpha})\cup (p^{1/\alpha},1)}ity\d\nu(y)\Bigr) \, .
$$
Its logarithm is
$$\displaylines{\quad
  -\alpha q\int_{-1}^{-q^{1/\alpha}} it y|y|^{-\alpha-1}\d y
  -\alpha p\int_{p^{1/\alpha}}^1 it y^{-\alpha}\d y 
  \hfill\cr\hfill
  \eqalign{
  {}={}& \alpha q\int_{q^{1/\alpha}}^1 it z^{-\alpha}\d z
         -\alpha p\int_{p^{1/\alpha}}^1 it y^{-\alpha}\d y \cr
  {}={}& {\alpha\over \alpha-1} itq (1-q^{1/\alpha-1}) 
        -{\alpha\over 1-\alpha} itp (1-p^{1/\alpha-1}) \cr
  {}={}&-it\delta \, . \cr
  }
  \cr}
$$
So, indeed,
$$
  \E\exp\bigl( itR(1)\bigr) = e^{it\delta}\E\exp\bigl(itL(1)\bigr) \, .
$$

\bigskip

\noindent{\bf Conclusion.} I suggest that in the paper we put
$$
  {\d\nu\over\d\lambda}(x)
  =p\alpha x^{-\alpha-1}\One_{(0,\infty)}(x)x^{-\alpha-1}
  + q\alpha |x|^{-\alpha-1}\One_{(-\infty,0)}(x) \, .
$$
and
$$\displaylines{
  \E\exp\bigl(itL(1)\bigr)
  =\exp\Bigl(\int_{(-\infty,-q^{1/\alpha})\cup (p^{1/\alpha},^\infty)} 
   (e^{ity}-1)\d\nu(y) 
  \hfill\cr\hfill
    +\int_{(-q^{1/\alpha},p^{1/\alpha})}(e^{ity}-1-ity)\d\nu(y)
  \Bigr) \, .
  \qquad\cr}
$$
Then, add a remark that if instead of using the centering 
$$
  \E X\One_{(F^\leftarrow(1/n),F^\leftarrow(1-1/n))}(X)
$$ 
we use
$$
  \E X\One_{(-F_*^\leftarrow(1-1/n),F_*^\leftarrow(1-1/n))}(X)
$$ 
the process $L$ should be replaced by a stable L\'evy process $\tilde L$ with
the nicer characteristic function
$$\displaylines{
  \E\exp\bigl(it\tilde L(1)\bigr)
  =\exp\Bigl(\int_{|y|>1} (e^{ity}-1)\d\nu(y) 
  \hfill\cr\hfill
    +\int_{|y|<1}(e^{ity}-1-ity)\d\nu(y)
  \Bigr) \, .
  \cr}
$$
However, the normalization that we have is the one usually used in empirical
process and the most natural one given the technique that we use.

And perhaps add the difference of centering amount to a shift of the limiting
process given by
$$
  \tilde L= L+ 
  {\alpha\over 1-\alpha} \Bigl( p(1-p^{(1/\alpha)-1})
   -q(1-q^{(1/\alpha)-1})\Bigr) \, .
$$


\vfill\eject
\def\DATE{{March 25, 2010}}
\pageno=1
\snumber=0
\equanumber=0
\def\prevs{B.\the\snumber }
\def\preveq{({\rm B}.\the\equanumber)}

\noindent Dear Bill,

\bigskip

\noindent
here are some calculations related to the process $L$ which show that the 
assertion that $L$ has a scaling property (made in the introduction) is not 
quite correct (but can be easily fixed).

\bigskip

We start with $L^+=\limeps L_\epsilon^+$ which we have after the proof of 
Lemma \DeltaConverge\ and before the statement of Lemma \lemmaCentering, and
recall that
$$
  L_\epsilon^+(s)=\int_{[0,s]\times (1,\infty)} x\d N(v,x)
  +\int_{[0,s]\times(\epsilon^{1/\alpha},1)} x\d (N-\E N)(v,x)
$$
where $N$ is a Poisson process with
$$
  \E N\bigl([0,t]\times [x\infty)\bigr)= t x^{-\alpha} \, .
$$

It is clear that $L_\epsilon^+$ has independent increments, because of the 
independence of the Poisson process on nonoverlaping intervals.

\bigskip

\noindent
{\bf Proof of  $L^+(t)-L^+(s)\eqd L^+(t-s)$.}
\medskip

Let $0<s<t<1$. 
$$\eqalign{
  & \hskip -15pt L_\epsilon^+(t)-L_\epsilon^+(s) \cr
  &{}= \int_{[s,t]\times(1,\infty)} x\d N(v,x) 
    +\int_{[s,t]\times [\epsilon^{1/\alpha},1]} x \d (N-\E N)(v,x) \cr
  &{}\eqd L_\epsilon^+(t-s) \, . \cr
  }
$$
because the Poisson process $N$ is homogenous on its first component and 
its mean measure is a product one.\hfill\qed

\bigskip

\noindent
{\bf Calculation of $\E e^{i\lambda L^+(t)}$.} (cf.\ Lemma \LPlusLMeasure)

\medskip

Set 
$$
  f(v,x)= x\One_{[0,t]\times[\epsilon^{1/\alpha},\infty)}(v,x)
$$
and
$$
  g(v,x)=x\One_{[0,t]\times [\epsilon^{1/\alpha},1]}(v,x) \, .
$$
We have
$$
  L_\epsilon^+(t) = \int f\d N -\int g\d \E N \, .
$$
Therefore,
\hfuzz=7pt
$$\eqalign{
  \E e^{i\lambda L_\epsilon^+(t)}
  &{}=\exp\Bigl( \int e^{i\lambda f}-1\d \E N - i\lambda\int g\d \E N\Bigr) \cr
  &{}=\exp\Bigl( \int_{[0,t]\times[\epsilon^{1/\alpha},\infty)} 
    (e^{i\lambda x}-1)\d v \alpha x^{-\alpha-1} \d x \cr
  &\hskip 104pt  -i\lambda\int_{[0,t]\times [\epsilon^{1/\alpha},1]} x \d v \alpha 
    x^{-\alpha-1} \d x \Bigr) \cr
  \noalign{\vskip 3pt}
  &{}=\exp\Bigl( t\int_{(\epsilon^{1/\alpha},\infty)}(e^{i\lambda x}-1)
    \d\nu^+(x)
    -i\lambda t\int_{[\epsilon^{1/\alpha},1]} x\d\nu^+(x) \Bigr) \cr
  &{}=\exp\Bigl( t\int_{[\epsilon^{1/\alpha},1]} e^{i\lambda x}-1-i\lambda x
    \One_{[\epsilon^{1/\alpha},\infty)}(x) \d \nu^+(x)\Bigr) \, .\cr
  }
$$
\hfuzz=0pt
Therefore,
$$
  \E e^{i\lambda L^+(t)}
  =\exp\Bigl( t\int_{[0,1]} e^{i\lambda x}-1-i\lambda t
    \One_{[0,\infty)}(x) x \d \nu^+(x)\Bigr) \, .
  \eqno{\equa{AppA}}
$$

As it should, we see that
$$
  \E e^{i\lambda L^+(t)} = \bigl( \E e^{i\lambda L^+(1)}\bigr)^t \, .
$$

\bigskip

\noindent{\bf Selfsimilarity of $L^+$}

Consider the integral in \AppA\ and make the change of variable 
$x=t^{1/\alpha}y$. We have
$$
  \d\nu^+(x) 
  = {\d\nu^+\over\d\lambda}(x)\d x
  = \alpha  (t^{1/\alpha}y)^{-\alpha-1} t^{1/\alpha} \d y
  = t^{-1}\d\nu^+(y) \, .
$$
Therefore,
$$\displaylines{\quad
  t \int_0^\infty \bigl(e^{i\lambda x}-1-i\lambda x\One_{[0,1]}(x)\bigr)
  \d\nu^+(x) \hfill\cr\qquad
  {}=\int_0^\infty \bigl(e^{i\lambda t^{1/\alpha}y}-1-
        i\lambda t^{1/\alpha}y\One_{[0,1]}(t^{1/\alpha}y)\bigr)
        \d\nu^+(y) 
  \hfill\cr\qquad
  {}=\int_0^\infty \bigl(e^{i\lambda t^{1/\alpha}y}-1-
        i\lambda t^{1/\alpha}y\One_{[0,1]}(y)\bigr)
        \d\nu^+(y) 
  \hfill\cr\hfill
  {} -\int_0^\infty  i\lambda t^{1/\alpha}y
         \bigl( \One_{[0,1]}(t^{1/\alpha}y)-\One_{[0,1]}(y)\bigr)
         \d\nu^+(y) \, .
  \qquad\equa{DisplayA}\cr}
$$
But
$$\eqalign{
  \One_{[0,1]}(t^{1/\alpha}y)-\One_{[0,1]}(y)
  &{}=\One_{[0,t^{-1/\alpha}]}(y)-\One_{[0,1]}(y) \cr
  &{}=\One_{(1,t^{-1/\alpha}]}(y) \, .\cr
  }
$$
Therefore, the last line of \DisplayA\ is
$$
 -i\lambda t^{1/\alpha} \int_1^{t^{-1/\alpha}} \alpha y^{-\alpha} \d y
  =
  \cases{ i\lambda {\alpha\over\alpha-1} (t^{1/\alpha}-t) 
           & if $\alpha\not=1$,\cr
         - i\lambda t\log t & if $\alpha=1$.\cr}
$$
So we have
$$
  \E e^{i\lambda L^+(t)} 
  = 
  \cases{\E e^{i\lambda t^{1/\alpha}L^+(1) -i\lambda {\alpha\over\alpha-1}
          (t^{1/\alpha}-t)} 
              & if $\alpha\not=1$, \cr
          \E e^{i\lambda tL^+(1)-i\lambda t\log t} & if $\alpha=1$.\cr}
$$
It follows that
$$\cases{
  L^+(t)+{\alpha\over\alpha-1} (t^{1/\alpha}-t)\eqd t^{\alpha}L^+(1)
    & if $\alpha\not=1$, \cr
  L^+(t)+t\log t\eqd t L^+(1) & if $\alpha=1$. \cr}
  \eqno{\equa{LPlusSelfSimilar}}
$$
So the process
$$
  \cases{ L^++{\alpha\over\alpha-1}(\Id^{1/\alpha}-\Id) & if $\alpha\not=1$,\cr
          L^++\Id\log & if $\alpha=1$\cr}
$$
is self-similar.

\bigskip

\noindent
{\bf Relation between $L^+$ and $L^-$.} 

\medskip

We define $L^-$ has a L\'evy process with
$$
  \E e^{i\lambda L^-(t)}
  = \exp\Bigl( t\int e^{i\lambda y}-1-i\lambda y\One_{[-1,0]}(y)\d\nu^-(y)
  \Bigr) \, .
$$

The first thing to do is to check that this indeed defines a process. For this
tarting with \AppA, we have
$$\eqalign{
  \E e^{-i\lambda L^+(t)}
  &{}=\exp\Bigl( t\int_0^\infty \bigl( e^{-i\lambda x}-1+i\lambda x
     \One_{[0,1]}(x)\bigr)\alpha x^{-\alpha-1}\d x\Bigr) \cr
  &{}=\exp\Bigl( t\int_{-\infty}^0 \bigl( e^{i\lambda y}-1-i\lambda y
     \One_{[-1,0]}(y)\bigr)\alpha (-y)^{-\alpha-1}\d y\Bigr) \cr
  &{}=\exp\Bigl( t\int_{-\infty}^0 \bigl( e^{i\lambda y}-1-i\lambda y
     \One_{[-1,0]}(y)\bigr)\d\nu^- (y)\Bigr) \cr
  }
$$
Therefore, $L^-$ is indeed well defined and
$$
  -L^+\eqd L^- \, .
  \eqno{\equa{LinkLPlusLMinus}}
$$ 
In particular, \LPlusSelfSimilar\ yields
$$\cases{
  L^-(t)-{\alpha\over\alpha-1}(t^{1/\alpha}-t)\eqd t^{1/\alpha} L^-(1) 
    & if $\alpha\not=1$, \cr
  L^-(t)-t\log t \eqd t L^-(1) & if $\alpha=1$.\cr
  }
  \eqno{\equa{LMinusSelfSimilar}}
$$

\bigskip

\noindent
{\bf Calculation of $\E L^+(1)$ when $\alpha>1$.}

Starting with \AppA, we have, as $\lambda$ tends to $0$,
$$
  \E e^{i\lambda L^+(t)}
  = \exp\Bigl( t\int_0^\infty i\lambda x\One_{[1,\infty)}(x)\d\nu^+(x)
  +o(\lambda)\Bigr) \, .
$$
But
$$
  \int_1^\infty x\d\nu^+(x)
  = \alpha\int_1^\infty x^{-\alpha}\d x
  ={\alpha\over\alpha-1} \, .
$$
Therefore, as $\lambda$ tends to $0$,
$$\eqalign{
  \E e^{i\lambda L^+(t)}
  &{}=\exp\Bigl( ti\lambda {\alpha\over\alpha-1} +o(\lambda)\Bigr) \cr
  &{}=1+i\lambda t{\alpha\over\alpha-1} +o(\lambda) \cr
  }
$$
This implies
$$
  \E L^+(t)= t{\alpha\over\alpha-1} \, ,
$$
and, similarly,
$$
  \E L^-(t)=-t{\alpha\over\alpha-1} \, .
$$
Since (we will check this relation later in these notes)
$$
  L=p^{1/\alpha}L^++q^{1/\alpha}L^-
  \eqno{\equa{LRep}}
$$
we have
$$
  \E L(t) = (p^{1/\alpha}-q^{1/\alpha}) t {\alpha\over\alpha-1} \, . 
$$

\bigskip

\noindent
{\bf Selfsimilarity of $L$}

\medskip

From \LPlusSelfSimilar, \LMinusSelfSimilar\ and \LRep, we deduce
$$\eqalign{
  L(t)
  &{}\eqd p^{1/\alpha}\Bigl( t^{1/\alpha}L^+(1) -{\alpha\over\alpha-1} 
    (t^{1/\alpha}-t)\Bigr) \cr
  &\hskip 60pt{}+ q^{1/\alpha} \Bigl( t^{1/\alpha}L^-(1) +{\alpha\over\alpha-1} 
    (t^{1/\alpha}-t)\Bigr) \cr
  &{}=t^{1/\alpha} L(1) -{\alpha\over\alpha-1} (p^{1/\alpha}-q^{1/\alpha})
    (t^{1/\alpha}-t) \cr
  &{}=t^{1/\alpha} L(1) -{\alpha\over\alpha-1} (p^{1/\alpha}-q^{1/\alpha})
    t^{1/\alpha} + \E L(t) \, . \cr
  }
$$
It follows that, when $\alpha>1$,
$$
  L(t)-\E L(t)
  \eqd
  t^{1/\alpha} \bigl( L(1)-\E L(1)\bigr)
$$
and, more generally,
$$\left\{
  \eqalign{
  &{}L(t)-{\ds\alpha\over\ds\alpha-1}(p^{1/\alpha}-q^{1/\alpha}) t 
  \eqd
  t^{1/\alpha} 
  \Bigl( L(1)-{\ds\alpha\over\ds\alpha-1}(p^{1/\alpha}-q^{1/\alpha})\Bigr)\cr
  & \hskip 240pt\hbox{if $\alpha\not=1$} \cr
  & L(t)+(p-q)t\log t\eqd t L(1) \qquad\hbox{if $\alpha=1$.} \cr}
  \right.
$$

\bigskip

\noindent{\bf Making sure we wrote the limiting processes properly.}

\medskip

Write $\TT_n^+[b_n,F]$ for $\TT_n^+$ calculated with the sequence $(b_n)$ and
when the $X_i$ have distribution $F$, and similarly for $\TT_n^-[a_n,F]$
and the centerings $\ttT_n^+[b_n,F]$ and $\ttT_n^-[a_n,F]$.

We proved that $\TT_n^+-\ttT_n^+\to k\star\dot L^+$.

We have
$$
  \M_{-1}F(x)
  =\Prob\{\, -X\leq x\,\}
  =\Prob\{\, X\geq -x\,\}
  =1-F(-x-) \, . 
$$
Thus, if $F(x_n)=1/n$, then $\M_{-1}F(-x_n-)=1-1/n$ and $x_n\sim 
(\M_{-1}F)^\leftarrow(1-1/n)$. So we have
$$
  F^\leftarrow(1/n)\sim -(\M_{-1}F)^\leftarrow(1-1/n)
$$
as $n$ tends to infinity. This implies
$$\displaylines{\quad
  (\TT_n^--\ttT_n^-)[a_n,F](t)
  \hfill\cr\qquad
    \sim {-1\over (\M_{-1}F)^\leftarrow(1-1/n)} \sum_{1\leq i\leq n}
      k_n\Bigl( t-{i\over n}\Bigr)
      \bigl( X_i\One_{(-\infty,a_n)}(X_i)-\mu_n^-\bigr) 
  \hfill\cr\qquad
    {}= {1\over (\M_{-1}F)^\leftarrow(1-1/n)} \sum_{1\leq i\leq n}
      k_n\Bigl( t-{i\over n}\Bigr)
  \hfill\cr\hskip 120pt
      {}\times\bigl( -X_i\One_{(-a_n,\infty)}(-X_i)-\mu_n^+[\M_{-1}F,a_n]\bigr) 
  \hfill\cr\qquad
    {}=(\TT_n^+-\ttT_n^+)[-a_n,\M_{-1}F](t) \, . \hfill
    \equa{TMinusTPlus}\cr
  }
$$
Consider the L\'evy measures $\nu^+[\M_{-1}F]$ and $\nu^-[F]$. Under the tail
balance condition, $\nu^+[\M_{-1}F]=\nu^+[F]$. Thus \TMinusTPlus\ yields
Therefore, since $\nu^+[\M_{-1}F]=\nu^-[F]=\nu^+[F]$ under the tail
balance condition,
$$
  (\TT_n-\ttT_n^-)\to k\star \dot L^+
$$
Given \LinkLPlusLMinus, we have
$$
  (\TT_n-\ttT_n^-)\to -k\star\dot L^-
$$
(this is not what we wrote on p.38, l.7: there, we should have 
$(\TT_n^--\ttT_n^-,\TT_n^+-\ttT_n^+)$ converges 
to $(-k\star\dot L^-,k\star\dot L^+)$ where $L^-$ and $L^+$ are 
indeed independent).

Now, given \SnVn, we have
$$
  {\SS_n-s_n\over F_*\leftarrow(1-1/n)}
  \to p^{1/\alpha}L^++q^{1/\alpha}L^- \, .
$$
So
$$
  L=p^{1/\alpha}L^+ + q^{1/\alpha}L^- \, .
$$
(so there is a sign wrong on p.38, l.-13)

\bigskip

\noindent
{\bf Conclusion.} I propose the following changes:

p.3: remove l.-4,-5 and add the correct selfsimilarity property after the 
first sentence of p.4 and do a little fixing for p.4, l.2--5.

On p.38, spell out the arguement showing that $\TT_n^--\ttT_n^-$ converges
to $-L^-$ and fix p.38, l.-14 (the first $-$ on that line should be a $+$).

\bye


\Lemma{\label{LTwoCv}
  Let $(k_n)_{n\geq 1}$ be a sequence of functions linearly interpolated 
  between points $i/n$, $0\leq i\leq n$. Assume that

  \medskip

  \noindent (i) $\limn k_n(\lfloor nt\rfloor/n)=k(t)$ for almost every $t$;

  \medskip

  \noindent (ii) $k$ is in square integrable;

  \medskip

  \noindent (iii) $\limsupn n^{-1}\sum_{0\leq i\leq n} k_n^2(i/n)$ is finite;

  \medskip

  \noindent (iv) $\limn n^{-1}\sum_{1\leq i\leq n} 
  \Bigl( k_n\Bigl({\ds i\over\ds n}\Bigr)-k_n\Bigl({\ds i-1\over\ds n}\Bigr)
  \Bigr)^2=0$.

  \medskip

  Then $k_n$ converges to $k$ in $\LTwo$ and for any $s$ and $t$,
  $$
    \limn n^{-1}\sum_{1\leq i\leq n} k_n(s-i/n)k_n(t-i/n)=\int_0^1k(s-u)k(t-u)
    \d u \, .
  $$
}

\Proof We write $\calF(x)$ for the fractional part of a positive real 
number $x$, that is for $x-\lfloor x\rfloor$. Define
$$
  \Delta_n(u)
  =k_n\Bigl( {\lfloor nu\rfloor+1\over n}\Bigr) 
  -  k_n\Bigl( {\lfloor nu\rfloor\over n}\Bigr) \, ,
$$
so that, since $k_n$ is linearly interpolate,
$$
  k_n(u)=k_n\Bigl({\lfloor nu\rfloor\over n}\Bigr) + \calF(nu)\Delta_n(u) \, .
$$
Note that the function $\Delta_n$ is constant on intervals $[(i-1)/n,i/n)$. 

We have
$$\displaylines{\quad
  \int_{(i-1)/n}^{i/n} (k_n-k)^2(u) \d u
  \hfill\cr\hfill
  \leq \int_{(i-1)/n}^{i/n} \Bigl( k_n\Bigl({\lfloor nu\rfloor\over n}\Bigr)
  -k(u)\Bigr)^2 \d u
  + 2\int_{(i-1)/n}^{i/n} \Delta_n^2(u) \d u \, .
  \quad
  \equa{LTwoCvA}\cr}
$$
Assumption (i) ensures that $\bigl(k_n(\lfloor nu\rfloor/n) -k(u)\bigr)^2$ 
converges to $0$ almost everywhere. Moreover, this function is at most
$2k_n^2(\lfloor nu\rfloor/n)+2k^2(u)$, which are uniformly integrable thanks to
(iii) and (ii). It follows that
$$
  \limn \int_0^1\Bigl( k_n\Bigl({\lfloor nu\rfloor\over n}\Bigr)-k(u)\Bigr)^2
  \d u = 0 \, .
$$
Furthermore, assumption (iv) ensures that
$$
  \limn \int_0^1\Delta_n^2(u)\d u= 0 \, .
$$
Combined with \LTwoCvA, this yields the convergence of $k_n$ to $k$ in $\LTwo$.

To prove the second assertion of the lemma, we first show that for any fixed
$s$,
$$
  \limsupn n^{-1}\sum_{1\leq i\leq n} k_n(s-i/n)^2 < \infty \, .
  \eqno{\equa{LTwoCvB}}
$$
Indeed, we have
$$
  k_n^2(s-i/n)\leq 2k_n^2\Bigl({\lfloor ns\rfloor -i\over n}\Bigr)^2 
  +\Delta_n^2\Bigl(s-{i\over n}\Bigr) \, ,
$$
so that \LTwoCvB\ follows from (iii) and (iv). 

Now, let $u$ between $(i-1)/n$ and $i/n$. Then
$$
  \lfloor ns\rfloor - i
  = \lfloor ns-i\rfloor
  \leq \lfloor n(s-u)\rfloor
  \leq \lfloor ns-i+1\rfloor
  = \lfloor ns\rfloor - (i-1) \, .
$$
It follows that
$$
  |k_n(s-i/n)-k_n(s-u)| \leq 3\Delta_n(s-i/n) \, .
$$
Consequently, the absolute value of
$$\displaylines{\quad
  k_n(s-i/n)k_n(t-i/n)-k_n(s-u)k_n(t-u)
  \hfill\cr\noalign{\vskip 3pt}\hfill
  \eqalign{
    {}={}& k_n(s-i/n)\bigl(k_n(t-i/n)-k_n(t-u)\bigr) \cr
         & {}+ k_n(t-i/n)\bigl( k_n(s-i/n)-k_n(s-u)\bigr)\cr
         & {}+\bigl( k_n(t-u)-k_n(t-i/n)\bigr) \bigl(k_n(s-i/n)-k_n(s-u)\bigr)
    \cr}
  \quad\cr}
$$
is at most
$$\displaylines{\qquad
  3|l_n(s-i/n)\Delta_n(t-i/n)| + 3|k_n(t-i/n)\Delta_n(s-i/n)|
  \hfill\cr\hfill
  {}+9|\Delta_n(t-i/n)\Delta_n(s-i/n)| \, .
  \qquad\cr}
$$
This implies, after splitting the following integral over $[\,0,1\,]$ as a sum
of integrals over $[(i-1)/n,i/n)$ and using Cauchy-Schwartz' inequality, that
$$\displaylines{
  \Bigl| n^{-1}\sum_{1\leq i\leq n} k_n(s-i/n)k_n(t-i/n)
  -\int_0^1k_n(s-u)k_n(t-u)\d u\Bigr|^2
  \hfill\equa{LTwoCvC}\cr\hfill
  \eqalign{
    \leq{}& 3n^{-1} \sum_{1\leq i\leq n} k_n^2(s-i/n) 
            n^{-1}\Delta_n^2(t-i/n) \cr
          &{}+3 n^{-1} \sum_{1\leq i\leq n} k_n^2(t-i/n) 
            n^{-1}\sum_{1\leq i\leq n} \Delta_n^2(s-i/n) \cr
          &{}+9n^{-1}\sum_{1\leq i\leq n} \Delta_n^2(t-i/n) 
            9n^{-1}\sum_{1\leq i\leq n} \Delta_n^2(s-i/n) \, . \cr
    }
  \qquad\cr
  }
$$
It then follows from (iv) and \LTwoCvB\ that \LTwoCvC\ tends to $0$ as $n$
tends to infinity. Since $(k_n)$ converges to $k$ in $\LTwo$, we also have
$$
  \limn \int_0^1 k_n(s-u)k_n(t-u)\d u
  = \int k(s-u)k(t-u) \d u
$$
and this proves the last assertion of the lemma.\hfill\qed

\bye

\bigskip

We need to extend Lemma \ModulusUnKnB\ by removing the restriction that $ns$
and $nt$ are integers. This requires the additional assumption \HypKnModulus\
on the modulus
$\omega_{k_n}$. The following result deals only with even powers and will be
sufficient for our purpose.

\Lemma{\label{ModulusUn}
  Assume that \HypKnModulus\ is satisfied. For any integer $r$, there 
  exists a constant $c_r$
  such that for any $n$ large enough and any $s$ and $t$ in $[\,0,1\,]$,
  $$
    E\bigl( \MM_n(t)-\MM_n(s)\bigr)^{2r}\leq c_{2r} |t-s|^{(\theta\wedge 1) 2r}\, .
  $$
}

\Proof The result being obvious if $r$ vanishes, we assume that $r$ is at least
$1$.

\noindent{\it Step 1.} We consider a real number $s$ in $[\,0,1\,]$ and define
the positive integer $p=\lceil ns\rceil$. Since $k_n$ and hence $\MM_n$ are 
linearly interpolated,
$$
  \MM_n(p/n)-\MM_n(s)
  =(p-ns)\Bigl( \MM_n\Bigr( {p\over n}\Bigr)-\MM_n\Bigl({p-1\over n}\Bigr)\Bigr) 
  \, .
$$
Therefore, Lemma \ModulusUnKnB\ implies
$$
  \E\bigl( \MM_n(p/n)-\MM_n(s)\bigr)^{2r}\leq (ns-p)^{2r} c \omega_{k_n}(1/n)^{2r} \, .
  \eqno{\equa{ModulusUnA}}
$$
Using assumption \HypKnModulus, this yields a bound in 
$c(ns-p)^{2r}n^{-2r\theta}$. Since $ns-p$ is less than $1$ and $r$ is at least
$1$, we have $(ns-p)^{2r}\leq (ns-p)^{2r(\theta\wedge 1)}$. Thus \ModulusUnA\
is at most 
$$
  c(ns-p)^{2r(\theta\wedge 1)} n^{-2r(\theta\wedge 1)}
  = c(s-p/n)^{2r(\theta\wedge 1)} \, .
$$

\noindent{\it Step 2.} If there exist two integers $p$ and $q$ with
$$
  {p-1\over n} <s\leq {p\over n} \leq {q\over n} \leq t < {q+1\over n} \, ,
$$
we have
$$\displaylines{\quad
  \bigl( (\MM_n(t)-\MM_n(s)\bigr) ^{2r}
  \leq 3^r \Bigl( \bigl( \MM_n(t)-\MM_n(q/n)\bigr) ^{2r}
  \hfill\cr\hfill 
  + \bigl( \MM_n(q/n)-\MM_n(p/n)\bigr)^{2r}
  + \bigl( \MM_n(p/n)-\MM_n(s)\bigr)^{2r}\Bigr) \, .
  \quad\cr}
$$
Taking expectation, using Lemma \ModulusUnKnB\ assumption \HypKnModulus\ and 
our first step, we obtain
$$\displaylines{\qquad
  \E\bigl( (\MM_n(t)-\MM_n(s)\bigr) ^{2r}
  \hfill\cr\hfill
  \eqalign{
  {}\leq{}& c\Bigl( \Bigl( t-{q\over n}\Bigr)^{(\theta\wedge 1)2r} 
    +\Bigl( {q-p\over n}\Bigr)^{(\theta\wedge 1) 2r} 
    +\Bigl({p\over n}-s\Bigr)^{(\theta\wedge 1)2 r} \Bigr)\cr
  {}\leq{}& 3 c(t-s)^{(\theta\wedge 1)2 r} \, .\cr
  }\cr
  }
$$

\noindent{\it Step 3.} If $s$ and $t$ are in an interval 
$\bigl(p/n,(p+1)/n\bigr)$, Lemma \ModulusUnKnB\ implies
$$
  E\bigl( (\MM_n(t)-\MM_n(s)\bigr) ^{2r} \leq (t-s)^{2r} \omega_{k_n}^{2r}(1/n) \, ,
$$
and this proves the lemma.\hfill\qed

mpage -2 -c -o -M-100rl-80b-220t -t veraverbeke7.ps > toto.ps

\bye